\documentclass[10pt]{amsart}
\usepackage{amsmath}
\usepackage{xypic}
\usepackage{amssymb}
\usepackage{amsfonts}
\usepackage{amsthm}
\usepackage{enumerate}
\usepackage[all]{xy}
\usepackage[latin1]{inputenc}
\usepackage{graphicx}
\usepackage{latexsym}
\usepackage{accents}
\usepackage{mathdots}
\usepackage{mathrsfs}
\usepackage{color}
\input xy
\usepackage{tikz}
\usetikzlibrary{arrows}
\usetikzlibrary{decorations,arrows}
\usetikzlibrary{decorations.markings}
\usetikzlibrary{arrows.meta, bending}
\usetikzlibrary{positioning}

\makeatletter
    \newtheoremstyle{definition}
        {5pt}
        {3pt}
        {}
        {0pt}
        {\scshape}
        {.}
        {5pt}
        {\thmname{#1} \thmnumber{#2} \thmnote{[#3]}} 

\newtheoremstyle{theorems}
        {5pt}
        {3pt}
        {\itshape}
        {0pt}
        {\scshape}
        {.}
        {5pt}
        {\thmname{#1} \thmnumber{#2}\thmnote{[#3]}} 

\swapnumbers \theoremstyle{theorems}

\newtheorem{Theo}{Theorem}[subsection]
\newtheorem{Prop}[Theo]{Proposition}
\newtheorem{Cor}[Theo]{Corollary}
\newtheorem{Lemma}[Theo]{Lemma}

\newtheorem{Prop(BG)}[Theo]{Proposition (Bongartz-Gabriel)}
\newtheorem{Lemma(Asashiba)}[Theo]{Lemma(Asashiba)}
\newtheorem{Lemma(Gab)}[Theo]{Lemma(Gabriel)}
\newtheorem{Theo(Mil)}[Theo]{Theorem (Milicic)}

\theoremstyle{definition}
\newtheorem{Defn}[Theo]{Definition}

\newtheorem{Defn(Asashiba)}[Theo]{Definition (Asashiba)}

\theoremstyle{theorems}

\newcommand{\Z}{\mathbb{Z}}

\newcommand{\cA}{\mathcal{A}}

\def\La{\hbox{$\it\Lambda$}}

\newcommand{\mit}{\mathit}

\newcommand{\Hom}{{\rm Hom}}
\newcommand{\GHom}{{\rm GHom}}
\newcommand{\Ext}{{\rm Ext}}

\newcommand{\GMod}{{\rm GMod}}
\newcommand{\mmod}{{\rm mod}}
\newcommand{\soc}{{\rm soc \hspace{.4pt}}}
\newcommand{\rad}{{\rm rad \hspace{.4pt}}}

\def\Mod{\hbox{{\rm Mod}{\hskip 0.3pt}}}

\def\proj{\hbox{{\rm gproj}{\hskip 0.5pt}}}

\def\GrLa{{\rm GMod}\La}

\def\mf{\mathfrak}

\newcommand{\dt}{{\accentset{\hspace{2pt}\mbox{\large\bfseries .}}{}}}

\newcommand{\cdt}{\dt\hspace{2.5pt}}

\newcommand{\pdt}{{\hspace{.8pt}\dt\hspace{1.5pt}}}

\newcommand{\mk}{\mathfrak}

\newcommand{\mla}{\hspace{-.5pt} \langle \hspace{-0.4pt}}
\newcommand{\mra}{ \hspace{-.8pt} \rangle}

\newcommand{\nla}{ \hspace{-1pt} \langle \hspace{-0.4pt}}
\newcommand{\nra}{ \hspace{-.5pt} \rangle}

\newcommand{\sla}{\hspace{-.6pt} \langle \hspace{-.2pt}}
\newcommand{\sra}{\hspace{-.3pt}\rangle}

\newcommand{\tla}{\hspace{-.4pt} \langle \hspace{.4pt}}
\newcommand{\tra}{\hspace{-.0pt}\rangle}

\newcommand{\nid}{\hspace{-1.5pt}\mid\hspace{-1.5pt}}

\newcommand{\m}{\hspace{-.5pt}}

\begin{document}

\title[Graded algebras]{\sc Representation theory of graded algebras given by locally finite quivers
}

\author[Z. Lin]{Zetao Lin}

\author[S. Liu]{Shiping Liu 
}

\keywords{Algebras given by quivers with relations; graded algebras; graded modules; Nakayama functor; derived categories; Auslander-Reiten formula; almost split sequences; almost split triangles.}

\subjclass[2010]{16E35, 16G20, 16G70, 16W50.}



\address{Shiping Liu \\ D\'epartement de math\'ematiques, Universit\'e de Sherbrooke, Sherbrooke, Qu\'ebec, Canada.}
\email{shiping.liu@usherbrooke.ca }

\address{Zetao Lin \\ Department of Mathematical Sciences, Tsinghua University, 100084, Beijing, P. R. China.}
\email{zetao.lin@foxmail.com}

\begin{abstract}

This paper aims to study graded modules over a graded algebra $\La$ given by a locally finite quiver with homogeneous relations. By constructing a graded Nakayama functor, we discover a novel approach to establish Auslander-Reiten formulas,
from which we derive almost split sequences in the category of all graded $\La$-modules. 
In case $\La$ is locally left (respectively, right) bounded, the category of finitely presented graded modules and that of finitely copresented graded modules both have almost split sequences on the left (respectively, right). We shall also obtain existence theorems for almost split triangles in various derived categories of graded $\La$-modules. In case $\La$ is locally bounded, an indecomposable complex 
in the bounded derived category of finite dimensional graded modules is the starting (respectively, ending) term of an almost split triangle if and only if it has a finite graded projective reso\-lution (respectively, injective coresolution); and consequently, this bounded derived category has almost split triangles on the right (respectively, left) if and only if every graded simple module is of finite graded projective (respectively, injective) dimension. Finally, we specialize to the existence of almost split sequences and almost split triangles for graded representations of any locally finite quiver.  

\vspace{-15pt}

\end{abstract}

\maketitle

\section*{Introduction}


Graded algebras play an essential role in many domains such as commutative algebra, Lie theory, algebraic geometry and algebraic topology; see, for example, \cite{AuR0,
EKo, Har, Mac, Mat}. 
The representation theory of locally finite dimensional graded algebras have been studied thoroughly by numerous researchers; see, for example, \cite{AuR, GoG, Mar, RMV2}.
Motivated by the application of the covering technique; see \cite{BaL, BoG, 
Gre}, this paper aims to develop the representation theory of graded algebras given by locally finite quivers with homogeneous relations.


Almost split sequences in abelian categories, introduced by Auslander and Reiten; see \cite{AuR2}, and almost split triangles in triangulated categories, later developed by Happel; see \cite{Ha1}, provide a powerful tool for understanding these categories. 
In the classical setting, the existence of almost split sequences for graded modules was studied first by Gordon and Green; see \cite{GoG}, then by Auslander and Reiten; see \cite{AuR0}, and later by Martinez-Villa; see \cite{RMV2}. In the locally finite dimensional case, Martinez-Villa established an Auslander-Reiten formula 
by the classical approach of tensor product and adjunction isomorphism, and he obtained an existence theorem for almost split sequences ending with finitely presented graded modules in the category of locally finite dimensional graded modules. In this paper, we shall provide a novel approach to establish an Auslander-Reiten formula for finitely presented graded modules and a generalized Auslander-Reiten formula for finitely copresented graded modules, which enable us to obtain an existence theorem in the category of all graded modules for almost split sequences ending with finitely presented graded modules or starting with finitely copresented graded modules.


Determining which categories have almost split sequences or almost split triangles has long been a key research topic; see, for example, \cite{
GaR, LiM, LiN, SLPC, RMV2}. It is known that the category of finitely generated graded modules over a graded artin algebra or a graded order which is an isolated singularity has almost split sequences; see \cite{AuR, GoG}. Our results yield many interesting abelian categories of graded modules having almost split sequences on one or two sides, and derived categories of graded modules having almost split triangles on one or two sides. The content of the paper is outlined section by section as follows.

In Section 1, we shall lay down the foundation of the paper. In Section 2, we study some generalities about the category $\GrLa$ of all unitary graded left modules over a graded algebra $\La$ given by a locally finite quiver with homogeneous relations. 
Of fundamental importance, there exists a duality between the piecewise finite dimensional graded modules over $\La$ and those over its opposite; see (\ref{gr-duality}), and the categories ${\rm gproj}\La$ and ${\rm ginj}\La$ of finitely generated graded projective modules and of finitely cogenerated graded injective modules are Hom-finite and Krull-Schmidt; see (\ref{proj-KS}). As in the classical graded setting, $\GrLa$ has enough projective objects and enough injective objects; see (\ref{GM-proj}) and (\ref{GM-pi}). 


In Section 3, we study the existence of almost split sequences for graded $\La$-modules. We first construct a Nakayama functor from ${\rm gproj}\La$ to $\GrLa$; see (\ref{N-Functor}), which leads to an Auslander-Reiten formula for finitely presented graded $\La$-modules and a generalized Auslander-Reiten formula for finitely copresented graded modules; see (\ref{AR-formu}). These enable us to establish an existence theorem in ${\rm GMod}\La$ for almost split sequences ending with finitely presented graded modules or starting with finitely copresented graded modules; see (\ref{AR-sequence}). 
In case $\La$ is locally left (respectively, right) bounded, the categories of finitely presented graded modules and of finitely copresented graded modules have almost split sequences on the left (respectively, right); see (\ref{ass_pi}). And in case $\La$ is locally bounded, the category of finite dimensional graded modules has almost split sequences; see (\ref{rlb-ass}).

In Section 4, we study almost split triangles in various derived categories of graded $\La$-modules. Indeed, the graded Nakayama functor ensures the existence in the derived categories of almost split triangles ending with bounded complexes over ${\rm gproj}\La$ or starting with bounded complexes over ${\rm ginj}\La$; see (\ref{GAST-1}). In case $\La$ is locally left and right noetherian, an indecomposable complex of finitely generated (respectively, cogenerated) graded modules is the ending (respectively, starting) term of an almost split triangle if and only if it has a finite graded projective resolution (respectively, injective coresolution); see (\ref{GAST-lb}). And the bounded derived category of finite dimensional graded $\La$-modules has almost split triangles on the right (respectively, left) if and only if every graded simple $\La$-module is of finite graded projective (respectively, injective) dimension.


In Section 5, we study the existence of almost split sequences and almost split triangles for graded representations of a locally finite quiver $Q$. In case $Q$ is strongly locally finite, this has been done for ungraded representations; see \cite{BLP}. Our results say that the abelian category of finitely presented (respectively, copresented) graded representations has almost split sequences on the left (respectively, right) if and only if 
$Q$ has no infinite path with a starting (respectively, end) point; see (\ref{ass_rep}). And the bounded derived category of finitely presented graded representations has almost split triangles if and only if so does the bounded derived category of finitely copresented graded representations if and only if $Q$ has no infinite path; see (\ref{GAST-3_hered}).

\section{Preliminaries}

\noindent The objective of this section is to lay down the foundation of this paper. We shall fix some terminology and notation which will be used throughout this paper and collect some preliminary results. 

\vspace{-6pt}

\subsection{\sc Linear algebra.} Throughout this paper, let $k$ denote a commutative field. All tensor products will be over $k$. Given a set $\mathcal S$, the $k$-space spanned by $\mathcal S$ will be written as $k\mathcal S$. We shall write ${\rm Mod}\hspace{.4pt}k$ for the category of all $k$-spaces and ${\rm mod}\hspace{.4pt}k$ for the category of finite dimensional $k$-spaces. We shall make a frequent use of the exact functor $D=\Hom_k(-, k): {\rm Mod}\hspace{.4pt}k\to {\rm Mod}\hspace{.4pt}k$, which restricts to a duality $D: {\rm mod}\hspace{.4pt}k\to {\rm mod}\hspace{.4pt}k$. The following statement is well-known.

\begin{Lemma}\label{Tensor}

Given $U, V; M, N\in {\rm Mod}\hspace{.4pt}k$, there exists a $k$-linear map \vspace{-3pt}
$$\rho: \Hom_k(U, V) \otimes \Hom_k(M, N)\to \Hom_k(U\otimes M, V\otimes N): f\otimes g\mapsto \rho(f\otimes g) \vspace{-3pt} $$
such that $\rho(f\otimes g)(u\otimes m)=f(u) \otimes g(m)$ for $u\in U$ and $m\in M$, which is natural in all variables. Moreover, $\rho$ is an isomorphism in case $U, V\in {\rm mod}k$ or $M, N\in {\rm mod}k.$


\end{Lemma}


\noindent {\sc Remark.} In case $\rho$ is an isomorphism, we shall identify $f\otimes g$ with $\rho(f\otimes g)$.

\smallskip

As an immediate consequence of Lemma \ref{Tensor}, we have the following statement.

\begin{Cor}\label{Cor 1.2}

Let $U$ and $V$ be $k$-vector spaces. 

\begin{enumerate}[$(1)$]

\vspace{-2pt}

\item There exists a natural $k$-linear map $\sigma: DU \otimes V \to \Hom_k(U, V)$ in such a way that $\sigma(f\otimes v)(u)=f(u)v$, for all $f\in DU$, $v\in V$ and $u\in U$. Moreover, $\sigma$ is an isomorphism in case $U$ or $V$ is finite dimensional. 

\vspace{1pt}

\item  There exists a natural $k$-linear map $\theta: DV \otimes DU\to D(V\otimes U)$ such that $\theta(f\otimes g)(v\otimes u)= f(v)g(u)$, for all $f\in DU$, $g\in DV;$ $u\in U$ and $v\in V$. Moreover, $\theta$ is an isomorphism in case $U$ or $V$ is finite dimensional. 
\end{enumerate}
\end{Cor}

\vspace{-8pt}

\subsection{\sc Quivers} Let $Q=(Q_0, Q_1)$ be a quiver, where $Q_0$ is the set of vertices and $Q_1$ is the set of arrows between vertices. Given an arrow $\alpha: x\to y$ in $Q_1$, we call $x$ the {\it starting point} and $y$ the {\it end point} of $\alpha$; and write $s(\alpha)=x$ and $e(\alpha)=y$. For each vertex $x\in Q_0$, one associates a {\it trivial path} $\varepsilon_x$ with $s(\varepsilon_x)=e(\varepsilon_x)=x$. A {\it path} of positive length $n$ is a sequence $\rho=\alpha_n \cdots \alpha_1$, where $\alpha_i\in Q_1$ such that $s(\alpha_{i+1})=e(\alpha_i)$, for all $1\le i< n$. One says that $Q$ is {\it locally finite} provided, for any $x\in Q_0$, that the set $Q_1(x,-)$ of arrows $\alpha$ with $s(\alpha)=x$ and the set $Q_1(-,x)$ of arrows $\beta$ with $e(\beta)=x$ are both finite. Moreover, $Q$ is {\it strongly locally finite} provided, for  any $x,y\in Q_0$, that 
the set $Q(x,y)$ of paths from $x$ to $y$ is finite.



The {\it opposite} quiver $Q^{\rm o}$ of $Q$ is defined in such a way that $(Q^{\rm o})_0=Q_0$ and $(Q^{\rm o})_1=\{\alpha^{\rm o}: y\to x \mid \alpha: x\to y \in Q_1\}$. A non-trivial path $\rho=\alpha_n\cdots \alpha_1$ in $Q(x, y)$, where $\alpha_i\in Q_1$, corresponds to a non-trivial path $\rho^{\rm o}=\alpha_1^{\rm o} \cdots \alpha_n^{\rm o}$ in $Q^{\rm o}(y,x)$. For convenience, the trivial path in $Q^{\rm o}$ at a vertex $x$ will be identified with the trivial path in $Q$ at $x$.


\subsection{\sc Algebras given by quivers with relations} In this paper, a $k$-algebra does not necessarily have an identity. 
Let $Q=(Q_0, Q_1)$ be a locally finite quiver. We shall write $kQ$ for the path algebra of $Q$ over $k$ and $kQ^+$ for the two-sided ideal in $kQ$ generated by $Q_1$. A two-sided ideal in $kQ$ is called a {\it relation ideal} if it is contained in $(kQ^+)^2$. An element in $kQ$ is called {\it homogeneous} if it is a linear combination of paths of the same length, and a relation ideal in $kQ$ is called {\it homogeneous} if it is generated by some homogeneous elements.


Let $\La=kQ/R$, where $R$ is a relation ideal of $kQ$. We fix some notation for $\La$, which will be used throughout this paper. Write $\bar{\hspace{-1pt}\gamma}=\gamma+R\in \La$ for $\gamma\in kQ$, and $e_x = \bar{\varepsilon}_x$ for $x\in Q_0$. Then,
$\{e_x \nid x\in Q_0\}$ is a complete set of pairwise orthogonal idempotents in $\La$. Note that the opposite algebra of $\La$ is $\La^{\rm o}=kQ^{\rm o}/R^{\rm o},$ where $R^{\rm o}= \{\rho^{\rm o} \mid \rho\in R\}$. We shall write $\bar{\gamma\hspace{1.8pt}}^{\hspace{-1pt}\rm o}=\gamma^{\rm o}+R^{\rm o}$ for $\gamma\in kQ$, but $e_x=\varepsilon_x+R^{\rm o}$ for $x\in Q_0$. In this way, we have an algebra anti-isomorphism $\La\to \La^{\rm o}: \bar{\gamma\hspace{1.5pt}} \to \bar{\gamma\hspace{1.5pt}}^{\hspace{-1pt}\rm o}$.


Noetherian algebras play an important role in commutative algebra and algebraic geometry. In our context, noetherianness will be replaced by local noetherianness as follows: $\La$ is called {\it locally left noetherian} if the left $\La$-modules $\La e_x$ with $x\in Q_0$ are noetherian; {\it locally right noetherian} if the right $\La$-modules $e_x\La$ with $x\in Q_0$ are noetherian. 
As examples, recall that $\La$ is called {\it special multi-serial} provided, for any $\alpha\in Q_1$, that there exists at most one $\beta\in Q_1$ such that $\beta \alpha \notin R$ and at most one $\gamma\in Q_1$ such that $\alpha \gamma \notin R$. In this case, for any $x\in Q_0$, both $\sum_{\alpha\in Q_1(x,-)} \La \bar\alpha$ and $\sum_{\beta\in Q_1(-, x)} \bar\beta \La$ are finite sums of uniserial modules; see \cite{VHW, GrS}. Then, it follows that $\La e_x$ and $e_x\La$ are noetherian. Hence, $\La$ is locally left and right noetherian.

Finally, we shall say that $\La$ is {\it locally left bounded} if the $\La e_x$ with $x\in Q_0$ are finite dimensional, {\it locally right bounded} if the $e_x\La$ with $x\in Q_0$ are finite dimensional, and {\it locally bounded} if $\La$ is locally left and right bounded; compare \cite[(2.1)]{BoG}. Clearly, a locally left or right bounded algebra is locally left or right noetherian, respectively.


\subsection{\sc Additive categories} Throughout this paper, all categories are additive $k$-categories in which morphisms are composed from the right to the left. All functors between additive $k$-categories are additive. Let $\cA$ be an additive $k$-category. A full subcategory of $\cA$ is called {\it strictly full} if it is closed under isomorphisms. An object in $\cA$ is called {\it strongly indecomposable} if it has a local endomorphism algebra. One says that $\cA$ is {\it Hom-finite} if all morphisms spaces in $\cA$ are finite dimensional over $k$, and {\it Krull-Schmidt} if $\cA$ is nonzero such that every nonzero object is a finite direct sum of strongly indecomposable objects. In case $\cA$ is Hom-finite, it is well-known that $\cA$ is Krull-Schmidt if and only if all idempotents in $\cA$ split; see \cite[(1.1)]{SLPC}. 

A morphism $f: X\to Y$ in $\cA$ is called {\it left minimal} provided that every morphism $g:Y\to Y$ such that $g f=f$ is an automorphism, and  {\it right minimal} provided that every morphism $g:X\to X$ such that $fg=f$ is an automorphism. Applying Corollary 1.4 in \cite{HKMS} and its dual, we obtain the following well-known statement.

\begin{Prop}\label{minimal_mor}

Let $\cA$ be a Krull-Schmidt $k$-category. 

\begin{enumerate}[$(1)$]

\vspace{-2pt}

\item A nonzero morphism $f: X\to Y$ in $\cA$ is left minimal if and only if $pf\ne 0$, for any nonzero retraction $p: Y\to N.$

\item A nonzero morphism $f: X\to Y$ in $\cA$ is right minimal if and only if $fq\ne 0$, for any nonzero section $q: M\to X.$

\end{enumerate}
    
\end{Prop}

Finally, a morphism $f: X\to Y$ in $\cA$ is called {\it left almost split} if $f$ is not a section and every non-section morphism $g:X\to Z$ factors through $f$; and {\it minimal left almost split} if it is left minimal and left almost split. Dually, one defines {\it right almost split morphisms} and {\it minimal right almost split morphisms} in $\cA$; see \cite{AuR2}. 

\subsection{\sc Exact categories} Let $\cA$ be an exact $k$-category, that is an extension-closed full subcategory of an abelian $k$-category $\mk A$. Given objects $X, Y\in \cA$, we write $\Ext^i_{\mathcal{A}}(X, Y)=\Ext^i_{\mk A}(X, Y)$ for all integers $i\ge 0$. And one says that $\cA$ is {\it Ext-finite} if $\Ext^i_{\mathcal{A}}(X, Y)$ is finite dimensional for all $X, Y\in \cA$ and $i\ge 0$.
An object $P$ in $\cA$ is called 
{\it Ext-projective} if every short exact sequence $\xymatrixcolsep{18pt}\xymatrix{0\ar[r] &X\ar[r] &Y\ar[r] & P\ar[r] &0}$ in $\cA$ splits. The 
{\it Ext-injective objects} in $\cA$ are defined dually. If $\cA$ is abelian, then the projective objects and the injective objects in $\cA$ coincide with the Ext-projective objects and the Ext-injective objects, respectively. 

An epimorphism $f: X\to Y$ in $\cA$ is called {\it superfluous} if a morphism $g: M\to X$ in $\cA$ such that $f\circ g$ is an epimorphism is an epimorphism, and a monomorphism $f: X\to Y$ is called {\it essential} if a morphism $h: Y\to N$ in $\cA$ such that $h\circ f$ is a monomorphism is a monomorphism. Let $X$ be an object in $\cA$. A superfluous epimorphism $f: P\to X$ with $P$ projective in $\cA$ is called a {\it projective cover} of $X$ in $\cA$, and an essential monomorphism $g: X\to I$ with $I$ injective in $\cA$ is called an {\it injective envelope} of $X$ in $\cA$.  The following statement is well-known; see \cite[(3.4)]{HKKS}.

\begin{Lemma}\label{pc-ie-min}

Let $\cA$ be an exact $k$-category.

\begin{enumerate}[$(1)$]

\vspace{-2pt}

\item An epimorphism $f: P\to X$ with $P$ projective in $\cA$ is a projective cover of $X$ in $\cA$ if and only if $f$ is right minimal.

\item A monomorphism $g: X\to I$ with $I$ injective in $\cA$ is an injective envelope of $X$ in $\cA$ if and only if $g$ is left minimal.

\end{enumerate} \end{Lemma}



\vspace{-2pt}

Recall that a short exact sequence 
$\xymatrix{0\ar[r] &X\ar[r]^f &Y\ar[r]^g &Z\ar[r] &0}$ in $\cA$ is called an {\it almost split sequence} if $f$ is minimal left almost split and $g$ is minimal right almost split. In this case, one calls $X$ the {\it starting term} and $Z$ the {\it ending term}, and we write  $X=\tau Z$ and $Z=\tau^-\!X$; see \cite{AuR2}. 


We shall say that $\cA$ {\it has almost split sequences on the left} if every strongly indecomposable and non Ext-injective object is the starting term of an almost split sequence, and $\cA$ {\it has almost split sequences on the right} if every strongly indecomposable and non Ext-projective object is the ending term of an almost split sequence, and finally, $\cA$ {\it has almost split sequences} if it has almost split sequences on the left and on the right. 


\subsection{\sc Almost split triangles} \vspace{-4pt} Let $\mathcal{T}$ be a triangulated $k$-category with translation functor $[1]$. An exact triangle $\xymatrix{X \ar[r]^f &Y \ar[r]^g & Z \ar[r]^-\delta &X[1]}$ in $\mathcal A$ is called {\it almost split} if $f$ is minimal left almost split and $g$ is minimal right almost split; see \cite{Ha1}. In this case, one calls $X$ the {\it starting term} and $Z$ the {\it ending term}, and we write $X=\tau Z$ and $Z=\tau^-\!X$.
We say that $\mathcal{T}$ {\it has almost split triangles on the right} if every strongly indecomposable object is the ending term of an almost split triangle, $\mathcal{T}$ {\it has almost split triangles on the left} if every strongly indecomposable object is the starting term of an almost split triangle, and $\mathcal{T}$ {\it has almost split triangles} if it has almost split triangles on the right and on the left.

\subsection{\sc Derived categories} Let $\cA$ be a strictly full additive subcategory of an abelian category $\mk A$. We denote by $C(\cA)$ the additive category of complexes over $\cA$ with shift functor $[1]$, and by $C^+(\cA)$, $C^-(\cA)$ and $C^b(\cA)$ the full additive subcategories of $C(\cA)$ of bounded-below complexes, of bounded-above complexes and of bounded complexes, respectively. Given a complex $M^\cdt\in C(\mk A)$, a quasi-isomorphism $f^\pdt: P^\pdt\to M^\cdt$ in $C(\mk A)$ with $P^\pdt$ a complex of projective objects is called a {\it projective resolution} of $\m M^\cdt$, and a quasi-isomorphism $g^\cdt: M^\cdt\to I^\pdt$ in $C(\mk A)$ with $I^\cdt$ a complex of injective objects is called a {\it injective coresolution} of $M^\cdt$.


Fix $*\in \{b, +, -\}$. Endowed with the induced shift functor $[1]$, the quotient $K^*(\cA)$ of $C^*(\cA)$ modulo the null-homotopic morphisms is a triangulated category with exact triangles given by the mapping cones of morphisms; see \cite[(III.2.1.1)]{Mil}. Clearly, $K^*(\cA)$ is a full triangulated subcategory of $K^*(\mk A)$. A morphism in $K^*(\cA)$ is called a {\it quasi-isomorphism} if it is a quasi-isomorphism in $K^*(\mk A)$. 
By the same argument used in \cite[(III.3.1.1), (III.3.1.2)]{Mil}, we see that the quasi-isomorphisms in $K^*(\cA)$ form a localizing class compatible with the triangulation of $K^*(\cA).$ Thus,
the localization $D^*(\cA)$ of $K^*(\cA)$ at quasi-isomorphisms is a triangulated category; see \cite[(II.1.6.1)]{Mil}, called a {\it derived cateogry} of $\cA$. We shall say that $\cA$ has {\it enough $\mk A$-projective objects} if every object $X$ in $\cA$ admits an epimorphism $f: P\to X$ in $\mk A$, where $P\in \mathcal A$ is projective in $\mk A$. 

\begin{Prop}\label{Der_add}

Let $\cA$ be a strictly full additive subcategory of an abelian cate\-gory $\mk A$. If $\cA$ has enough $\mk A$-projective objects, then $D^b\m(\cA)$ can be regraded as a full triangulated subca\-tegory of $D^b\m(\mk A)$.

\end{Prop}

\noindent{\it Proof.} Clearly, $K^b(\cA)$ is a full triangulated subcategory of $ K(\mk A).$ The inclusion functors $q: K^b(\cA)\to K(\mk A);$ $i: K^b(\cA)\to K^b(\mk A)$ and $j: K^b(\mk A)\to D(\mk A)$ induce triangle-exact functors $q^D: D^b(\cA)\to D(\mk A);$ $i^D: D^b(\cA)\to D^b(\mk A)$ and $j^D: D^b(\mk A)\to D(\mk A)$ such that $q^D=j^D\circ i^D$. It is well-known that $j^D$ is fully faithful; see \cite[(III.3.4.5)]{Mil}. If $\cA$ has enough $\mk A$-projective objects, then $q^D$ is fully faithful; see \cite[(1.11)]{BaL}, and so is $i^D$. The proof of the proposition is completed.


\section{\sc Categories of graded modules} 

The objective of this section is to study generalities concerning graded modules over graded algebras given by locally finite quivers with homogeneous relations. 
The results obtained in this section will be needed not only in the following sections of this paper but also in future study of Koszul algebras given by locally finite quivers.


\vspace{2pt}

Throughout this section let $\La=kQ/R$, where $Q$ is a locally finite quiver and $R$ is a homogeneous relation ideal of $kQ$. Then $\La$ is a positively graded $k$-algebra with grading $\La=\oplus_{i\ge 0}\La_i$, where $\La_i=\{\bar{\gamma\hspace{1.6pt}} \mid \gamma\in kQ_i\}$. For convenience, set $\La_i=0$ for $i<0$. Write $J=\oplus_{i\ge 1} \La_i$, which is a graded two-sided ideal of $\La$. 
The opposite algebra $\La^{\rm o}$ is also positively graded as $\La^{\rm o}=\oplus_{i\ge 0} \La_i^{\rm o}$, where $\La_i^{\rm o}=\{\bar{\gamma\hspace{1.8pt}}^{\hspace{-1pt}\rm o} \mid \gamma\in kQ_i\}$, for all $i\ge 0$. Note that this grading for the opposite algebra is different from the classical one; see \cite[(1.2.4)]{NCF2}.

\subsection{\sc Graded modules.} A left $\La$-module $M$ is {\it unitary} if \vspace{1pt} $M=\sum_{x\in Q_0}e_xM$ and {\it graded} if $M=\oplus_{i\in \mathbb{Z}} M_i$, where the $M_i$ are $k$-spaces such that $\La_jM_i\subseteq M_{i+j}$, for all $i,j\in \Z$. 
Let $M$ be a graded unitary left $\La$-module. Then $M=\oplus_{i\in \mathbb{Z}; x\in Q_0}\, M_{i}(x)$ as a $k$-space, where 
$M_{i}(x)=e_xM_i,$ called the $(i,x)$-{\it piece} of $M$. Given $u\in e_y\La_je_x$, we shall write $M(u): M_i(x)\to M_{i+j}(y)$ for the $k$-linear map given by the left multiplication by $u$. An element $m\in M$ is called {\it homogeneous} of degree $i$ if $m\in M_i$ and {\it pure} if $m\in M_i(x)$ for some $(i,x)\in\mathbb{Z}\times Q_0$. A $\La$-submodule $N$ of $M$ is called {\it graded} if $N=\sum_{i\in \mathbb Z}(M_i\cap N).$ In this case, $N$ is graded as $N=\oplus_{i\in \mathbb Z} N_i,$ where $N_i=M_i\cap N$, such that if $m=\sum_{(i,x)\in \mathbb{Z}\times Q_0}m_{i,x}\in N$ \vspace{1pt} with $m_{i,x}\in M_i(x)$, then $m_{i,x}\in N$ for all $(i,x)\in \Z\times Q_0.$ 



\vspace{1pt}

Let $M, N$ be graded unitary left $\La$-modules. A $\La$-linear morphism $f: M\to N$ is {\it graded} if $f(M_{i})\subseteq N_{i}$ for all $i\in \Z$. In this case, $f$ restricts to $k$-linear maps $f_i: M_i\to N_i$ and $f_{i,x}: M_i(x)\to N_i(x)$ such that $f=\oplus_{i\in \Z} f_i =\oplus_{(i,x)\in \Z\times Q_0}f_{i,x}$. 
Conversely, given $k$-linear maps $f_{i,x}: M_i(x)\to N_i(x)$ with $(i, x)\in \Z\times Q_0$, the $k$-linear map $f=\oplus_{(i,x)\in \Z\times Q_0}f_{i,x}: M\to N$ is graded $\La$-linear if and only if $uf_{i,x}(m)=f_{i+j, y}(um)$, for $m\in M_i(x)$ and $u\in e_y\La_j e_x$ with $i,j\in \Z$ and $x,y\in Q_0$. 




The graded unitary left $\La$-modules together with the graded $\La$-linear morphisms form an
abelian $k$-category, which will be written as $\GrLa$. The morphism spaces and the extension groups in $\GrLa$ will be written respectively as ${\rm GHom}_\mathit\Lambda(M, N)$ and ${\rm GExt}^i_{\mathit \Lambda}(M, N)$. Moreover, we put ${\rm GEnd}_\mathit\Lambda(M)={\rm GHom}_\mathit\Lambda(M, M).$ The following statement is evident; compare \cite[page 20]{NCF2}.

\vspace{-.5pt}

\begin{Lemma}\label{pro-cop}

Let $\La=kQ/R$ be a graded algebra with $Q$ a locally finite quiver. If $\{M_{\sigma}\}_{\sigma\in\mathit\Sigma}$ is a family of modules in $\GrLa$, then

\begin{enumerate}[$(1)$]

\vspace{-2pt}

\item there exists a direct sum $M=\oplus_{\sigma\in \it\Sigma} M_{\sigma},$ 
defined by $M_{i}(x)=\oplus_{\sigma\in \it\Sigma} (M_\sigma)_{i}(x)$ for all $(i, x)\in \Z\times Q_0$. 


\item there exists a product $N=\Pi_{\sigma\in \it\Sigma} M_{\sigma},$ 
defined by
$N_{i}(x)=\Pi_{\sigma\in \it\Sigma} (M_\sigma)_{i}(x),$ for all $(i, x)\in \Z\times Q_0$. 

\end{enumerate}\end{Lemma}



Let $M\in\GrLa$ with $M=\oplus_{i\in \Z}M_i=\oplus_{i\in \Z; x\in Q_0}M_i(x)$. One says that $M$ is {\it bounded above} if $M_{i}=0$ for $i\gg 0,$  {\it bounded below} if $M_{i}=0$ for $i\ll 0$ and {\it bounded} if $M_{i}=0$ for all but finitely many $i\in \Z.$ The full subcategories of $\GrLa$ of bounded below modules and of bounded above modules will be written as 
${\rm GMod}^+\hspace{-3.6pt}\La$ and ${\rm GMod}^{-}\hspace{-3pt}\La$, respectively. Moreover, 
$M$ is called {\it locally finite dimensional} if $M_i$ is finite dimensional for all $i\in\mathbb{Z}$ and {\it piecewise finite dimensional} if $M_i(x)$ is finite dimensional for all $(i, x)\in\mathbb{Z}\times Q_0.$ We shall denote by ${\rm gmod}\La$ 
the full subcategory of $\GrLa$ of piecewise finite dimensional modules.


Let $V\in {\rm Mod}k$. Setting $(M\otimes V)_i=M_i\otimes V$, we obtain 
a graded module $M\otimes V=\oplus_{i\in \mathbb{Z}} (M\otimes V)_i\in \GrLa$.
Let $s\in  \Z$. The {\it grading} $s$-{\it shift} $M\sla s\sra$ of $M$ is defined by $M\sla s\sra_{i}=M_{i+s}$ for all $i\in \Z$. 
For a morphism $f: M\to N$ in $\GrLa$, the {\it grading} $s$-{\it shift} $f\nla s\nra:M\nla s\nra\to N\nla s\nra$ of $f$ is defined by $f\nla s\nra_i=f_{i+s},$ for all $i\in \Z$. 
It is clear that $(M\otimes V)\sla s\sra=M\sla s\sra\otimes V,$ for all $s\in \Z$ and $V\in \Mod k$. 

\subsection{\sc The duality $\mk D$} In the classical graded setting, there exists a duality for locally finite dimensional graded modules given by applying componentwise the functor $D$; see \cite[page 70]{Mar}. In our setting, we shall apply the functor $D$ piecewise in order to obtain a duality for piecewise finite dimensional graded modules. 


Given $M\in \GrLa$, we define $\mf{D}M = \oplus_{(i,x)\in \mathbb{Z}\times Q_0}\, D(M_{-i}(x))\in \GrLa^{\rm o}$, whose left $\La$-multiplication is such, for $\varphi\in D(M_{-i}(x))$ and $u\in e_x\La_je_y$, that $u^{\rm o}\cdot \varphi=\varphi \circ M(u),$ that is, $(u^{\rm o} \m \cdot \m \varphi)(m)=\varphi(u m),$ for $m\in M_{-i-j}(y).$ 
In particular, $(\mf{D}M)_i=\oplus_{x\in Q_0} D(M_{-i}(x))$ and $(\mf{D}M)_i(x)=D(M_{-i}(x))$ for $i\in \mathbb{Z}$ and $x\in Q_0$. Given a morphism $f: M\to N$ in $\GrLa$, we define a morphism $\mf{D}f: \mf{D}N\to \mf{D}M$ in $\GrLa^{\rm o}$ by setting $(\mf{D}f)_{i,x}=D(f_{-i, x}),$ for all $(i, x)\in \Z\times Q_0$. This clearly yields a contravariant functor $\mf{D} : \GrLa \to \GrLa^{\rm o}.$ 

\vspace{-.5pt}

\begin{Lemma}\label{mkD}

Let $\La=kQ/R$ be a graded algebra with $Q$ a locally finite quiver. Consider $M\in \GrLa$ and $V\in{\rm Mod}k$.

\begin{enumerate}[$(1)$]

\vspace{-2pt}

\item Given 
$s\in \Z$, we have $\mf{D}(M\nla s\nra)=(\mf{D} M)\nla -s\sra$.

\item 
There exists a natural monomorphism $\rho: M\to \mk D^2\!M$ in $\GrLa$, 
which is an isomorphism in case $M\in {\rm gmod}\La.$

\item There exists a binatural morphism 
$\theta: \mf{D}M\otimes DV \to \mf{D}(M\otimes V)$ in $\GrLa^{\rm o}$, 
which is an isomorphism in case $M\in {\rm gmod}\La$ or $V\in{\rm mod}k$.

\end{enumerate} \end{Lemma}

\noindent{\it Proof.} Statement (1) is evident. For $(i,x)\in \Z\times Q_0$, we have a canonical $k$-linear monomorphim $\rho_{i,x}: M_{i}(x)\to D^2(M_{i}(x))=(\mk D^2M)_i(x)$. 
Given $m\in M_i(x)$, $u\in e_y\La_je_x$ and $f\in D(M_{i+j}(y))$, we have \vspace{-3pt}
$$\rho_{i+j,y}(um)(f)=f(um)=(u^{\rm o}\cdot f)(m)=\rho_{i,x}(m)(u^{\rm o} \cdot f)=(u\cdot \rho_{i,x}(m))(f). \vspace{-3pt} $$ 

That is, $\rho_{i+j,y}(um)=u\, \rho_{i,x}(m)$. Thus, $\rho=\oplus_{(i,x)\in \Z\times Q_0} \rho_{i,x}: M\to \mk D^2\hspace{-1pt}M$ is a monomorphism in $\GrLa$, which is clearly natural in $M$. If $M\in {\rm gmod}\La$, then $M_i(x)\in {\rm mod} k$, and hence, $\rho_{i,x}$ is a $k$-linear isomorphism, for all $(i,x)\in \Z\times Q_0.$ That is, $\rho$ is an isomorphism. This establishes Statement (2). 
 
Next, given $(i, x)\in\Z\times Q_0,$ we have $(\mf{D}M\otimes DV)_i(x)=D(M_{-i}(x))\otimes DV$ and $(\mf{D}(M\otimes V))_{i}(x)=D(M_{-i}(x)\otimes V).$ Let $\theta_{i,x}: D(M_{-i}(x))\otimes DV\to D(M_{-i}(x)\otimes V)$ be the $k$-linear map as defined in 
Corollary \ref{Cor 1.2}(2). 
Given $g\in D((M_{-i}(x))$, $f\in DV$ and $u\in e_x\La_j e_y$, it is easy to verify that 
\vspace{-2pt} $$(u^{\rm o}\cdot\theta_{i,x}(g\otimes f))(m\otimes v)=\theta_{i+j,y}(u^{\rm o}\, (g\otimes f))(m\otimes v),
\mbox{ for } m\in M_{-i-j}(y), v\in V.\vspace{-2pt}$$

That is, $u^{\rm o}\cdot \theta_{i,x}(g\otimes f)=\theta_{i+j,y}(u^{\rm o}\, (g\otimes f))$. And consequently, we have
a morphism $\theta=\oplus_{(i,x)\in \mathbb{Z}\times Q_0}\theta_{i,x}: \mf{D}M\otimes DV \to \mf{D}(M\otimes V)$ in $\GrLa^{\rm o}$. It is a routine verification that $\theta$ is natural in $M$ and $V$. Finally, if $M\in {\rm gmod}\La$ or $V\in{\rm mod}k$ then, by Corollary \ref{Cor 1.2}(2), $\theta_{i,x}$ is a $k$-linear isomorphism for all $(i,x)\in \Z\times Q_0.$ That is, $\theta$ is an isomorphism in $\GrLa^{\rm o}$. The proof of the proposition is completed.


\vspace{2pt}

As a consequence of Lemma \ref{mkD}(2), we obtain our promised duality as follows.

\begin{Prop}\label{gr-duality}

Let $\La=kQ/R$ be a graded algebra with $Q$ a locally finite quiver. The contravariant functor 
$\mk D: {\rm GMod}\La \to {\rm GMod}\La^{\hspace{-.3pt} \rm o}\hspace{-2pt}$ is exact and restricts to a duality $\mf{D} : {\rm gmod}\La \to {\rm gmod}\La^{\hspace{-.3pt} \rm o}\hspace{-2pt}$.
\end{Prop}



The following statement says that $\mk D$ converts direct sums into direct products.

\begin{Prop}\label{cop-pro}

Let $\La=kQ/R$ be a graded algebra with $Q$ a locally finite quiver. Given $M_\sigma\in \GrLa$ with $\sigma\in \Sigma$, we have $\mk D(\oplus_{\sigma\in\mathit\Sigma} M_{\sigma})\cong \Pi_{\sigma\in\mathit\Sigma}\hspace{.6pt}\mk D(M_{\sigma}).$ 
    
\end{Prop}

\noindent{\it Proof.} Let $M_\sigma\in \GrLa$ with $\sigma\in \Sigma$. Write $M=\oplus_{\sigma\in\mathit\Sigma}M_{\sigma}$ and $N=\Pi_{\sigma\in\mathit\Sigma}\mk D(M_{\sigma}).$ 
Fix $(i,x)\in\Z\times Q_0$. Then $M_{i}(x) \m =\m \oplus_{\sigma\in\mathit\Sigma}(M_{\sigma})_{i}(x)$ and $N_{i}(x) \m = \m \Pi_{\sigma\in\mathit\Sigma}D(\m(M_{\sigma})_{-i}(x)\m)$. For each $\sigma\in \Sigma$, denote by $q_{\hspace{.4pt}\sigma}:(M_\sigma)_{i}(x)\to M_i(x)$ the canonical injection. Then, we have a canonical $k$-linear isomorphism \vspace{-2pt}
$$\Phi_{i,x}:(\mk{D}M)_i(x)=D(\oplus_{\sigma\in\mathit\Sigma}(M_{\sigma})_{-i}(x))\to \Pi_{\sigma\in\mathit\Sigma}D((M_{\sigma})_{-i}(x))=N_{i}(x) \vspace{-2pt}$$ such that $\Phi_{i,x}(f)=(f\circ q_{\hspace{.4pt}\sigma})_{\sigma\in\mathit\Sigma}$ for all $f\in (\mk{D}M)_i(x)$. Given $f\in (\mk{D}M)_i(x)$ and $u \m\in\m e_x\La_j e_y$, it is easy to verify that $u^{\rm o} \cdot \Phi_{i,x}(f)=\Phi_{i+j,y}(u^{\rm o}\cdot f)$. This yields an isomorphism $\Phi=\oplus_{(i,x)\in \mathbb{Z}\times Q_0} \Phi_{i,x}:\mk{D}M\to N$ in $\GrLa.$ The proof of the proposition is completed.

\subsection{\sc Graded projective modules.} A projective object in $\GrLa$ is called {\it graded projective.} For each $a\in Q_0$, we put $P_a=\La e_a=\oplus_{(i,x)\in \mathbb{Z}\times Q_0} e_x\La_ie_a$. Since $Q$ is a locally finite, $P_a$ is locally finite dimensional with $(P_a)_i=0$ for all $i<0$. To describe the graded morphisms starting from these modules, 
we fix some notation. Let $M\in \GrLa$. 
Given a pure element $m\in M_s(a)$ with $s\in \Z$ and $a\in Q_0$, the right multiplication by $m$ yields a graded $\La$-linear morphism $M[m]: P_a\sla -s \sra\to M.$ By definition, $M[m](e_a)=e_am=m.$

\begin{Prop}\label{proj-mor}
Let $\La=kQ/R$ be a graded algebra with $Q$ a locally finite quiver. Consider $P_a\sla -s \sra$ with $(s, a)\in \Z\times Q_0$ and $M\in \GrLa$. Then we have a natural $k$-linear isomorphism 
$ \eta : \GHom_{\mathit\Lambda}(P_a\sla -s \sra, M) \to  M_{s}(a): f\mapsto f(e_a),$ \vspace{1pt} whose inverse is given by  $\varphi: M_{s}(a) \to \GHom_{\mathit\Lambda}(P_a\sla -s \sra, M): m \mapsto M[m].$ 

\end{Prop}

\noindent{\it Proof.} Let $f\in \GHom_{\mathit\Lambda}(P_a\sla -s \sra, M).$ Observing that $e_a\in P_a\sla -s \sra_s(a)$, we see that $m=f(e_a)\in M_s(a)$ such that $f=M[m]$. It is easy to verify that $\eta^{-1}=\varphi$. The proof of the proposition is completed. 


\vspace{2pt}

The following statement is an immediate consequence of Proposition \ref{proj-mor}. 

\vspace{-.5pt}

\begin{Cor}\label{G-proj}
    
Let $\La=kQ/R$ be a graded algebra with $Q$ a locally finite quiver. Then, $P_a\sla -s \sra\otimes V$ is a graded projective module, for $(s, a)\in \Z\times Q_0$ and $V\in {\rm Mod}\hspace{.4pt}k$.

\end{Cor}


In the sequel, we shall denote by ${\rm GProj}\La$ the strictly full additive subcategory of $\GrLa$ generated by the $P_a\sla -s \sra\otimes V$ with $(s, a)\in \Z\times Q_0$ and $V\in {\rm Mod}\hspace{.4pt}k$, and by ${\rm gproj}\La$ the one generated by the $P_a\sla -s\mra$ with $(s, a)\in \Z\times Q_0$. We shall describe the morphsims in ${\rm GProj}\,\La$; compare \cite[(7.6)]{BaL}. Given $u \in e_a\La_{s-t}e_b=P_b\hspace{.3pt}\tla -t \tra_{s}$, in order to simplify the notation, we shall write the right multiplication by $u$ as \vspace{-1pt}$$P[u]: P_a\sla -s \sra\to P_b \hspace{.3pt}\tla -t \tra: v\mapsto vu.\vspace{-2pt}$$ 
Note that this notation does not distinguish $P[u]$ from its grading shifts. 

\vspace{-.5pt}

\begin{Prop}\label{rqz-pm}

Let $\La=kQ/R$ be a graded algebra with $Q$ a locally finite quiver. Consider $P_{a}\sla -s \sra \otimes V$ and $P_b \hspace{.3pt} \tla -t \tra \otimes W$ with $(s, a), (t, b) \in \Z \times Q_0 $ and $V, W\in \Mod k$. We have a $k$-linear isomorphism \vspace{-2pt}
$$\varphi: e_a\La_{s-t}e_b \otimes\m \Hom_k(V, W) \m\to\m \GHom_{\mathit\Lambda}(P_{a}\sla -s \sra \m\otimes V, P_b \hspace{.3pt} \sla -t \sra \otimes W)\m:\m  u\otimes f \m\mapsto\m P[u] \otimes f.\vspace{-2pt}$$ 

\end{Prop}

\noindent{\it Proof.} Clearly, we have a $k$-linear map $\varphi$ as stated in the proposition. Choose a $k$-basis $\{u_1, \ldots, u_n\}$ of $e_a\La_{s-t}e_b$. Consider $\omega\in {\rm Ker}(\varphi)$. Then, $\omega=\sum_{i=1}^n u_i\otimes f_i$, where $f_i\in \Hom_k(V, M)$. Given $v\in V$, we have $\varphi(\omega)(e_a\otimes v)=\sum_{i=1}^n u_i\otimes f_i(v)=0$, and hence, $f_i(v)=0$, for all $1\le i\le n$. Hence, $\omega=0.$ So, $\varphi$ is a monomorphism.

On the other hand, let $f\in \GHom_{\mathit\Lambda}(P_{a}\sla -s \sra \m\otimes\m V, P_b \hspace{.3pt} \tla -t \tra \m\otimes\m W)$. Given $v\in V$, observing that $e_a\otimes v\in P_{a}\sla -s \sra_{s}\otimes V$, we have $f(e_{a}\otimes v)=\sum_{i=1}^n u_i\otimes w_{i, v},$ for some unique $w_{i,v}\in W$. This yields $k$-linear maps $f_i: V\to W: v\mapsto w_{i,v}$, for $i=1, \ldots, n$, such that $f=\varphi(\sum_{i=1}^n u_i\otimes f_i)$. The proof of the proposition is completed.


\smallskip

The following statement is well-known in case $\La$ has an identity; see \cite[(2.2)]{NCF2}.

\begin{Prop}\label{GM-proj}

Let $\La=kQ/R$ be a graded algebra with $Q$ a locally finite quiver. Then ${\rm GMod}\La$ has enough projective objects.

\end{Prop}

\noindent{\it Proof.} Let $M\in {\rm GMod}\La$. For  $(i,x)\in\Z\times Q_0$, the multiplication map yields a graded morphism $f_{i,x}:P_x \,\nla -i \nra\otimes M_{i}(x)\to M$. Consider $P=\oplus_{(i,x)\in \Z\times Q_0}P_x\sla -i\sra\otimes M_{i}(x)$ with canonical injections $q_{s,a}:P_a\sla -s\sra\otimes M_{s}(a)\to P,$ for $(s,a)\in \Z\times Q_0$. Then, we have a graded morphism $f: P\to M$ such that $f\circ q_{s,a}=f_{s,a}$ for all $(s,a)\in\Z\times Q_0$. 
Clearly, $f: P\to M$ is an epimorphism with $P$ graded projective. The proof of the proposition is completed.

\subsection{\sc Graded injective modules.} An injective object in $\GrLa$ is called a {\it graded injective.} Given $a\in Q_0$, write $P_a^{\hspace{.4pt} \rm o}=\La^{\rm o} e_a\in \proj\La^{\rm o}$. Applying the duality $\mk D$, we obtain $I_a=\mf{D}P^{\hspace{.4pt}\rm o}_a\in {\rm gmod}\La$ with $(\m I_a\m )_{i}(x)=D(e_x \La_{-i}^{\rm o} e_a)$, for $(i,x)\in \Z\times Q_0$. Note that $I_a$ is locally finite dimensional with $(I_a)_{i}=0$ for $i>0$. Given $f \!\in\! (\m I_a\m )_{i}(x)$ and $u\in e_y \La_j e_x$, by definition, $u f\in (I_a)_{i+j}(y)=D(e_y \La_{-i-j}^{\rm o} e_a)$ such that
\vspace{-1.5pt} $$(u f)(v^{\rm o})=f(u^{\rm o} v^{\rm o}), \mbox{ for all } 
v\in e_a \La_{-i-j} e_y. \vspace{-1.5pt} $$ 
So, $I_a(u)\!=\!D(P_a^{\rm o}(u^{\rm o}))\m:\m (\m I_a\m )_i(x) \!\to \! (\m I_a\m )_{i+j}(y)$, where $P_a^{\rm o}(u^{\rm o})\m:\m P_a^{\rm o}(y)_{-i-j} \!\to \! P_a^{\rm o}(x)_{-i}$ is the left multiplication by $u^{\rm o}$.

\begin{Prop}\label{proj-inj}

Let $\La=kQ/R$ be a graded algebra with $Q$ a locally finite quiver. Consider $M\in \GrLa$ and $I_a\tla s \sra\otimes V$ with  $(s, a)\in  \Z \times Q_0$ and $V\in \Mod k$. Then, we have a natural $k$-linear isomorphism $$\psi:\GHom_{\mathit\Lambda}(M, I_a\tla s \sra\otimes V)\to \Hom_{k}(M_{-s}(a), V).$$

\end{Prop}

\noindent{\it Proof.} First, we have a $k$-linear isomorphism $\theta_a\!:\m\Hom_{k}(e_{a}\La^{\rm o}_{0}e_{a}, V) \!\to\! V\!\!:\! g \!\mapsto \! g(e_{a}).$ Given
$(i, x)\in \Z\times Q_{0}$, by Corollary \ref{Cor 1.2}(1), we have a $k$-linear isomorphism 
$$\sigma_{i,x} : I_a \tla s \sra_i(x)\otimes V=D(e_{x}\La^{\rm o}_{-i-s}e_{a})\otimes V\to \Hom_{k}(e_{x}\La^{\rm o}_{-i-s}e_{a}, V)$$ so that $\sigma_{i,x}(h\otimes v)(u^{\rm o})=h(u^{\rm o})v$, for $h\in D(e_{x}\La^{\rm o}_{-i-s}e_{a})$, $u\in e_a\La_{-i-s}e_x$ and $v\in V$.
Further, given any morphism $f : M\to I_a \tla s \sra\otimes V$ in $\GrLa$, we have a $k$-linear map $f_{-s,a}: M_{-s}(a)\to I_a\tla s \sra_{-s}(a)\otimes V$. This yields a natural $k$-linear map
$$\psi: \GHom_{\mathit\Lambda}(M, I_a\tla s \sra\otimes V)\to \Hom_{k}(M_{-s}(a), V) : f\mapsto \theta_a \circ\sigma_{-s,a}\circ f_{-s,a}. \vspace{-2pt}$$

Suppose that $\psi(f)=0$. Fix $(i, x)\in \Z\times Q_0$ and $m\in M_{i}(x)$. We may write $f_{i,x}(m)=\sum_{j=1}^{r}h_{j}\otimes v_{j}$, where $h_j \in D(e_x\La^{\rm o}_{-i-s}e_{a})$ and the $v_j$ are $k$-linearly independent in $V$. If $u\in e_{a}\La_{-i-s}e_{x}$, then $f_{-s,a}(u m)=u f_{i,x}(m)=\textstyle\sum_{j=1}^{r}u h_j\otimes v_j.$ Observing that $u h_j\in I_a\sla s \sra_{-s}(a)$, we obtain
\vspace{-3pt}
$$\textstyle 0 =\psi(f)(um)  = \sum_{j=1}^{r}\sigma_{-s,a}(u h_j\otimes v_j)(e_{a}) = \sum_{j=1}^{r} (u h_j)(e_a) v_j = \sum_{j=1}^{r}h_j(u^{\rm o})v_j. \vspace{-2pt} $$
Since the $v_j$ are $k$-linearly independent, $h_j(u^{\rm o})=0$ for $j=1, \ldots, r$. Hence, $h_j=0$, for $j=1, \ldots, r$. In particular, $f_{i,x}(m)=0$. Thus, $f_{i,x}=0$ for all $(i,x)\in \Z\times Q_0$. That is, $f=0$. So, $\psi$ is a monomorphism.

Consider now a $k$-linear map $g: M_{-s}(a)\to V$. Given $(i,x)\in \Z\times Q_0$, we shall define a $k$-linear map $f_{i, x}: M_i(x)\to I_a \tla s \sra_i(x)\otimes V$.
For any $m\in M_i(x)$, we have a $k$-linear map $g_{i,x}(m): e_{x}\La^{\rm o}_{-i-s}e_{a}\to V$ such that $g_{i,x}(m)(u^{\rm o})= g(u m)$ for all $u\in e_a\La_{-i-s}e_x.$ 
This yields a $k$-linear map $f_{i,x}: M_i(x)\to I_{a}\sla s \sra_i(x)\otimes V$, sending $m$ to $\sigma_{i,x}^{-1}(g_{i,x}(m)).$ In other words, $\sigma_{i,x}(f_{i, x}(m))=g_{i,x}(m),$ for all $m\in M_i(x).$  

Consider $v\in e_y\La_j e_x $ and $m\in M_i(x)$. Given $u\in e_{a}\La_{-i-j-s}e_y$, we obtain \vspace{-4pt} $$\sigma_{i+j,y}(f_{i+j,y}(v m))(u^{\rm o})=g_{i+j,y}(v m)(u^{\rm o})=
g(u v m)= g_{i,x}(m)((u v)^{\rm o}).\vspace{-2pt}$$

On the other hand, $\sigma_{i,x}^{-1}(g_{i, x}(m))=\sum_{p=1}^{r}h_p\otimes v_p$, for some $h_p\in D(e_{x}\La_{-i-s}^{\rm o} e_{a})$ and $v_p\in V$. Thus, $v f_{i,x}(m) = \sum_{p=1}^r (v h_p) \otimes v_p$ with
$v h_p \in D(e_{y}\La_{-i-j-s}^{\rm o} e_a)$. So \vspace{-3pt}
$$\textstyle \sigma_{i+j,y}(\m v f_{i,x}(m)\m)(\m u^{\rm o}\m) \!=\! \sum_{p=1}^s \! (\m v h\m_p)(u^{\rm o}\m) v\m_p \!=\! \sum_{p=1}^s\! \sigma_{i,x}(h\m_p\otimes v\m_p)((u v)^{\rm o}\m) \!=\! g_{i,x}(m\m)((u v)^{\rm o}\m).\vspace{-3pt}$$
Thus, $\sigma_{i+j,y}(v f_{i,x}(m)) = \sigma_{i+j,y}(f_{i+j,y}(vm))$. Hence, $f_{i+j,y}(v m)=v f_{i,x}(m)$. This yields a morphism $f=(f_{i, x})_{(i,x)\in \Z\times Q_0}: M\to I_{a}\sla s \sra\otimes V$ in $\GrLa$ such that $\psi(f)=g$. The proof of the proposition is completed.

\smallskip

The following statement is an immediate consequence of Proposition \ref{proj-inj}. 

\begin{Cor}\label{G-inj}
    
Let $\La=kQ/R$ be a graded algebra with $Q$ a locally finite quiver. Then  $I_a\mla s\mra\otimes V$ is  graded injective for any $(s, a)\in \Z\times Q_0$ and $V\in {\rm Mod}\hspace{.4pt}k$.

\end{Cor}

\smallskip

We denote by ${\rm GInj}\La$ the strictly full additive subcategory of $\GrLa$ generated by the $I_a\mla s\mra\otimes V$ with $(s, a)\in \Z\times Q_0$ and $V\in {\rm Mod}\hspace{.4pt}k$, and by ${\rm ginj}\La$ the strictly full additive subcategory of ${\rm GInj}\La$ generated by the $I_a\sla -s\mra$ with $(s, a)\in \Z\times Q_0$. To describe the morphisms in ${\rm GInj}\La$, we need to introduce some notation. Given $u \in e_a\La_{t-s} e_b$ with $s,t\in \Z$ and $a,b\in Q_0$, the right multiplication by $u^{\rm o}$ yields a graded $\La^{\rm o}$-linear morphism $P[u^{\rm o}]:\m P_b^{\rm o}\tla -t \tra \!\to\! P_a^{\rm o}\tla -s \tra$. Applying the duality $\mk D: {\rm gmod}\La^{\rm o}\to {\rm gmod}\La$, we obtain a morphism $I[u]=\mf{D}(P[u^{\rm o}]): I_a \sla s\sra  \to I_b\tla t \tra$
in $ {\rm GInj}\La$. Note that this notation does not distinguish $I[u]$ from its grading shifts. 

\begin{Prop}\label{Inj-Mor-2}

Let $\La=kQ/R$ be a graded algebra with $Q$ a locally finite quiver.
Consider $I_{a}\sla s \sra \otimes V$ and $I_b \hspace{.3pt} \tla t \tra \otimes W,$ for some $(s, a),(t, b)\in \Z\times Q_0$ and $V, W\in \Mod k$. Then, we have a $k$-linear isomorphism 
\vspace{-2pt}$$\phi: e_a\La_{t-s}e_b \m\otimes\m \Hom_k(V, W) \!\to \! \GHom_{\mit\Lambda}(I_{a}\sla s \sra \m\otimes\m V, I_b \hspace{.3pt} \tla t \tra \m\otimes\m W)\m: u \hspace{-1.5pt} \otimes\! f \mapsto I[u] \otimes\! f.\vspace{-2pt}$$ 

\end{Prop}

\noindent{\it Proof.} First, since $e_a\La_{t-s}e_b$ is finite dimensional, we have a $k$-linear isomorphism $\eta:e_a\La_{t-s}e_b\to D^2(e_b\La^{\rm o}_{t-s}e_a)$ such that $\eta(u)(g)=g(u^{\rm o})$ for all $u\in e_a\La_{t-s}e_b$ and $g\in D(e_b\La^{\rm o}_{t-s}e_a)$. 
Moreover, by Lemma \ref{Tensor}, we have a $k$-linear isomorphism $$\rho:D^2(e_b\La^{\rm o}_{t-s}e_a)\m\otimes\m \Hom_k(V, W)\to \Hom_k(D(e_b\La^{\rm o}_{t-s}e_a)\m\otimes\m V,W):\varphi\otimes f\mapsto \rho(\varphi\otimes f)$$ such that $\rho(\varphi\otimes f)(g\otimes v)=\varphi(g)f(v)$, for all $g\in D(e_b\La^{\rm o}_{t-s}e_a)$ and $v\in V$. Further, as did in the proof of Proposition \ref{proj-inj}, we consider two $k$-linear isomorphisms 
$\theta_b: \Hom_{k}(e_{b}\La^{\rm o}_{0}e_{b}, V)\to W: g\mapsto g(e_{b})$ and 
$$\sigma_{-t, b} : D(e_{b}\La^{\rm o}_{0}e_{b})\otimes W\to \Hom_{k}(e_{b}\La^{\rm o}_{0}e_{b}, W): 
g\otimes w\mapsto \sigma_{-t, b}(g\otimes w)$$ such that $\sigma_{-t,b}(g\otimes w)(e_b)=g(e_b)w$. 
Since $(I_a\sla s\sra\otimes V)_{-t}(b)=D(e_b\La^{\rm o}_{t-s}e_a)\otimes V,$ we obtain a $k$-linear isomorphism
$$\psi:\GHom_{\mit\Lambda}(I_{a}\sla s \sra \m\otimes\m V, I_b \hspace{.3pt} \tla t \tra \m\otimes\m W) \! \to \!\Hom_k(D(e_b\La^{\rm o}_{t-s}e_a)\m\otimes\m V,W)\! : \! h \m \mapsto \m \theta_b\circ \sigma_{-t,b}\circ h_{-t,b}.$$ 
In view of the above $k$-linear isomorphisms, we obtain a $k$-linear isomorphism $$\phi=\psi^{-1} \!\!\circ \m \rho \m\circ\m (\eta\otimes {\rm id}):e_a\La_{t-s}e_b \m\otimes\m \Hom_k(V, W) \!\to \! \GHom_{\mit\Lambda}(I_{a}\sla s \sra \m\otimes\m V, I_b \hspace{.3pt} \tla t \tra \m\otimes\m W).\vspace{-2pt}$$ 

Now, given $u\in e_a\La_{t-s}e_b$ and $f\in \Hom_k(V, W)$, it is a routine verification that $(\rho\circ (\eta\otimes {\rm id}))(u\otimes f)=\psi(I[u]\otimes f)$, that is, $\phi(u\otimes f)=I[u]\otimes f$. 
The proof of the proposition is completed.


\vspace{2pt}

In order to state a dual statement of Proposition \ref{proj-mor}, for each $a\in Q_0$, we denote by $e_a^\star$ the $k$-linear map in $(I_a)_0=D(ke_a)$ such that $e_a^\star(e_a)=1.$

\begin{Lemma}\label{soc-extend}

Let $\La=kQ/R$ be a graded algebra with $Q$ a locally finite quiver. Consider $M \in \GrLa$ and $I_a\sla s\sra$ with $(s,a)\in \Z\times Q_0$. Given $m\in M_{-s}(a)$, we have a graded morphism $f: M\to I_a\sla s\sra$ such that $f(m)=e_a^\star.$

\end{Lemma}

\noindent{\it Proof.} Fix $m\in M_{-s}(a)$ with $(s,a)\in \Z\times Q_0$. Considering the $k$-linear isomorphism $\theta_a:\Hom_{k}(e_a\La^{\rm o}_0e_a, k)\to k:g\mapsto g(e_a)$, \vspace{1pt} by Proposition \ref{proj-inj}, we have a $k$-linear isomorphism
$\psi:\GHom_{\mathit\Lambda}(M, I_a\tla s \sra)\to \Hom_{k}(M_{-s}(a), k):f\mapsto \theta_a\circ f_{-s,a}.$

Consider $h\in \Hom_k(M_{-s}(a), k)$ such that $h(m)=1$. Then, $\psi(f)=h$ for some 
$f\in \Hom_{\it\Lambda}( M, I_a\sla s\sra)$. So, $f_{-s,a}(m)(e_a)=\theta_a(f_{-s,a}(m))=\psi(f)(m)
=h(m)=1.$ Hence, $f_{-s,a}(m)=e^\star_a$. That is, $f(m)=e^\star_a$. The proof of the lemma is completed.



\smallskip

The following statement is well-known in case $\La$ has an identity; see \cite[(2.2)]{NCF2}.

\begin{Prop}\label{GM-pi}

Let $\La=kQ/R$ be a graded algebra with $Q$ a locally finite quiver. Then, ${\rm GMod}\La$ has enough injective objects.

\end{Prop}

\noindent{\it Proof.} Let $M\in {\rm GMod}\La$. Considering $\mk DM\in {\rm GMod}\La^{\rm o}$, we have a graded $\La^{\rm o}$-linear epimorphism \vspace{-1pt} $g:P^{\rm o}\to \mk DM,$ where $P^{\rm o}=\oplus_{(i,x)\in \Z\times Q_0}P^{\rm o}_x\sla -i\sra\otimes D(M_{i}(x))$; see (\ref{GM-proj}). Applying the exact functor $\mk D$ yields a graded monomorphism $\mk Dg:\mk D^2M\to \mk DP^{\rm o}.$ And by
Proposition \ref{gr-duality}(1), we obtain a graded monomorphism $h: M\to \mk DP^{\rm o}$. Now, 
we deduce from Proposition \ref{cop-pro} and Lemma \ref{mkD}(2) that \vspace{-1pt} $$\mk D P^{\rm o} \cong \Pi_{(i,x)\in \Z\times Q_0}\mk D(P^{\rm o}_x\sla -i\sra\otimes D(M_{i}(x)))\cong 
\Pi_{(i,x)\in \Z\times Q_0}I_x\sla i\sra\otimes D^2(M_{i}(x)),\vspace{-1pt}$$ 
which is graded injective. The proof of the proposition is completed.

\vspace{-4pt}

\subsection{\sc Graded semisimple modules.} A nonzero module in $\GrLa$ is called {\it graded simple} if it contains exactly two graded submodules. For each $a\in Q_0$, we put $S_a=P_a/Je_a$, which is clearly graded simple.

\begin{Prop}\label{graded-sim}

Let $\La=kQ/R$ be a graded algebra with $Q$ a locally finite quiver. Then, a module $S\in \GrLa$ is graded simple if and only if $S\cong S_{a}\sla i\sra$ for some $i\in \Z$ and $a\in Q_0$.

\end{Prop}

\noindent{\it Proof.} Let $S\in \GrLa$ be graded simple. Choose some $0\ne m \in S_i(a)$ with $(i, a)\in \Z\times Q_0$. Then, $S=\La m$ and $Jm=0$. By Proposition \ref{proj-mor}, we have a graded epimorphism $p: P_a\sla -i \sra \to S$ such that $p(e_a)=m$. This induces a graded epimorphism $\bar p: S_a\nla -i \nra = P_a\sla -i \sra / (JP_a)\sla -i\sra \to S.$ Since $S_a\nla -i \nra$ is graded simple, $S\cong S_a\mla -i\mra.$ The proof of the proposition is completed.

\smallskip

A nonzero module in $\GrLa$ is called  {\it graded semisimple} if it is a sum of graded simple modules. They can be characterized as follows.

\begin{Prop}\label{semi-rad}

Let $\La=kQ/R$ be a graded algebra with $Q$ a locally finite quiver. A nonzero module $M\in \GrLa$ is graded semisimple if and only if $JM=0;$ if and only if $\mk D M$ is graded semisimple.

\end{Prop}

\noindent{\it Proof.} Let  $M\in \GrLa$ be nonzero. If $M$ is graded semisimple, then $JM=0$ by Proposition \ref{graded-sim}. If $JM=0$, then $M=\oplus_{(i,a)\in \Z\times Q_0} M_i(a)$, where $M_i(a)$ is a  graded semisimple submodule of $M$. This proves the first equivalence. 

If $JM=0$, by definition, 
$J^{\rm o} \cdot \mk DM=0$, and hence, $\mk D M$ is graded semisimple. Suppose that $J M\ne 0$, say $um\ne 0$ for some $m\in M_i(x)$  with $(i,x)\in \Z\times Q_0$ and $u\in e_y\La_je_x$ with $j\ge 1$ and $y\in Q_0$. Then, $f(um)\ne 0$ for some $f\in D(M_{i+j}(y))$, that is, $(u^{\rm o} \cdot f )(m)=f(um)\ne 0$.
Thus, $J^{\rm o} \cdot \mk DM\ne 0$. Therefore, $\mk D M$ is not graded semisimple. The proof of the proposition is completed.



\vspace{-2pt}

\subsection{\sc Graded radical.} Let $M\in \GrLa$. A graded submodule of $M$ is called {\it graded maximal} if it is maximal among the graded submodules of $M.$ We shall describe all graded maximal submodules of $M$. An element $m\in M$ is called a \textit{top-element} if $m \in M_n(a) \backslash JM$ for some $(n,a)\in \Z\times Q_0$. In this case, we can find a $k$-subspace $L_{n,a}$ of $M_{n}(a)$, containing $M_{n}(a) \cap JM$, such that $M_{n}(a)=L_{n,a} \oplus k m$. Setting $L_{i,x}=M_i(x)$ for $(i,x)\in \Z\times Q_0$ with $(i, x)\ne (n, a),$ we obtain a $k$-subspace $L(m)=\oplus_{(i,x)\in \Z\times Q_0} L_{i,x}$ of $M$, which is of codimension one such that $m\not\in L(m)$. 

\begin{Lemma}\label{graded max}

Let $\La=kQ/R$ be a graded algebra with $Q$ a locally finite quiver. Consider a module $M$ in $\GMod\La$. A graded submodule $L$ of $M$ is graded maximal if and only if $L=L(m)$ for some top-element $m\in M;$ and in this case, $JM\subseteq L.$

\end{Lemma}

\noindent{\it Proof.} Let $M\in \GMod\La$. Consider a top-element $m\in M_n(a)$ with $(n,a)\in \Z\times Q_{0}$. We claim that $L(m)$ is a $\La$-submodule of $M$. Otherwise, there exist some $m'\in L_{i,x}$ and $u\in e_y\La_je_x$ with $i,j\in \Z$ and $x,y\in Q_0$ such that $um' \in M_{i+j}(y) \backslash L_{i+j, y}.$ By definition, $(i+j, y)=(n, a)$.  If $j=0$, then $(i,x)=(n,a)$ and $u\in e_a\La_0e_a=ke_a$, so $um'\in L_{i,x}=L_{i+j, y}$, absurd. If $j>0$, then $um'\in M_n(a) \cap JM\subseteq L_{n, a}=L_{i+j, y}$, 
a contradiction. This establishes our claim. Being of codimension one, $L(m)$ is graded maximal in $M$. Since $M_{n}(a) \cap JM\subseteq L_{n,a}$, we see that $JM\subseteq L(m)$. 

Let $L$ be a graded maximal submodule of $M$. Then, we have a graded simple module $M/L=\oplus_{(i,x)\in \mathbb{Z}\times Q_0}(M_{i}(x)+L)/L$. By Proposition \ref{graded-sim}, 
$M/L\cong S_{a}\mla -n\mra$ for some $(n,a)\in \Z\times Q_0$. Therefore, $(M_n(a)+L)/L=k(m+L)$ for some top element $m\in M_n(a)\backslash L_n(a)$, and $L_{i}(x)=M_{i}(x)$ for all $(i,x)\in \Z\times Q_0$ with $(i, x)\ne (n, a)$. As a consequence, $M_n(a)=L_n(a)+km$ and $M_{n}(a)  \cap JM \subseteq L_n(a)$. Since  $m\not\in L_n(a)$, we see that $M_n(a)=L_n(a)\oplus km$. In view of the above construction, $L=L(m)$. The proof of the lemma is completed.


\vspace{2pt}

The {\it graded radical} $\rad M$ of $M$ is the intersection of all graded maximal submodules of $M$. 
A graded submodule $N$ of $M$ is called {\it graded superfluous} in $M$ if $N+L\ne M$ for any proper graded submodule $L$ of $M$. 

\begin{Prop}\label{gr-rad}

Let $\La=kQ/R$ be a graded algebra with $Q$ a locally finite quiver. If $M\in \GMod \La$, then ${\rm rad}M=JM$, which contains all graded superfluous submodules of $M$.

\end{Prop}

\noindent{\it Proof.} Let $M\in \GMod \La$. By Lemma \ref{graded max}, $JM\subseteq \rad M.$ Consider $m\in M\backslash J M$. Write $m=\sum_{(i,x)\in \mathbb{Z}\times Q_0} m_{i,x}$, where $m_{i,x}\in M_i(x)$. Then, $m_{n,a} \not\in JM$ for some $(n, a)\in \Z\times  Q_0$. By Lemma \ref{graded max}, we have a graded maximal submodule 
$L(m_{n,a})$ of $M$ with $m_{n,a}\notin L(m_{n,a})$. Then, $m\not\in L(m_{n,a})$, and hence, $m\notin {\rm rad} M$. So, ${\rm rad} M= JM$. Suppose that $N$ is a graded superfluous submodule of $M$ with $N\not\subseteq {\rm rad}M$. Then,  $N\not\subseteq L$, for some graded maximal submodule $L$ of $M$. Thus $N+L=M,$ a contradiction. The proof of the proposition is completed.


\vspace{2pt}

\noindent{\sc Remark.} In case $Q$ is finite, it is known that $\rad M=JM$ for modules $M$ in ${\rm GMod}^-\La$; see, for example, \cite[Page 70]{Mar}.


\vspace{2pt}

As an immediate consequence of Proposition \ref{gr-rad}, we obtain the following statement, which is well-known in case $Q$ is finite. 

\begin{Cor}\label{proj_rad}

Let $\La=kQ/R$ be a graded algebra with $Q$ a locally finite quiver. Then, $\rad(_{\mathit\Lambda}\La)=J$ and $\rad P_a=Je_a$ for all $a\in Q_0$.
    
\end{Cor}





We have a sufficient condition for $\rad M$ to be graded superfluous in $M$, which is known in case $Q$ is finite; see \cite[Page 70]{Mar}. 

\begin{Prop}\label{bounded-b}

Let $\La=kQ/R$ be a graded algebra with $Q$ a locally finite quiver. If $M\in \GMod^+\!\!\La$, then ${\rm rad} M$ is graded superfluous in $M$. 

\end{Prop}

\noindent{\it Proof.} Let $M\in \GMod^+\!\!\La$. Assume that ${\rm rad}M+N=M$, where $N$ is a graded submodule of $M$ with $N\ne M$. We may find a minimal $s$ such that $N_s\ne M_s$. Choose $m\in M_s(a) \backslash N_s(a)$ for some $a\in Q_0$. By the minimality of $s$, we see that $M_s(a)\cap JM\subseteq N_s(a).$ In particular, $m$ is a top-element. It is easy to see that there exists a $k$-subspace $L_{s,a}$, containing $N_s(a)$, of $M_s(a)$ such that $M_s(a)=L_{s,a}\oplus km$. Since $M_s(a)\cap JM\subseteq L_{s,a}$, by Lemma \ref{graded max}, we may construct a graded maximal submodule $L(m)$ of $M$. Since $N_s(a)\subseteq L_{s,a}$, we have $N\subseteq L(m)$, and consequently, $M=L(m),$ absurd. The proof of the proposition is completed.


\vspace{2pt}

\noindent{\sc Example.} Let $M$ be a graded module over $k[x]$, which is illustrated as follows: 

\vspace{3pt} 

$\xymatrix{&\cdots \ar[r] & v_{-n} \ar[r]^{x} & \cdots  \ar[r]  & v_{-2} \ar[r]^{x} & v_{-1} \ar[r]^{x} & v_0 & \ar[l]_{x} u_{-1},}$ 

\vspace{2pt} 

\noindent where $x\cdot v_0=0$. By Proposition \ref{gr-rad}, $\rad M=k \{\hspace{.5pt} \ldots, v_{-n}, \ldots, v_{-2}, v_{-1}, v_0 \hspace{.5pt} \} \ne M.$ Observe that $M=\rad M+N$, where $N=k\langle u_{-1}, v_0\rangle$ is a graded submodule of $M$. Thus, $\rad M$ is not graded superfluous in $M$.


\vspace{2pt}

Given $M\in \GrLa$, we put ${\rm top}M=M/\rad M$, called the {\it graded top} of $M$. 
The following statement is known in case $Q$ is finite; see \cite[Page 70]{Mar}.

\begin{Cor}\label{gra-JR}

Let $\La=kQ/R$ be a graded algebra with $Q$ a locally finite quiver. If $M\in \GMod^+\hspace{-3pt}\La$ is nonzero, then ${\rm top}M$ is graded semisimple.

\end{Cor}

\noindent{\it Proof.} Let $M\in \GMod^+\hspace{-3pt}\La$ be nonzero. By Proposition \ref{bounded-b}, $\rad M$ is graded superfluous in $M$. In particular, ${\rm top}M\ne 0$. Since $\rad M=JM$; see (\ref{gr-rad}),  ${\rm top}M$ is graded semisimple by Proposition \ref{semi-rad}. The proof of the corollary is completed. 


\vspace{2pt}

\noindent{\sc Remark.} Corollary \ref{gra-JR} includes the graded version of Nakayama Lemma, which is known for positively graded algebras with an identity; see \cite[(2.9.2)]{NCF2}. 

\vspace{-2pt}

\subsection{\sc Finitely generated modules}  
Note that every finitely generated module in $\GMod \La$ is generated by finitely many pure elements. This fact leads to the following notion.

\begin{Defn}\label{top-basis}

Let $\La\!=\!kQ/R$ be a graded algebra with $Q$ a locally finite quiver. Given $M \m\in \m \GrLa$, a set $\{m_1, \ldots, m_r\}$ of pure elements in $M$ is called a {\it top}-{\it \hspace{-.4pt}basis} if $\{m_1+{\rm rad}M, \ldots, m_r+{\rm rad}M\}$ is a $k$-basis of ${\rm top} M$ and $M=\La m_1+ \cdots + \La m_r.$ 

\end{Defn}

\begin{Prop}\label{fg-iff}

Let $\La=kQ/R$ be a graded algebra with $Q$ a locally finite quiver. Given $M\in \GMod\La$, the following statements are equivalent$\,:$

\begin{enumerate}[$(1)$]

\vspace{-2pt}

\item $M$ is finitely generated$\,;$

\item $M$ admits a finite top-basis$\,;$

\vspace{.5pt}
    
\item $M$ is bounded below and ${\rm top}M$ is finite dimensional.

\end{enumerate}
\end{Prop} 

\noindent{\it Proof.} Let $M\in \GrLa$ be nonzero. Assume that Statement (1) holds. Clearly, $M\in {\rm GMod}^+\hspace{-3.5pt}\La$. By Corollary \ref{gra-JR}, ${\rm top}M$ is graded semisimple. Being finitely generated,  ${\rm top}M$ is finite dimensional; see (\ref{graded-sim}). Thus, Statement (3) holds.

Suppose that Statement (3) holds. We may choose pure elements $m_1, \ldots, m_r$ in $M$ such that 
$\{m_1 + \rad M, \ldots, m_r + \rad M\}$ is a $k$-basis of ${\rm top}M$. In particular,
$M/\rad M=(\sum_{i=1}^r \La m_i+\rad M)/\rad M.$ \vspace{.5pt} Since $\rad M$ is graded superfluous in $M$; see (\ref{bounded-b}),  
$M=\sum_{i=1}^r \La m_i$. 
Thus, Statement (2) holds, and so does Statement (1). The proof of the proposition is completed.

\vspace{-2pt}

\subsection{\sc Graded projective cover} A superfluous epimorphism in $\GrLa$ is called {\it graded superfluous}, and a projective cover of a module in $\GrLa$ is called a {\it graded projective cover.}

\begin{Lemma}\label{proj_cov}

Let $\La=kQ/R$ be a graded algebra with $Q$ a locally finite quiver. 

\begin{enumerate}[$(1)$]
   
\item An epimorphism $f:M\to N$ in $\GrLa$ is graded superfluous if and only if ${\rm Ker}(f)$ is graded superfluous in $M;$ and in this case, $f^{-1}(\rad N)=\rad M$. 

\item An epimorphism $f: P \m\to \m M$ in ${\rm GMod}^{-\hspace{-3pt}}\La$ with $P$ graded projective is a graded projective cover of $M$ if and only if ${\rm Ker}(f) \subseteq \rad P.$

\end{enumerate}
\end{Lemma}

\noindent{\it Proof.} Statement (1) is easily adapted from the non-graded setting; see, for exa\-mple, \cite[(5.15)]{FWAF}. Consider an epimorphism $f: P \m\to \m M$ in ${\rm GMod}^{-\hspace{-3pt}}\La$ with $P$ graded projective. By Propositions \ref{gr-rad} and \ref{bounded-b}, $\rad P$ is the largest superfluous graded submodule of $P$. So, every graded submodule of $\rad P$ is superfluous in $P$. Now, Statement (2) follows from Statement (1). The proof of the lemma is completed.

\smallskip

\noindent{\sc Example.} Given $a\in Q_0$, the canonical projection $p_a: P_a\to S_a$ is a graded projective cover of $S_a$.

\vspace{2pt} 

We are ready to construct a graded projective cover for every finitely generated graded module; compare \cite[(1.1)]{LiM}.

\begin{Prop}\label{p-cover}

Let $\La=kQ/R$ be a graded algebra with $Q$ a locally finite quiver. A module $M\in {\rm GMod}\La$ admits a graded projective cover \vspace{-2pt}  $$f \m : \m P_{a_{\hspace{-.8pt}1}} \m \sla -s_{\hspace{-.6pt}1}\nra \m\oplus \cdots \oplus\m P_{a_{\m r\m}} \nla -s_{\hspace{-.6pt}r} \nra \m \to \! M: e_{a_i}\mapsto m_i \vspace{-2pt} $$ if and only if $\{m_1,\dots, m_r\}$ with $m_i\in M_{s_i}(a_i)$ is a top-basis for $M.$

\end{Prop}

\noindent{\it Proof.} Let $\{m_1,\dots,m_r\}$ be a top-basis for $M$, where $m_i\in M_{s_i}(a_i)$ with $s_i\in \Z$ and $a_i\in Q_0$. In view of Proposition \ref{proj-mor}, we obtain a graded epimorphism \vspace{-2pt}$$f: P_{a_{\hspace{-.8pt}1}} \m \sla -s_{\hspace{-.6pt}1}\nra \m\oplus \cdots \oplus\m P_{a_{\m r\m}} \nla -s_{\hspace{-.6pt}r} \nra \to M: e_{a_i}\mapsto m_i.\vspace{-2pt}$$

Since $\{m_1+\rad M, \ldots, m_r+\rad M\}$ is $k$-linearly independent, it follows that ${\rm Ker}(f)\subseteq \rad(\oplus_{i=1}^r P_{a_{\hspace{-.8pt}1}} \m \sla -s_{\hspace{-.6pt}1}\nra).$
So 
$f$ is a graded projective cover of $M$; see (\ref{proj_cov}).

Suppose that $M$ has a graded projective cover as stated in the proposition. Then, $M=\sum_{i=1}^r \La m_i$ and ${\rm top}M = \sum_{i=1}^r k(m_i+\rad M).$ Let $\sum_{i=1}^r \lambda_i(m_i+\rad M)=0$, where $\lambda_i\in k$. Since $\rad P=f^{-1}(\rad M)$; see (\ref{proj_cov}), there exists some $u\in \rad P$ such that $f(u)=\sum_{i=1}^r\lambda_im_i=\sum_{i=1}^rf(\lambda_ie_{a_i})$. \vspace{1pt} Since ${\rm Ker}(f)\subseteq \rad P$; see (\ref{proj_cov}), we have $\sum_{i=1}^r\lambda_ie_{a_i}\in \oplus_{i=1}^r \rad (P_{a_{\hspace{-.8pt}i}} \m \sla -s_{\hspace{-.6pt}i}\nra \m)$. So, $\lambda_i=0$, for $i=1,\dots,r$. 
That is, $\{m_1,\dots,m_r\}$ is a top-basis for $M$. 
The proof of the proposition is completed.




\vspace{3pt}

A module $M$ in ${\rm gmod}\La$ is called {\it finitely presented} if it admits a {\it graded projective presentation} over ${\rm gproj}\La$, that is an exact sequence  \vspace{-5pt}
$$\xymatrix{P^{-1} \ar[r]^{d^{-1}} &P^{0} \ar[r]^{d^{\hspace{.4pt}0}}  & M \ar[r] & 0 }\vspace{-5pt} $$ in ${\rm gmod}\La$, where $P^0, P^{-1}\in {\rm gproj}\La.$ Such a graded projective presentation is called {\it minimal} if $d^{-1}$ and $d^{\hspace{.4pt}0}$ are both right minimal. Applying Proposition \ref{p-cover} and Schanuel's Lemma, we obtain the following statement.

\begin{Lemma}\label{mini-fpres}

Let $\La=kQ/R$ be a graded algebra with $Q$ a locally finite quiver. Then, every finitely presented module in ${\rm gmod}\La$ admits a minimal graded projective presentation over ${\rm gproj}\La$, which is unique up to isomorphism.

\end{Lemma}



\subsection{\sc Graded socle.} Let $M\in \GrLa$. The {\it graded socle} ${\rm soc} M$ of $M$ is the sum of all graded simple submodules of $M$. 
A graded submodule $N$ of $M$ is called {\it graded essential} in $M$ if $N\cap L\ne 0$ for any nonzero graded submodule $L$ of $M.$ 

\begin{Lemma}\label{sim-soc}

Let $\La=kQ/R$ be a graded algebra with $Q$ a locally finite quiver. If $M\in \GMod\La$, then $\soc M$ is contained in every essential graded submodule of $M$ such that $({\rm soc}\hspace{.3pt} M)_{i}(x)=\{m\in M_i(x) \mid Jm=0\}$ for all $(i,x)\in \Z\times Q_0$.

\end{Lemma}

\noindent{\it Proof.} Let $M\in \GrLa$. Assume that $L$ is an essential graded submodule of $M$. If $S$ is a graded simple submodule of $M$, then $S=L\cap S\subseteq L$. Thus, ${\rm soc}M\subseteq L.$ The second part of the statement follows from Proposition \ref{semi-rad}. The proof of the lemma is completed. 


In general, $\soc M$ is not necessarily graded essential in $M$. Nevertheless, we have the following sufficient condition for this to happen.

\begin{Lemma}\label{ess-soc}

Let $\La=kQ/R$ be a graded algebra with $Q$ a locally finite quiver. If $M \in  {\rm GMod}^-\hspace{-2.5pt}\La$, then ${\rm soc} M$ is graded essential in $M$.

\end{Lemma}

\noindent{\it Proof.} Let $M \in {\rm GMod}^-\hspace{-3pt}\La$. Consider  a nonzero graded submodule $N$ of $M$. Choose $0\ne m\in N_i$ for some integer $i$. Since $N\in {\rm GMod}^-\hspace{-3pt}\La$, there exists some $j\ge 0$ such that $\La_j m\ne 0$ but $\La_{j+1} m=0$. By Lemma \ref{sim-soc}, $ \La_jm\subseteq {\rm soc}M$. The proof of the lemma is completed.

\smallskip

\noindent{\sc Example.} Consider $\La=kQ/R$, where \vspace{-6pt} $$\begin{tikzpicture}[-{Stealth[inset=0pt,length=4.5pt,angle'=35,round,bend]},scale=.4]
\node (a) at (-4.5,0) {$Q:$};
 \draw (0,0.25) arc (15:345:0.8);
    \node (1) at (0:0.3){$1$};
    \node[left] at (-1.5,0.1){$\alpha$}; 
\draw (0.65,0)--node[above]{\scalebox{0.88}{\( \beta \)}} (2.8,0) node[right]{$2$};
 \end{tikzpicture}\vspace{-3pt} $$ and $R=\langle\hspace{1pt}\beta\alpha\hspace{0.7pt}\rangle$.
Then, $P_1=k\{ e_1, \bar\beta, \bar\alpha, \bar \alpha^2\m, \cdots\}$ with $\soc P_1=k\bar\beta$. Observe that $L=k\{ \bar\alpha, \bar \alpha^2\m, \cdots \}$ is a graded submodule of $P_1$ such that $L\cap \soc P_1=0$. So, $\soc P_1$ is not graded essential in $P_1$. 




\vspace{2pt}

The next statement describes the graded socle for modules in ${\rm ginj}\La$. 

\begin{Cor}\label{soc_ia}

Let $\La=kQ/R$ be a graded algebra with $Q$ a locally finite quiver. If $a\in Q_0,$ then $\soc I_a=ke_a^\star$, which is graded essential in $I_a$.

\end{Cor}

\noindent{\it Proof.} Fix $a\in Q_0$. Note that
$(I_a)_0=D(e_a\La_0^{\rm o}e_a)=k e^\star_a.$ Since $(I_a)_i=0$ for $i>0$, by Lemma \ref{sim-soc}, $e_a^\star\in \soc I_a$. So, $(\soc I_a)_0=k e^\star.$ Consider $0\ne f\in (I_a)_{-i}(x) 
$ for some $i>0$ and $x\in Q_0$. Then, $f(u^{\rm o})\ne 0$ for some $u\in e_a \La_i e_x$, that is, $(u\cdot f)(e_a)\ne 0$. By Lemma \ref{sim-soc}, $f\notin \soc I_a$. Thus, $\soc I_a=(\soc I_a)_0=ke_a^\star$. Moreover, by Lemma \ref{ess-soc}, $\soc I_a$ is graded essential in $I_a$. The proof of the corollary is completed.

\vspace{-0pt}

\subsection{\sc Finitely cogenerated modules} A module $M$ in $\GrLa$ is called {\it finite\-ly cogenerated} if $\soc M$ is finitely generated and graded essential in $M$. 

\begin{Defn}\label{soc-basis}

Let $\La=kQ/R$ be a graded algebra with $Q$ a locally finite quiver. Given $M\in \GrLa$, a set $\{m_1, \ldots, m_r\}$ of pure elements in $M$ is called a {\it soc}-{\it basis} if $\soc M$ has $\{m_1, \ldots, m_r\}$ as a $k$-basis and  
is graded essential in $M.$

\end{Defn}

Finitely cogenerated graded modules are characterized as follows.

\begin{Prop}\label{fcg-iff}

Let $\La=kQ/R$ be a graded algebra with $Q$ a locally finite quiver. Given $M\in \GrLa$, the following statements are equivalent$\,:$

\begin{enumerate}[$(1)$]

\vspace{-1pt}

\item $M$ is finitely cogenerated$\,;$

\item $M$ admits a finite soc-basis$\,;$

\item $M$ is bounded above and  $\soc M$ is finite dimensional.

\end{enumerate}

\end{Prop}

\noindent{\it Proof.} Let $M\in \GrLa$ be nonzero. Assume that Statement (1) holds. Then, $\soc M$ is finitely generated and graded semisimple. By Proposition \ref{graded-sim}, ${\rm soc}M$ is finite dimensional. Thus, $\soc M$ has a $k$-basis 
$\{m_1, \dots, m_r\}$, where $m_i\in M_{s_i}(a_i)$ with $(s_i, a_i)\in \Z\times Q_0$. Set $s=s_1+\dots+s_r$. Then $M_j \cap \hspace{.5pt}\soc M=0$ for $j>s$. Suppose that there exists $0\ne m\in M_p$ for some $p>s$. Since $\soc M$ is graded essential in $M$, there exists some $u\in \La_t$ with $t\ge 0$ such that $0\ne um\in  M_{t+p}\cap\soc M$, a contradiction. Hence, Statement (3) holds.

Suppose that Statement (3) holds. Then, $\soc M$ has a $k$-basis 
$\{m_1, \dots, m_r\}$, where the $m_i$ are pure elements in $M$. Since $\soc M$ is graded essential in $M$, by Lemma \ref{ess-soc}, $\{m_1, \dots, m_r\}$ is a soc-basis for $M$. Hence, Statement (2) holds. The proof of the proposition is completed. 

\smallskip

As an immediate consequence of Proposition \ref{fcg-iff} and Corollary \ref{soc_ia}(3), we obtain 
the following statement.

\begin{Cor}\label{inj_fcg}

Let $\La=kQ/R$ be a graded algebra with $Q$ a locally finite quiver. Then, every module in ${\rm ginj}\La$ is finitely cogenerated, and every graded submodule of a finitely cogenerated graded module is finitely cogenerated. 

\end{Cor}



\vspace{-3pt}

\subsection{\sc Graded injective envelope.} An essential monomorphism in $\GrLa$ is called {\it graded essential}, and an injective envelope of a module in $\GrLa$ is called a {\it graded injective envelope}.  

\begin{Lemma}\label{inj_env}

Let $\La=kQ/R$ be a graded algebra with $Q$ a locally finite quiver. 

\begin{enumerate}[$(1)$]

\item A monomorphism $f:M\to N$ in $\GrLa$ is graded essential if and only if ${\rm Im}(f)$ is graded essential in $N;$ and in this case, $\soc N=f(\soc M)$. 

\item A monomorphism $f: M\to I$ in ${\rm GMod}^+\hspace{-3pt}\La$ with $I$ graded injective is a graded injective envelope of $M$ if and only if ${\rm soc}I\subseteq {\rm Im}(f)$.

\end{enumerate}
\end{Lemma}

\noindent{\it Proof.} Statement (1) is easily adapted from the non-graded setting; see, for example, \cite[(5.13)]{FWAF}. By Lemmas \ref{sim-soc} and \ref{ess-soc}, $\soc I$ is the smallest essential graded submodule of $I$. Thus, every graded submodule of $I$ containing $\soc I$ is graded essential in $I$. Now, Statement (2) follows from Statement (1). The proof of the lemma is completed.

\smallskip

\noindent{\sc Example.} Given $a\in Q_0$, by Lemma \ref{soc-extend}, we have a graded monomorphism $q_a: S_a\to I_a$, sending $e_a+Je_a$ to $e_a^\star$. By Corollary \ref{soc_ia}, ${\rm Im}(q_a)=\soc I_a$, and by Lemma \ref{inj_env}, $q_a$ is a graded injective envelope of $S_a$. 

\smallskip

We are ready to construct a graded injective envelope for every finitely cogene\-rated graded module.

\begin{Prop}\label{i-envelope}

Let $\La=kQ/R$ be a graded algebra with $Q$ a locally finite quiver. A module $M$ in ${\rm GMod}\La$ admits a graded injective envelope $$g\m:\! M \!\to\! I_{a_1} \m\mla s_1 \mra \hspace{.3pt}\oplus \hspace{.3pt} \cdots \hspace{.3pt} \oplus \hspace{.3pt} I_{a_r} \m \mla s_r\mra:m_i\mapsto e_{a_i}^\star$$ if and only if $\{m_1, \dots, m_r\}$ with $m_i \in M_{-s_i}(a_i)$ is a soc-basis for $M.$

\end{Prop}

\noindent{\it Proof.} Suppose that $\{m_1, \dots, m_r\}$ is a soc-basis for $M,$ where $m_i \in M_{-s_i}(a_i)$. Then, $\soc M=km_1\oplus\cdots \oplus km_r.$ By Lemma \ref{soc-extend}, we have a graded monomorphism $$q: \soc M\to I_{a_1} \sla s_1 \sra \hspace{.3pt}\oplus \hspace{.3pt} \cdots \hspace{.3pt} \oplus \hspace{.3pt} I_{a_r}\sla s_r\sra=I:m_i\mapsto e_{a_i}^\star.$$ 
By Corollary \ref{soc_ia}, $\soc I=ke_{a_1}^\star\!\!\oplus\cdots \oplus ke_{a_r}^\star={\rm Im}(q).$ Since the inclusion map $h:\soc M\to M$ is graded essential; see (\ref{inj_env}), we have a graded monomorphism $g: M\to I$ such that $g\circ h=q$. Since ${\rm soc} I={\rm Im}(q) \subseteq {\rm Im}(g),$ by Lemma \ref{inj_env}(2), $g$ is a graded injective envelope of $M$.

Suppose that $g: M\to I_{a_1} \m\mla s_1 \mra \hspace{.3pt}\oplus \hspace{.3pt} \cdots \hspace{.3pt} \oplus \hspace{.3pt} I_{a_r} \m \mla s_r\mra=I$ is a graded injective envelope of $M$. By Lemma \ref{inj_env}, $ke_{a_1}^\star\!\!\oplus \cdots \oplus ke_{a_r}^\star=g(\soc M),$ where $e_{a_i}^\star\in I_{a_i} \m \mla s_i\mra_{-s_i}$. Thus, $e_{a_i}^\star=g(m_i)$ for some $m_i\in (\soc M)_{-s_i}(a_i)\subseteq M_{-s_i}(a_i)$. Since $g$ is a monomorphism, $M$ is bounded above and $\{m_1,\dots,m_r\}$ is a $k$-basis of $\soc M$. By Lemma \ref{ess-soc}, $\soc M$ is graded essential in $M$. 
The proof of the proposition is completed.


\smallskip

A module $M$ in ${\rm gmod}\La$ is called {\it finitely copresented} if it admits a {\it graded injective copresentation} over ${\rm ginj}\La$, that is an exact sequence  \vspace{-5pt}
$$\xymatrixcolsep{20pt}\xymatrix{0\ar[r] & M\ar[r]^{d^{\hspace{.4pt}0}} & I^0 \ar[r]^{d^1} & I^1,}\vspace{-5pt}$$ in ${\rm gmod}\La$ with $I^0, I^{1}\in {\rm ginj}\La.$ Such a graded injective copresentation is called {\it minimal} if $d^{\hspace{.4pt}^0}$ and $d^1$ are both left minimal. Applying Corollary \ref{inj_fcg}, Proposition \ref{i-envelope} and the dual of Schanuel's Lemma, we obtain the following statement.

\begin{Lemma}\label{mini-fpres}

Let $\La=kQ/R$ be a graded algebra with $Q$ a locally finite quiver. Every finitely copresented module in ${\rm gmod}\La$ admits a minimal graded injective copresentation over ${\rm ginj}\La$, which is unique up to isomorphism.

\end{Lemma}

\vspace{-5pt}

\subsection{\sc Krull-Schmidt subcategories} In this subsection, we provide several Hom-finite Krull-Schmidt subcategories of $\GrLa$, which will play an important role in our later study of almost split sequences and almost split triangles. 

\begin{Lemma}\label{f-GHom} Let $\La=kQ/R$ be graded algebra with $Q$ a locally finite quiver. If $P\in {\rm gproj}\La$ and $I\in {\rm ginj}\La$, then  $\GHom_{\mit\Lambda}(P,M)$ and $\GHom_{\mit\Lambda}(M,I)$ are finite dimensional, for all $M\in {\rm gmod}\La$.

\end{Lemma}

\noindent{\it Proof.} We shall only prove the first part of the statement. Let $P\in {\rm gproj}\La$. Then $P\cong \oplus_{i=1}^rP_{a_i}\sla -s_i\sra$, where $(s_i,a_i)\in \Z\times Q_0$. Given $M\in {\rm gmod}\La$, by Proposition \ref{proj-mor}, $\GHom_{\mit\Lambda}(P,M) 
\cong \oplus_{i=1}^rM_{s_i}(a_i).$ 
The proof of the lemma is completed.


\vspace{2pt}

The following statement exhibits some particular feature of the graded setting.

\begin{Lemma}\label{proj-KS}

Let $\La=kQ/R$ be a graded algebra with $Q$ a locally finite quiver. 

\begin{enumerate}[$(1)$]

\vspace{-7pt}

\item The category ${\rm gproj}\La$ is Hom-finite Krull-Schmidt and contains all finitely gene\-rated graded projective modules in $\GrLa$. 

\item  The category ${\rm ginj}\La$ is Hom-finite Krull-Schmidt and contains all finitely cogene\-rated graded injective modules in $\GrLa$. 

\end{enumerate}

\end{Lemma}

\noindent{\it Proof.} We shall only prove Statement (1). By Lemma \ref{f-GHom}, ${\rm gproj}\La$ is Hom-finite. Given $s\in\Z$ and $a\in Q_0$, by Proposition \ref{rqz-pm}, ${\rm GEnd}_\mit\Lambda(P_{a} \sla -s\mra)\cong ke_a$. Thus, ${\rm gproj}\La$ is Krull-Schmidt. If $M\in \GrLa$ is finitely generated and graded projective, then it has a graded projective cover $f: P\to M$ with $P\in{\rm gproj}\La$; see (\ref{p-cover}), and hence, $M\cong P$. 
The proof of the lemma is completed.

\vspace{3pt}

We write ${\rm gmod}^{+,\hspace{.6pt} b\hspace{-2.5pt}}\La$, ${\rm gmod}^{-,\hspace{.6pt} b\hspace{-2.5pt}}\La$ and ${\rm gmod}^{\hspace{.4pt}b}\hspace{-3pt}\La$ for the full subcategories of $\GrLa$ of finitely gene\-rated modules, of finitely cogenerated modules, and of finite dimensional modules, respectively. 
Clearly, they are all subcategories of ${\rm gmod}\La$. 

\begin{Lemma}\label{fg-KrS}

Let $\La=kQ/R$ be a graded algebra with $Q$ a locally finite quiver. 

\begin{enumerate}[$(1)$]

\vspace{-6.5pt}

\item The restricted functor $\mk D: {\rm gmod}^{+,\hspace{.5pt}b\hspace{-2.5pt}}\La\to {\rm gmod}^{-,\hspace{.5pt}b\hspace{-2.5pt}}\La^{\rm o}$ is a duality.

\vspace{.5pt}

\item Both ${\rm gmod}^{+,\hspace{.6pt} b\hspace{-2.5pt}}\La$ and 
${\rm gmod}^{-,\hspace{.6pt} b\hspace{-2.5pt}}\La$ are Hom-finite Krull-Schmidt extension-closed subcategories of $\GrLa$, whose intersection is ${\rm gmod}^{\hspace{.4pt}b}\hspace{-3pt}\La$.

\end{enumerate}

\end{Lemma}

\noindent{\it Proof.} (1) Let $M\in {\rm gmod}^{+,\hspace{.5pt}b\hspace{-2.5pt}}\La$. Then, we have a graded epimorphism $f:P\to M$ with $P\in{\rm gproj}\La$. Applying the duality $\mk D: {\rm gmod}\La \to {\rm gmod}\La^{\rm o}$; see (\ref{gr-duality}), we obtain a graded monomorphism $\mk Df:\mk DM\to \mk DP$ with $\mk DP$ in ${\rm ginj}\La^{\rm o}$. By Corollary \ref{inj_fcg}, $\mk DM \in {\rm gmod}^{-,\hspace{.5pt}b\hspace{-2.5pt}}\La^{\rm o}$. Dually, if $N\in {\rm gmod}^{-,\hspace{.5pt}b\hspace{-2.5pt}}\La^{\rm o}$, then $\mk D N\in {\rm gmod}^{+,\hspace{.5pt}b\hspace{-2.5pt}}\La$. 

\vspace{1pt}

(2) Clearly ${\rm gmod}^{+,\hspace{.6pt} b\hspace{-2.5pt}}\La$ is closed under direct summands and extensions, and by Lemma \ref{f-GHom}, it is Hom-finite. So, ${\rm gmod}^{+,\hspace{.6pt} b\hspace{-2.5pt}}\La$ is Krull-Schmidt. Then, by Statement (1), ${\rm gmod}^{-,\hspace{.6pt} b\hspace{-2.5pt}}\La$ is Hom-finite Krull-Schmidt and extension-closed in $\GrLa$. Finally, let $M\in {\rm GMod}\La$ be finitely generated and finitely cogenerated. By Propositions \ref{fg-abel} and \ref{fcg-iff}, $M$ is bounded. Since the modules in ${\rm gproj}\La$ are locally dimensional, so are those in ${\rm gmod}^{+,\hspace{.6pt} b\hspace{-2.5pt}}\La$. As a consequence, $M\in {\rm gmod}^{\hspace{.4pt}b}\hspace{-3pt}\La$.  The proof of the lemma is completed.



\vspace{2pt}

In view of Lemma \ref{fg-KrS}(2), ${\rm gmod}^{+,\hspace{.6pt} b\hspace{-2.5pt}}\La$ and 
${\rm gmod}^{-,\hspace{.6pt} b\hspace{-2.5pt}}\La$ are exact $k$-categories, which are not abelian in general. A module $M\in \GrLa$ is called {\it noetherian} if every graded submodule of $M$ is finitely generated. Note that this is equivalent to $M$ being noetherian as a ungraded $\La$-module; see \cite[(5.4.7)]{NCF2}. 

\begin{Prop}\label{fg-abel}

Let $\La=kQ/R$ be a graded algebra with $Q$ a locally finite quiver. Then ${\rm gmod}^{+,\hspace{.6pt} b\hspace{-2.5pt}}\La \hspace{1pt}$ or ${\rm gmod}^{-,\hspace{.6pt} b\hspace{-2.5pt}}\La\hspace{.5pt}$ is abelian if and only if $\La$ is locally left or right noetherian$\hspace{.4pt};$ and in this case, ${\rm gmod}^{+,\hspace{.6pt} b\hspace{-2.5pt}}\La$ or ${\rm gmod}^{-,\hspace{.6pt} b\hspace{-2.5pt}}\La$ is Ext-finite, respectively.

\end{Prop}

\noindent{\it Proof.} Suppose that ${\rm gmod}^{+,\hspace{.6pt} b\hspace{-2.5pt}}\La$ is abelian. Since ${\rm gmod}^{+,\hspace{.6pt} b\hspace{-2.5pt}}\La$ is closed under graded quotients, it is closed under graded submodules. In particular, $\La$ is locally left noetherian.
Conversely, suppose that $\La$ is locally left noetherian. Then, the modules in ${\rm gproj}\La$ are noetherian, and so are those in ${\rm gmod}^{+,\hspace{.6pt} b\hspace{-2.5pt}}\La$. Therefore, ${\rm gmod}^{+,\hspace{.6pt} b\hspace{-2.5pt}}\La$ is abelian. As a consequence, every module in ${\rm gmod}^{+,\hspace{.6pt} b\hspace{-2.5pt}}\La$ admits a graded projective resolution over ${\rm gproj}\La$. In view of Lemma \ref{f-GHom}, we see that ${\rm gmod}^{+,\hspace{.6pt} b\hspace{-2.5pt}}\La$ is Ext-finite. 

Finally, $\La$ is locally right noetherian if and only if $\La^{\rm o}$ is locally left noetherian. In view of the duality $\mk D: {\rm gmod}^{+,\hspace{.6pt} b\hspace{-2.5pt}}\La^{\rm o} \to {\rm gmod}^{-,\hspace{.6pt} b\hspace{-2.5pt}}\La$; see (\ref{fg-KrS}), we see that the second part of the statement holds. The proof of the proposition is completed. 





\vspace{3pt}

Next, we shall study the full subcategories ${\rm gmod}^{+,\hspace{.6pt} p\hspace{-2.5pt}}\La$ and ${\rm gmod}^{-,\hspace{.6pt} i\hspace{-2.5pt}}\La$ of $\GrLa$ of finitely presented modules and of finitely copresented modules, respectively.

\begin{Lemma}\label{fpres}

Let $\La=kQ/R$ be a graded algebra with $Q$ a locally finite quiver. 

\begin{enumerate}[$(1)$]

\vspace{-6pt}

\item The duality $\mk D: {\rm gmod}\La\to {\rm gmod}\La^{\rm o}$ restricts to two mutually quasi-inverse functors 
$\mk D: {\rm gmod}^{+,\hspace{.5pt}p\hspace{-2.5pt}}\La \to {\rm gmod}^{-,\hspace{.5pt}i\hspace{-2.5pt}}\La^{\rm o}$ and $\mk D: {\rm gmod}^{-,\hspace{.5pt}i\hspace{-2.5pt}}\La^{\rm o} \to {\rm gmod}^{+,\hspace{.5pt}p\hspace{-2.5pt}}\La$.

\vspace{1pt}

\item Both ${\rm gmod}^{+,\hspace{.6pt} p \hspace{-2.5pt}}\La$ and 
${\rm gmod}^{-,\hspace{.6pt} i\hspace{-2.5pt}}\La$ are Hom-finite Krull-Schmidt extension-closed subcategories of $\GrLa$, whose intersection is ${\rm gmod}^{\hspace{.4pt}b}\hspace{-3pt}\La$.

\end{enumerate}\end{Lemma}

\vspace{-1.5pt}

\noindent{\it Proof.} (1) \vspace{-1pt} Given $M\in {\rm gmod}^{+,\hspace{.6pt} p\hspace{-2.5pt}}\La,$ it admits a graded projective presentation 
$\xymatrixcolsep{18pt}\xymatrix{\hspace{-2pt} P^{-1} \ar[r]& P^{0} \ar[r] &M \ar[r] &0,} \vspace{-4pt}$ where $P^{-1}\!, P^0\in {\rm gproj}\La.$ Applying $\mf{D}$ yields a graded injective copresentation \vspace{-1.5pt} $\xymatrixcolsep{18pt}\xymatrix{\! 0\ar[r] &\mk DM\ar[r]&\mk DP^0\ar[r]&\mk DP^{-1}}\hspace{-3pt}$ with $\mk DP^0\hspace{-2pt}, \hspace{.4pt} \mk DP^{-1} \!\m\in\m {\rm ginj}\La^{\rm o}\!.$ That is, $\mk DM\in {\rm gmod}^{-,\hspace{.6pt} i\hspace{-2.5pt}}\La^{\rm o}$. Dually, if $N\in {\rm gmod}^{-,\hspace{.6pt} i\hspace{-2.5pt}}\La^{\rm o}$, then $\mk D N\in {\rm gmod}^{+,\hspace{.6pt} p\hspace{-2.5pt}}\La.$ In view of Proposition \ref{gr-duality}, we have a duality $\mk D: {\rm gmod}^{+,\hspace{.5pt}p\hspace{-2.5pt}}\La \to {\rm gmod}^{-,\hspace{.5pt}i\hspace{-2.5pt}}\La^{\rm o}$.

(2) By Lemma \ref{fg-KrS}, ${\rm gmod}^{+,\hspace{.6pt} p\hspace{-2.5pt}}\La$ is Hom-finite, and by Proposition 2.1 in \cite{Aus2}, it is extension-closed in $\GrLa$. Assume that $M\in {\rm gmod}^{+,\hspace{.6pt} p\hspace{-2.5pt}}\La$ with $M=M^1\oplus M^2$. By Lemma \ref{mini-fpres}, $M$ admits a graded projective cover $f: P\to M$ with ${\rm Ker}(f)$ finitely generated. Being finitely generated, $M^i$ has a graded projective cover $f^i: P^i\to M^i$ for $i=1, 2$. Then, ${\rm Ker}(f)\cong {\rm Ker}(f^1) \oplus {\rm Ker}(f^2)$. In particular, ${\rm Ker}(f^i)$ is finitely generated, and hence, $M^i\in {\rm gmod}^{+,\hspace{.6pt} p\hspace{-2.5pt}}\La$ for $i=1, 2$. Thus, ${\rm gmod}^{+,\hspace{.6pt} p\hspace{-2.5pt}}\La$ is closed under direct summands. So, ${\rm gmod}^{+,\hspace{.6pt} p\hspace{-2.5pt}}\La$ is Krull-Schmidt. Then, by Statement (1), ${\rm gmod}^{-,\hspace{.6pt} i\hspace{-2.5pt}}\La$ is also Hom-finite Krull-Schmidt and extension-closed in $\GrLa$.

Finally, assume that $M\in {\rm gmod}^{\hspace{.5pt} b\hspace{-2.5pt}}\La.$ Let $t\in \Z$ be such that $M_i= 0$ for $i\ge t$. By Proposition \ref{p-cover}, $M$ has a graded projective cover $f:P\to M$ with $P\in {\rm gproj}\La$. Write $L={\rm Ker}(f)$. Then $L=\oplus_{i\in \Z} L_i$, where $L_i\subseteq P_i$ for all $i\in \Z$, and $L_i=P_i$ for all $i\ge t.$ Given $i>t$, we see that $L_i=P_i =J_{i-t}P_t\subseteq \rad L.$ This implies that ${\rm top}L=\oplus_{i\le t} (L_i+\rad L)/\rad L$. On the other hand, since $P$ is locally finite dimensional, 
$\oplus_{i\le t}P_i$ is finite dimensional, and so is $\oplus_{i\le t}L_i$. Thus, ${\rm top}L$ is finite dimensional. Since $L\in {\rm GMod}^+\hspace{-3.5pt}\La$, by Proposition \ref{fg-iff}(3), $L$ is finitely generated. So, $M\in {\rm gmod}^{+,\hspace{.6pt} p\hspace{-2.5pt}}\La$. Dually, $M\in  {\rm gmod}^{-,\hspace{.6pt} i\hspace{-2.5pt}}\La.$ The proof of the lemma is completed.

\vspace{3pt}

It is evident that the projective objects in ${\rm gmod}^{+,\hspace{.6pt}p\hspace{-2.5pt}}\La$ are the modules in ${\rm gproj}\La;$ and the injective objects in ${\rm gmod}^{-,\hspace{.6pt}i\hspace{-2.5pt}}\La$ are the modules in ${\rm ginj}\La$. 

\begin{Prop}\label{fdi_pres} 

Let $\La=kQ/R$ be a graded algebra with $Q$ a locally finite quiver. The finitely cogenerated injective objects in ${\rm gmod}^{+,\hspace{.6pt}p\hspace{-2.5pt}}\La$ are the finite dimensional modules in ${\rm ginj}\La;$ and the finitely generated projective objects in ${\rm gmod}^{-,\hspace{.6pt}i\hspace{-2.5pt}}\La$ are the finite dimensional modules in ${\rm gproj}\La.$ 

\end{Prop} 

\noindent{\it Proof.} We shall only prove the first part of the statement. Given $M\in \GrLa$ and  $n\in \Z$, we see that $M_{\ge n}=\oplus_{i\ge n}M_i$ is a graded submodule of $M$. Let $L$ be an  injective object in ${\rm gmod}^{+,\hspace{.6pt}p\hspace{-2.5pt}}\La$, which is finitely cogenerated. Being an essential monomorphism in $\GrLa$; see (\ref{inj_env}), the inclusion map $j: \soc L \to L$ is an injective envelope of $\soc L$ in ${\rm gmod}^{+,\hspace{.6pt}p\hspace{-2.5pt}}\La$. On the other hand, by Proposition \ref{i-envelope}, $\soc L$ admits an injective envelope $q: \soc L\to I$ in $\GrLa$ with $I\in {\rm ginj}\La$. Since $L$ is finite dimensional; see (\ref{fpres}), $L=L_{\ge n}$ for some $n\in \Z$. Fix arbitrarily $t<n$. Since 
$\soc I\cong \soc L,$ which is generated in degrees $\ge n$, we see that $\soc (I_{\ge t})=\soc I$. Hence, $q$ co-restricts to a graded essential monomorphism $q_{\ge t}: \soc L\to I_{\ge t}$ in $\GrLa$; see (\ref{inj_env}). Since $I_{\ge t}$ is finite dimensional and $L$ is injective in ${\rm gmod}^{+,\hspace{.6pt}p\hspace{-2.5pt}}\La$, there exists a graded morphism $f: I_{\ge t}\to L$ such that $j= f\circ q_{\ge t}$. Since $q_{\ge t}$ is graded essential, $f$ is a monomorphism. And since $L_t=0$, we have $I_t=0$. This shows that $I$ is finite dimensional. In particular, $I\in {\rm gmod}^{+,\hspace{.6pt}p\hspace{-2.5pt}}\La.$ Therefore, $L\cong I.$ The proof of the proposition is completed. 

\smallskip

Applying Proposition \ref{fdi_pres},  we obtain the following interesting statement. 

\begin{Prop}\label{lobo-abel}

Let $\La=kQ/R$ be a graded algebra with $Q$ a locally finite quiver. The following statements hold.

\vspace{-2pt}

\begin{enumerate}[$(1)$]

\item Every $S_x$ with $x\in Q_0$ has an injective envelope in ${\rm gmod}^{+,\hspace{.6pt}p\hspace{-2.5pt}}\La$ if and only if $\La$ is locally right bounded if and only if ${\rm gmod}^{-,\hspace{.6pt}i\hspace{-2.5pt}}\La={\rm gmod}^{\hspace{.4pt}b}\hspace{-3pt}\La$.

\item Every $S_x$ with $x\in Q_0$ has a projective cover in ${\rm gmod}^{-,\hspace{.6pt}i\hspace{-2.5pt}}\La$ if and only if $\La$ is locally left bounded  if and only if ${\rm gmod}^{+,\hspace{.6pt}p\hspace{-2.5pt}}\La={\rm gmod}^{\hspace{.4pt}b}\hspace{-3pt}\La$. 

\end{enumerate}\end{Prop}

\noindent{\it Proof.} We shall only prove Statement (1). Clearly, $\La$ is locally right bounded if and only if $\La^{\rm o}$ is locally left bounded, or equivalently, ${\rm ginj}\La \subseteq {\rm gmod}^{\hspace{.4pt}b}\hspace{-3pt}\La$. By Lemma \ref{fpres}, the last condition is equivalent to ${\rm gmod}^{-,\hspace{.6pt}i\hspace{-2.5pt}}\La=
{\rm gmod}^{\hspace{.4pt}b}\hspace{-3pt}\La$. In this case, it is evident that $I_x$ is the  injective envelope of $S_x$ in ${\rm gmod}^{+,\hspace{.6pt}p\hspace{-2.5pt}}\La$, for every $x\in Q_0$. 

Next, suppose that $S_x$ has an injective envelope $j_x: S_x\to L_x$ in ${\rm gmod}^{+,\hspace{.6pt}p\hspace{-2.5pt}}\La$, for every $x\in Q_0$. As argued in the proof of Proposition \ref{fdi_pres}, we see that $L_x\cong I_x$. Thus, $I_x$ is finite dimensional. So, $\La$ is locally right bounded. The proof of the proposition is completed.

\section{Graded almost split sequences}

The objective of this section is to study the existence of almost split sequences for graded modules. 
We shall first construct a graded Nakayama functor, which allows us to establish directly a graded Auslander-Reiten formula for finitely presented graded modules and a generalized Auslander-Reiten formula for finitely copresented graded modules. From these formulas we derive two existence theorems for almost split sequences in the category of all graded modules, one for finitely presented graded modules and one for finitely copresented graded modules. Finally, we shall study when the category of finitely presented graded modules and that of finitely copresented graded modules have almost split sequences.

\subsection{\sc Graded transpose} 
In the locally finite dimensional graded case, Marinez-Villa has introduced the transpose of a graded $\La$-module; see \cite[(1.4)]{RMV2}. In our setting, we need to take more caution in the construction. We start with defining a contravariant functor $(-)^t: {\rm GMod}\La\to {\rm GMod}\La^{\rm o}$ as follows. 
Given $M\in \GrLa$, we define $M^t=\oplus_{i\in \mathbb{Z}}(M^t)_i\in \GrLa^{\rm o}$, where $(M^t)_i=\oplus_{x\in Q_0}{\rm GHom}_{\mit\Lambda}(M\tla\!-i\tra, P_x)$, as follows. Given $\varphi\in {\rm GHom}_{\mit\Lambda}(M\tla\!-i\tra, P_x)$ and $u\in e_x \La_j e_y$, 
considering the graded morphisms $\varphi\tla\m -j\tra: M\tla\m-i\m-\m\m\m j\tra\to P_x\tla- j\tra$ and $P[u]: P_x\tla\!-\m j\tra\to P_y$, the right multiplication by $u$, we set 
$u^{\rm o} \cdot \varphi = P[u] \hspace{-1pt} \circ \hspace{-1pt} \varphi\tla\!-\m j\tra \in {\rm GHom}_{\mit\Lambda}(M\tla\!-i\!- \!\m j\tra, P_y).$
In particular, we have $M^t_i(x)={\rm GHom}_{\mathit\Lambda}(M\tla\!-i\tra, P_x)$ for all $(i,x)\in \Z\times Q_0$. 
Given a morphism $f: M\to N$ in $\GrLa$, setting $(f^t)_{i, x}=\GHom_{\mit\Lambda}(f\langle\!-i\tra, P_x)$ for all $(i,x)\in \Z\times Q_0$, we obtain a morphism $f^t: N^t \to M^t$ in $\GrLa^{\rm o}$. 

\begin{Lemma}\label{transpose-1}

Let $\La=kQ/R$ be a graded algebra with $Q$ a locally finite quiver. 

\begin{enumerate}[$(1)$]

\vspace{-1.5pt}

\item If $M\in \GrLa$ and $s\m\in\m \Z$, then $M\sla s\nra^t \cong M^t \sla -\m s\nra.$ 

\vspace{.5pt}

\item If $M\in \GrLa$ and $V\in {\rm mod}\hspace{.5pt}k$, then $(M \m \otimes\m V)^t \cong M^t \otimes\m DV.$

\end{enumerate}

\end{Lemma}

\noindent{\it Proof.} We shall only prove Statement (2). Let $M \!\in\m \GrLa$ and $V \!\in\m {\rm mod} \hspace{.5pt} k$. For any $x\in Q_0$, applying first the adjunction isomorphism; see \cite[(2.4.9)]{NCF2} and then Lemma \ref{Cor 1.2}(1), we obtain \vspace{-2pt}
$$\GHom_{\mathit\Lambda}(M\otimes V, P_x) \cong {\rm Hom}_k(V, \GHom_{\mathit\Lambda}(M, P_x)) \cong \GHom_{\mathit\Lambda}(M, P_x)\otimes DV.\vspace{-3pt}$$ 
Now, in view of the definition of $(-)^t$, we see that $(M \otimes V)^t\cong M^t\otimes DV$. The proof of the lemma is completed.

\smallskip

The following statement is essential for our later investigation.

\begin{Prop}\label{transpose}

Let $\La=kQ/R$ be a graded algebra with $Q$ a locally finite quiver. The contravariant functor $(-)^t: \hspace{-1pt} \GrLa\to \GrLa^{\rm o}$ is left exact and restricts to a duality $(-)^t: {\rm gproj}\La \to {\rm gproj}\La^{\rm o}$ such that $P_a^t\cong P^{\rm o}_a$ for all $a\in Q_0$. 

\end{Prop}

\noindent{\it Proof.} Since the functors $\GHom_{\mit\Lambda}(-, P_x)$ with $x\in Q_0$ are left exact, so is $(-)^t$. 
Fix $a\in Q_0$. Given $(i, x)\in \Z\times Q_0$, by Proposition \ref{rqz-pm}, we have a $k$-linear isomorphism \vspace{-7pt} 
$$f^a_{i, x}:  (P_a^{\rm o})_i(x)=e_x \La_i^{\rm o} e_a \to \GHom_{\mathit\Lambda}(P_a\tla\!-i\tra, P_x)=(P_a)^t_i(x): v^{\rm o}\to P[v].\vspace{-2pt} $$

It is easy to verify that
$f^a=\oplus_{(i,x)\in \Z\times Q_0} f^a_{i,x}: P_a^{\hspace{.3pt} \rm o} \to P_a^t$ is an isomorphism in $\GrLa^{\rm o}$. Similarly, we may construct an isomorphism $g^a: (P^o_a)^t \to (P_a^{\rm o})^{\rm o}=P_a$ in $\GrLa$. This yields an isomorphism $\zeta_a=g^a \circ (f^a)^t: P_a^{\hspace{.3pt} t \hspace{.3pt}  t} \to P_a$ in ${\rm gproj}\La$. 

Fix $u\in e_a \La_s e_b$. We consider the graded morphisms $P[u]: P_a\to P_b \sla s \sra$ and $P[u^{\rm o}]: P_b^{\rm o} \sla - \m s \sra \to P_a^{\rm o}$, the right multiplications by $u$ and $u^{\rm o}$, respectively. Given $v\in e_b\La_{i-s} e_x$, we have $P[uv]\m=\!P[v] \circ \hspace{-.4pt} P[u]$, that is, $$f^a_{\m i,x}\m(\m P[u^{\rm o}]_{i,x}\m (v^{\rm o})\m)\!=\!\GHom(\m P[u] \mla \m - \m i\rangle, \m P_x)(f^b_{\m i-s, x}(v^{\rm o})\m).$$ So, $f^a_{i,x}\circ P[u^{\rm o}]_{i,x}=P[u]^t_{i,x} \circ f^b \sla\m -s \sra_{i, x}$, and hence, $f^a \circ P[u^{\rm o}] = P[u]^t \circ f^b\sla \m - \m s \sra$. Similarly, $P[u] \circ g^a  = g^b\tla s\tra \circ P[u^{\rm o}]^t$. This implies that $P[u] \circ \zeta^a=\zeta^b\sla s \sra\circ P[u]^{\hspace{.3pt} t \hspace{.3pt} t}$. Since every morphism in $\Hom_{\mit\Lambda}(P_a, P_b\sla s \sra)$ is of the form $P[u]$; see (\ref{rqz-pm}), $\zeta^a$ is natural in $P_a$. 
It is easy to see that $\zeta_a$ extends to a natural isomorphism $\zeta_{_P}: P^{\hspace{.3pt} t \hspace{.3pt} t} \to P$ for every module $P\in {\rm gproj}\La$. Thus, ${\rm id}\cong (-)^t \circ (-)^t$. Similarly, ${\rm id}_{{\rm gproj}\tiny \La^{\rm o}}\cong (-)^t \circ (-)^t$.
The proof of the proposition is completed.

\vspace{2pt}

Recall that the exact category ${\rm gmod}^{+,\hspace{.6pt} p\hspace{-2.5pt}} \La$ is Hom-finite and Krull-Schmidt; see (\ref{fpres}). A morphism $f: M\to N$ in ${\rm gmod}^{+,\hspace{.6pt} p\hspace{-2.5pt}}
\La$ is called {\it radical} if it lies in the Jacobson radical of ${\rm gmod}^{+,\hspace{.6pt} p\hspace{-2.5pt}}\La$. The following statement is interesting.

\begin{Lemma}\label{LR-mini} Let $\La=kQ/R$ be a graded algebra with $Q$ a locally finite quiver. Consider $M\in {\rm gmod}^{+,\hspace{.6pt} p\hspace{-2.5pt}}\La$ with a finitely generated graded projective presentation \vspace{-4pt}  $$\xymatrix{P^{-1} \ar[r]^{d^{-1}} &P^{\hspace{.5pt}0} \ar[r]^{d^{\hspace{.5pt}0}} & M \ar[r] &0.}$$ 

\begin{enumerate}[$(1)$] 

\vspace{-6pt} 

\item If $M$ is indecomposable and not graded projective, then $d^{-1}$ is left minimal. 

\item The epimorphism $d^{\hspace{.5pt}0}$ is a graded projective cover if and only if $d^{-1}$ is  radical.

\vspace{1pt}

\end{enumerate}

\end{Lemma}
\noindent{\it Proof.} (1) Suppose that $M$ is indecomposable and not graded projective. Let $j: N\to P^{\hspace{.5pt}0}$ be the kernel of $d^{\hspace{.5pt}0}: P^0\to M.$
Since $N={\rm Im}(d^{-1})$, there exists a graded epimorphism $v: P^{-1}\to N$ such that 
$d^{-1}=j v$. Suppose that  $f d^{-1}=d^{-1}$ for some graded morphism $f:P^{\hspace{.5pt}0}\to P^{\hspace{.5pt}0}$. Since $v$ is an epimorphism, $fj=j.$ This yields a commutative diagram with exact rows \vspace{-3pt}
$$\xymatrix{0\ar[r] & N\ar[r]^j \ar@{=}[d] & P^{0} \ar[r]^{d^{\hspace{.5pt}0}} \ar[d]^f& M\ar[r]\ar[d]^g & 0\\
0\ar[r] & N\ar[r]^j  & P^{\hspace{.5pt}0}\ar[r]^{d^{\hspace{.5pt}0}} & M\ar[r] & 0}$$ in $\GrLa$. Assume that $g$ is not a graded automorphism. Since $M$ is indecomposable and ${\rm GEnd}_{\mit\Lambda}(M)$ is finite dimensional, $g^s=0$ for some $s\ge 1$. Thus, $d^{\hspace{.5pt}0} f^s=g^s d^{\hspace{.5pt}0}=0$. Therefore, $f^s=j h$, for some $h:P_0\to N$. So, $jhj=f^sj=j$, and hence, $hj ={\rm id}_{N}$. As a consequence, $M$ is graded projective, a contradiction. Thus, $g$ is a graded automorphism, and so is $f$. That is, $d^{-1}$ is left minimal.

(2) By Lemma \ref{proj_cov}(2), $d^{\hspace{.5pt}0}$ is a graded projective cover of $M$ if and only if ${\rm Im}(d^{-1})\subseteq \rad P^{\hspace{.5pt}0}.$ Since $P^{\hspace{.5pt}0}$ is graded projective, this is equivalent to $d^{-1}$ being radical. The proof of the lemma is completed. 

\vspace{2pt}

We are ready to define the graded transpose. Let $M\in {\rm gmod}^{+,\,p}\!\La$. By Lemma \ref{mini-fpres}, $M$ admits a minimal graded projective presentation \vspace{-4pt} $$\xymatrix{P^{-1} \ar[r]^{d^{-1}} &P^{0} \ar[r]^{d^{0}} &M \ar[r] &0} \vspace{-3pt}$$ over ${\rm gproj}\La$. Applying the left exact functor $(-)^t: \GrLa\to \GrLa^{\rm op}$; see (\ref{transpose}) yields an exact sequence \vspace{-4pt} $$\xymatrix{0\ar[r] & M^t \ar[r]^{(d^{0})^t} &(P^{0})^t \ar[r]^{(d^{-1})^t} &(P^{-1})^t \ar[r] &{\rm Coker}(d^{-1})^t \ar[r] &0} \vspace{-4pt} $$ in ${\rm gmod}\La^{\rm o}$. Write ${\rm Tr}M={\rm Coker}(d^{-1})^t$, called the {\it graded transpose} of $M$.

\begin{Lemma}\label{Tr0} Let $\La=kQ/R$ be a graded algebra with $Q$ a locally finite quiver. 

\begin{enumerate}[$(1)$]

\item If $M, N\in {\rm gmod}^{+,\hspace{.6pt} p\hspace{-2.5pt}}\La$, then
${\rm Tr}(M\oplus N)\cong {\rm Tr}M \oplus {\rm Tr}N.$

\item If $M\in {\rm gmod}^{+,\hspace{.6pt} p\hspace{-2.5pt}}\La$, then $M$ is graded projective if and only if ${\rm Tr}M=0.$

\end{enumerate}\end{Lemma}

\noindent{\it Proof.} Statement (1) follows from the fact that the functor $(-)^t$ is additive. Let $M\in {\rm gmod}^{+,\hspace{.6pt} p\hspace{-2.5pt}}\La$ with a minimal graded projective presentation 
\vspace{-3pt} $$\xymatrix{P^{-1} \ar[r]^{d^{-1}} &P^{0} \ar[r]^{d^{0}} &M \ar[r] &0} \vspace{-2pt}$$ over ${\rm gproj}\La$. The necessity of Statement (2) is evident. Suppose that ${\rm Tr}M=0$. Then $(d^{-1})^t$ is a retraction, and by 
Proposition \ref{transpose}, $d^{-1}$ is a retraction. In view of Lemma \ref{LR-mini}(2), $d^{-1}=0$. So, $M\cong P^0.$ The proof of the lemma is completed.

\vspace{2pt}

The following statement is well-known in the finite dimensional ungraded setting. Our approach is different and the proof is shorter; compare \cite[(IV.1.7)]{ARS}.

\begin{Prop}\label{tr-pro}
Let $\La=kQ/R$ be a graded algebra with $Q$ a locally finite quiver. Consider $M\in {\rm gmod}^{+,\hspace{.6pt} p\hspace{-2.5pt}}\La$ with a minimal graded projective presentation \vspace{-4pt} $$\xymatrix{P^{-1}\ar[r]^-{d^{-1}} &P^0\ar[r]^-{d^0} &M\ar[r] &0} \vspace{-2pt}$$ over ${\rm gproj}\La$.
If $M$ is indecomposable and not graded projective, then ${\rm Tr}M$ is indecomposable and not graded projective with a minimal graded projective presentation \vspace{-10pt} $$\xymatrix{(P^{0})^t \ar[r]^-{(d^{-1})^t} &(P^{-1})^t \ar[r]^c&{\rm Tr}M \ar[r] &0.} \vspace{-1pt}$$
\end{Prop}

\noindent{\it Proof.} Suppose that $M$ is indecomposable and not graded projective. By Lemma \ref{LR-mini}, $d^{-1}$ is left minimal and radical, and by Proposition \ref{transpose}, $(d^{-1})^t$ is right minimal and radical. Then, $c$ is a graded projective cover of ${\rm Tr}M$ by Lemma \ref{LR-mini}(2). So, the graded projective presentation of ${\rm Tr}M$ stated in the proposition is minimal. Now, applying the duality $(-)^t: {\rm gproj}\La^{\rm o} \to {\rm gproj}\La$ yields a commutative diagram  with exact rows and vertical isomorphisms \vspace{-3pt} $$\xymatrixcolsep{35pt}\xymatrix{P^{-1} \ar[r]^{d^{-1}}\ar[d]_{\cong} &P^{0} \ar[r]^{d^{0}}\ar[d]^{\cong}  &M \ar@{.>}[d]^{\cong} \ar[r] &0 \\
(P^{-1})^{tt} \ar[r]^{(d^{-1})^{tt}} &(P^{0})^{tt} \ar[r]^{(d^{0})^{tt}} &{\rm Tr}^{\hspace{.5pt}2\hspace{-1.5pt}} M \ar[r] &0.}$$ 

Assume that ${\rm Tr}M=X^1\oplus X^2$, where $X^1, X^2$ are non-zero. If $X^1$ or $X^2$ is graded projective, then the co-restriction of $(d^{-1})^t$ to a non-zero direct summand of $(P^{-1})^t$ is zero. So $(d^{-1})^t$ is not left minimal; see (\ref{minimal_mor}), and hence, $d^{-1}$ is not right minimal, a contradiction. Thus, $X^i$ is not graded projective, for $i=1, 2.$ By Lemma \ref{Tr0}, $M\cong {\rm Tr}^{\hspace{.5pt}2\hspace{-1.5pt}}M\cong {\rm Tr}X^1 \oplus {\rm Tr}X^2$ with ${\rm Tr}X^i\ne 0$ for $i=1,2$, a contradiction. The proof of the proposition is completed. 

\subsection{\sc The graded Nakayama functor.} Composing the contravariant functors $(-)^t$ and $\mf{D}$, we obtain two covariant functors $\nu=\mf{D} \circ (-)^t: \GrLa\to \GrLa$ and $\nu^-= (-)^t \circ \mf{D}: \GrLa\to \GrLa$. By Propositions \ref{gr-duality} and \ref{transpose}, they restrict to two functors $\nu \m:\m {\rm gproj}\La \to {\rm ginj}\La$ and $\nu^- \m: {\rm ginj}\La \m\to\m {\rm gproj}\La$ respectively.

\begin{Theo}\label{N-Functor}

Let $\La=kQ/R$ be a graded algebra with $Q$ a locally finite quiver.

\begin{enumerate}[$(1)$]

\vspace{-0pt}

\item The functors $\nu : {\rm gproj}\La \to {\rm ginj}\La$ \vspace{1pt} and $\nu^- \!: {\rm ginj}\La \m\to\m {\rm gproj}\La$ are mutually quasi-inverse such that $\nu(P_a\sla s\sra \otimes V) \cong  I_a\sla s \sra \otimes V,$ for $(s,a)\in \Z\times Q_0$ and $V \!\m \in \! {\rm mod}k\hspace{.2pt}$.

\vspace{.5pt}

\item Given $M\in \GrLa$ and $P\in {\rm gproj}\La$, we have a binatural $k$-linear isomorphism 
\vspace{-8pt} $$\Phi_{P,M}: \GHom_{\mathit\Lambda}(M, \nu P) \to D\GHom_{\mathit\Lambda}(\m P, M).$$


\end{enumerate}

\end{Theo}

\noindent{\it Proof.} (1) By Propositions \ref{gr-duality} and \ref{transpose}, 
the functors $\nu: {\rm gproj}\La\to {\rm ginj}\La$ and $\nu^- \m: {\rm ginj}\La \m\to\m {\rm gproj}\La$ are mutually quasi-inverse. Given $V\in {\rm mod}k$, by Lemmas \ref{transpose-1} and \ref{mkD}(2) and Proposition \ref{transpose}, $\nu (P_a\sla s\sra\otimes V) \cong \mf{D}(P_a^{\hspace{.4pt} \rm o} \sla -s\sra \otimes DV) \cong I_a\sla s\sra \otimes V.$ 

\vspace{.5pt}

(2) Consider $P_a\sla s\sra$ with $(s,a)\in \Z\times Q_0$ and $M \in \GrLa$. By Proposition \ref{transpose}, we obtain na\-tural graded isomorphisms $f^a\m\sla -s \sra: P^{\rm \hspace{.5pt} o}_a\sla -s \sra\to P^t_a\sla -s \sra$ in ${\rm gproj}\La^{\rm o}$ and $\mf{D}(f^a)\m\sla s \sra: (\nu P_a)\sla s \sra \to I_a\sla s \sra$ in ${\rm ginj}\La$. So, we have a $k$-linear isomorphism $$\rho^{\,s, a}_{\m_M}= \GHom_{\mathit\Lambda}(M, \mf{D}(f^a)\sla s \sra): \GHom_{\mathit\Lambda}(M, (\nu P_a)\sla s \sra) \to \GHom_{\mathit\Lambda}(M, I_a\sla s \sra),$$ which is natural in $M$ and $P_a\sla s \sra$. Next, by Proposition \ref{proj-mor}, we obtain a $k$-linear isomorphism $\eta^{s,a}_{\m_M}: \GHom_{\mathit\Lambda}(P_a\sla s \sra, M) \to M_{-s}(a)$, which is clearly natural in $M$ and $P_a\sla s \sra$. This yields a binatural $k$-isomorphism $$D(\eta^{s,a}_{\m_M}): DM_{-s}(a)\to D\GHom_{\mathit\Lambda}(P_a\sla s \sra, M).\vspace{-1pt}$$

Finally, we have a $k$-linear isomorphism $\theta_a: D(e_a \La^{\rm o}_0e_a)\to k: h\mapsto h(e_a)$. Applying Proposition \ref{proj-inj} for the case $V=k$, we get a $k$-linear isomorphism \vspace{-1pt} $$\psi^{s,a}_{\m_M}: \GHom_{\mathit\Lambda}(M, I_a\sla s \sra) \to D(M_{-s}(a)): g\mapsto \theta_a \circ g_{-s, a},\vspace{-2pt} $$ which is clearly natural in $M$. 
Fix $u\in e_a\La_{t-s}e_b$. Consider $P[u]: P_a\sla s \sra\to P_b\tla t \tra$ and $P[u^{\rm o}]: P^{\hspace{.4pt}\rm o}_b\tla -\m t\tra \to P^{\hspace{.4pt}\rm o}_a\sla -\m s\sra$, the right multiplications by $u$ and by $u^{\rm o}$ respectively. Setting $I[u]=\mk{D}(P[u^{\rm o}])$, we claim that \vspace{-4pt}
$$\xymatrixrowsep{25pt}\xymatrix{
\GHom_{\mathit\Lambda}(M, I_a\sla s \sra\hspace{-0.8pt}) \ar[r]^-{\psi^{s,a}_{\!_M}} \ar[d]_{\GHom_{\mathit\Lambda}(M\m, \hspace{.4pt} I[u])} & D(M_{-s}(a)) \ar[d]^{D(M(u))} \\
\GHom_{\mathit\Lambda}(M, I_b\sla t \sra\hspace{-0.8pt}) \ar[r]^-{\psi^{t,b}_{\!_M}} \!&\! D(M_{-t}(b)\hspace{-0.8pt})
} \vspace{-2pt} $$ commutes. Indeed, given $g\in \GHom_{\mathit\Lambda}(M, I_a\sla s \sra\hspace{-0.8pt})$  
and $m\in M_{-t}(b)$, it is a routine verification that \vspace{.5pt} 
$(\theta_a \circ g_{-s, a} \circ M(u))(m) = g_{-t,b}(m) (u^{\rm o}) = (\theta_b \circ I[u]_{-t, b}\circ g_{-t, b})(m).\vspace{.5pt}$ 
Since every morphism in $\Hom_{\mit\Lambda}(P^{\hspace{.4pt}\rm o}_b\tla -\m t\tra, P^{\hspace{.4pt}\rm o}_a\sla -\m s\sra)$ is of the form $P[u^{\rm o}]$; see (\ref{rqz-pm}),  $\psi^{s,a}_{\!_M}$ is natural in $P_a\sla s \sra$. Thus, we obtain a binatural $k$-linear isomorphism \vspace{-1pt}
$$\Phi_{\hspace{-1pt}P{\m_a}\sla s \sra,\, M}=D(\eta^{s, a}_{\m_M}) \circ \psi^{s, a}_{\m_M} \circ \rho^{s, a}_{\m_M}: \GHom_{\mathit\Lambda}(M, \nu P_a\sla s \sra) \to D\GHom_{\mathit\Lambda}(\m P_a\sla s \sra, M).\vspace{-1pt}$$
It is easy to see that $\Phi_{\hspace{-1pt}P{\m_a}\sla s \sra,\, M}$ extends to a binatural $k$-linear isomorphism $\Phi_{\hspace{-1pt}P, M}$ for ever $P\in {\rm gproj}\La$. The proof of the theorem is completed.

\vspace{2pt}

\noindent {\sc Remark.} By Theorem \ref{N-Functor},
the functor $\nu: {\rm gproj}\La\to \GrLa$ is a Nakayama functor as defined in \cite[(5.4)]{LiN}.

\subsection {\sc Graded Auslander-Reiten translations} 
Given $M\in {\rm gmod}^{+,\hspace{.6pt} p\hspace{-2.5pt}}\La$ and $N\in {\rm gmod}^{-,\hspace{.6pt} i\hspace{-2.5pt}}\La$, we put $\tau M= \mf{D} {\rm Tr} M$ and $\tau^-N = {\rm Tr} \hspace{.8pt} \mf{D}N,$ called the {\it right} and the {\it left Auslander-Reiten translate} of $M$ and $N$, respectively. 

\begin{Lemma}\label{AR-translation}

Let $\La=kQ/R$ be a graded algebra with $Q$ a locally finite quiver. 

\begin{enumerate}[$(1)$]

\item If $M \!\in \m {\rm gmod}^{+,\hspace{.6pt} p\hspace{-2.5pt}}\La$ is indecomposable not graded projective, then $\tau M \!\in \m {\rm gmod}^{-,\hspace{.6pt} i\hspace{-2.5pt}}\La$ is indecomposable  not graded injective such that $\tau^-(\tau M)\cong M$.

\item If $N \!\in \m {\rm gmod}^{-,\hspace{.6pt} i\hspace{-2.5pt}}\La$ is indecomposable not graded injective, then $\tau^-\! N \!\in \m {\rm gmod}^{+,\hspace{.6pt} p\hspace{-2.5pt}}\La$ is indecomposable not graded projective such that $\tau(\tau^-\! N)\cong N.$
    
\end{enumerate}

\end{Lemma}

\noindent{\it Proof.} We shall only prove Statement (2). Let $N\in {\rm gmod}^{-,\hspace{.6pt} i\hspace{-2.5pt}}\La$ be indecomposable and not graded injective with a minimal graded injective copresentation \vspace{-3pt}
$$\xymatrix{0\ar[r]& N\ar[r]^{d^0} & I^{0}\ar[r]^-{d^1} &I^1}$$ over ${\rm ginj}\La.$ Applying the duality $\mk D: {\rm gmod}\La\to {\rm gmod}\La^{\rm o}$ yields a minimal graded projective presentation \vspace{-6pt} 
$$\xymatrix{\mk D I^1\ar[r]^-{\mk D d^{1}} & \mk D I^0\ar[r]^-{\mk D d^0} & \mk D N\ar[r] &0}
\vspace{-1pt} $$ over
${\rm gproj}\La^{\rm o},$ where $\mk D N$ is indecomposable and not graded projective. By Proposition \ref{tr-pro}, ${\rm Tr} \hspace{.5pt} \mk D N$ is indecomposable and not graded projective
with a minimal graded projective presentation \vspace{-4pt} $$\xymatrix{\nu^-\!I^1\ar[r]^-{\nu^-\m(d^{1})} & \nu^-\!I^0\ar[r]^-{c} & \tau^-\! N\ar[r] &0}$$ over ${\rm gproj}\La$. So, 
$\tau^-\! N\in {\rm gmod}^{+,\hspace{.6pt} p\hspace{-2.5pt}}\La$. Dually, we deduce from Proposition \ref{tr-pro} and Theorem \ref{N-Functor}(1) a commutative diagram with exact rows and vertical isomorphism
\vspace{-8pt} $$\xymatrixrowsep{22pt}\xymatrixcolsep{32pt}\xymatrix{
0\ar[r] & N\ar[r]^{d^0}\ar@{.>}[d]_\cong & I^{0}\ar[r]^-{d^1} \ar[d]_\cong \hspace{-2pt} &I^1 \ar[d]^\cong \hspace{-5pt}\\
0\ar[r] & \tau(\tau^-N)\ar[r]^d & \nu(\nu^-\!I^{0})\ar[r]^-{\nu(\nu^-(\m d^1))} &\nu(\nu^-\!I^1)}$$ in 
${\rm gmod}^{-,\hspace{.6pt} i\hspace{-2.5pt}}\La$. The proof of the lemma is completed.

\subsection{\sc Graded Auslander-Reiten formulae} The classical approach to establish an Auslander-Reiten formula involves the tensor product and the adjunction isomorphism; see \cite[(I.3.4)]{Aus}, \cite[(VI.5.1)]{CEi}, \cite{HKra} and \cite[(1.6.1)]{RMV2}. 
We shall take a novel approach by using the graded Nakayama functor. The key ingredient is the following exact sequence for a finitely presented graded module $M$, which relates the functors ${\rm GHom}_\mathit{\Lambda}(-, \tau M)$ and $D{\rm GHom}_\mathit{\Lambda}(M, -)$ in a surprising way.

\begin{Lemma}\label{Naka-fcon}

Let $\La=kQ/R$ be a graded algebra \vspace{-1.5pt} with $Q$ a locally finite quiver. Consider a short exact sequence $\xymatrixcolsep{20pt}\xymatrix{0\ar[r] &X\ar[r]^f &Y\ar[r]^g &Z\ar[r] &0}$ in $\GrLa$. Given $M\in {\rm gmod}^{+,\hspace{.6pt} p\hspace{-2.5pt}}\La$, there exists an exact sequence of $k$-linear maps \vspace{-5pt}
$$\hspace{-28pt}\xymatrixcolsep{22pt}\xymatrix{0\ar[r] &\GHom_{\mit\Lambda}(\m Z\!,\tau \m M)\ar[r] &\GHom_{\mit\Lambda}(\m Y\!,\tau \! M)\ar[r]^{f^*} &\GHom_{\mit\Lambda}(\m X,\tau \! M)&}\vspace{-8pt}$$ 
$$\xymatrixcolsep{22pt}\xymatrix{{\hspace{2pt}}\ar[r] &D\GHom_{\mit\Lambda}(M,Z) \ar[r]^{Dg_*} &D\GHom_{\mit\Lambda}(M,Y) \ar[r]&D\GHom_{\mit\Lambda}(M,X) \ar[r] &0,}\vspace{-1pt}$$ where 
$f^*=\GHom_{\mit\Lambda}(f,\tau\! M)$ and $g_*=\GHom_{\mit\Lambda}(M, g)$.

\end{Lemma}

\noindent{\it Proof.} Consider $M\in {\rm gmod}^{+,\hspace{.6pt} p\hspace{-2.5pt}}\La$ with a minimal graded projective presentaion \vspace{-3pt}$$\xymatrix{P^{-1}\ar[r]^-{d^{-1}} &P^0 \ar[r]^-{d^{\hspace{.5pt}0}} & M\ar[r] & 0.}$$ 
In view of Proposition \ref{tr-pro}, we obtain a minimal graded injective copresentation \vspace{-3pt}
$$\xymatrix{0\ar[r] & \tau M \ar[r] & \nu P^{-1} \ar[r]^-{\nu d^{-1}} &\nu P^{0}.}$$ 
Fix $L\in {\rm GMod}\La$. We shall compute the kernel and the cokernel of the morphism $\GHom_{\mit\Lambda}(L, \nu d^{-1})$. Firstly, applying $\GHom_{\mit\Lambda}(L, -)$ to the minimal graded injective co-presentation of $\tau M$, we obtain an exact sequence \vspace{-3pt}$$\xymatrixcolsep{28pt}\xymatrix{\hspace{-3pt} 0\ar[r]\hspace{-1pt} & \hspace{-1pt} \GHom_{\mit\Lambda}(L, \tau M) \hspace{-1pt} \ar[r] & \hspace{-1pt} \GHom_{\mit\Lambda}(L, \nu P^{-1}) \ar[rr]^--{\GHom_{\mit\Lambda}(L,\,\nu d^{-1})} && \GHom_{\mit\Lambda}(L, \nu P^{0}).}\vspace{-3pt}$$ 
Secondly, applying $D\GHom_{\mit\Lambda}(-, L)$ to the minimal graded projective presentation of $M$, we deduce from Theorem \ref{N-Functor}(2) a commutative diagram with exact lower row and vertical isomorphisms \vspace{-3pt}
$$\xymatrixrowsep{18pt}\xymatrixcolsep{19pt}\xymatrix{
 \GHom_{\mit\Lambda}(L,\nu P^{-1}) \ar[d]_{\cong} \ar[rrr]^---{\GHom_{\mit\Lambda}(L,\,\nu d^{-1})}&&&\hspace{1pt}\GHom_{\mit\Lambda}(L,\nu P^0)\ar[d]_{\cong}\\
 D\GHom_{\mit\Lambda}(P^{-1}\!,L) \ar[rrr]^---{D\Hom_{\mit\Lambda}(d^{-1}, L)} \hspace{-1pt} &&&  \hspace{-1pt}D\GHom_{\mit\Lambda}(P^{0}\!,L) \ar[r]^-{D\m(d^{\hspace{.5pt}0})^*} 
 & D\GHom_{\mit\Lambda}(M,L) \ar[r] & 0,}$$ where $(d^{\hspace{.5pt}0})^*=\GHom_{\mit\Lambda}(d^{\hspace{.5pt}0},L)$.
This yields an exact sequence 
$$\xymatrixcolsep{19pt}\xymatrix{
\GHom_{\mit\Lambda}(L,\nu P^{-1}) \ar[rrr]^---{{\rm GHom}_{\mathit\Lambda}(L,\, \nu d^{-1})} &&&\GHom_{\mit\Lambda}(L,\nu P^0) \ar[r] & D\GHom_{\mit\Lambda}(M,L) \ar[r] & 0.}$$

Combining the above two exact sequences for each of $Z, Y$ and $X$, we obtain a commutative diagram with exact rows and exact columns \vspace{2pt} 
$$\xymatrixcolsep{19pt}\xymatrixrowsep{22pt}\xymatrix{&0 \ar[d] &0 \ar[d] &0 \ar[d]\\
0\ar[r] &\GHom_{\mit\Lambda}(Z,\tau M) \ar[r]\ar[d] &\GHom_{\mit\Lambda}(Y,\tau M) \ar[r]^{f^*}\ar[d] &\GHom_{\mit\Lambda}(X,\tau M) \ar[d]\\
0\ar[r] &\GHom_{\mit\Lambda}(Z,\nu P^{-1}) \ar[r]\ar[d]^{\GHom_{\mit\Lambda}(Z, \hspace{.5pt} \nu d^{-1})} &\GHom_{\mit\Lambda}(Y,\nu P^{-1}) \ar[r]\ar[d]^{\GHom_{\mit\Lambda}(Y, \hspace{.5pt} \nu d^{-1})} &\GHom_{\mit\Lambda}(X,\nu P^{-1}) \ar[r]\ar[d]^{\GHom_{\mit\Lambda}(X, \hspace{.5pt} \nu d^{-1})} &0\\
0\ar[r] &\GHom_{\mit\Lambda}(Z,\nu P^{0}) \ar[r]\ar[d] &\GHom_{\mit\Lambda}(Y,\nu P^{0}) \ar[r]\ar[d] &\GHom_{\mit\Lambda}(X,\nu P^{0}) \ar[r]\ar[d] &0\\
&D\GHom_{\mit\Lambda}(M,Z) \ar[r]^{Dg_*}\ar[d] &D\GHom_{\mit\Lambda}(M,Y) \ar[d]\ar[r] &D\GHom_{\mit\Lambda}(M,X) \ar[r]\ar[d]  &0,\\
&0 &0 &0
}$$ where the two middle rows are exact because $\nu P^{-1}$ and $\nu P^0$ are graded injective.
Using the Snake Lemma, we obtain the desired exact sequence stated in the lemma. The proof of the lemma is completed. 

\smallskip 

Similarly, for a finitely copresented graded module $M$, the right exact functor $D{\rm GHom}_\mathit{\Lambda}(-, M)$ and the left exact functor $D^2{\rm GHom}_\mathit{\Lambda}(\tau^-\hspace{-1.5pt}M, -)$ are nicely related.

\begin{Lemma}\label{Naka-fcond} 

Let $\La=kQ/R$ \vspace{-1.5pt} be a graded algebra with $Q$ a locally finite quiver. Consider a short exact sequence $\xymatrixcolsep{20pt}\xymatrix{0\ar[r] &X\ar[r]^f &Y\ar[r]^g &Z\ar[r] &0}$ in $\GrLa$. Given $M\in {\rm gmod}^{-,\hspace{.6pt}i\hspace{-2.5pt}}\La$, we have an exact sequence \vspace{-5pt} 
$$\hspace{-5pt}\xymatrixcolsep{22pt}\xymatrix{0\ar[r] & \! D^2\GHom_{\mit\Lambda}(\tau^-\hspace{-1.5pt} M\hspace{-1pt},\hspace{-1pt} X)\ar[r]^{} & \! D^2\GHom_{\mit\Lambda}(\tau^-\hspace{-1.5pt} M\hspace{-1pt},\hspace{-1pt} Y)\ar[r]^{D^2\hspace{-.5pt}(g_*)} & \! D^2\GHom_{\mit\Lambda}(\tau^-\hspace{-1.5pt} M\hspace{-1pt},\hspace{-1pt} Z) } \vspace{-5pt}$$ 
$$\xymatrixcolsep{22pt}\xymatrix{\ar[r] & D\GHom_{\mit\Lambda}(X,M) \ar[r]^{Df^*} & D\GHom_{\mit\Lambda}(Y,M) \ar[r]
&D\GHom_{\mit\Lambda}(Z,M) \ar[r] &0,}$$
where $f^*=\GHom_{\mit\Lambda}(f, M)$ and $g_*=\GHom_{\mit\Lambda}(\tau^-\hspace{-1pt}M, g)$.

\end{Lemma}

\noindent{\it Proof.} Since ${\rm gmod}^{-,\,i}\!\La$ is Krull-Schmidt, 
we may assume that $M\in {\rm gmod}^{-,\hspace{.6pt}i\hspace{-2.5pt}}\La$ is indecomposable and not graded injective. By Lemma \ref{AR-translation}(2), $\tau^- \! M\in {\rm gmod}^{+,\hspace{.6pt} p\hspace{-2.5pt}}\La$ with $M\cong \tau (\tau^-\! M)$. And by Lemma \ref{Naka-fcon}, we have an exact sequence \vspace{-6pt} 
$$\hspace{-30pt}\xymatrixcolsep{22pt}\xymatrix{0\ar[r] &\GHom_{\mit\Lambda}(Z, M)\ar[r]
&\GHom_{\mit\Lambda}(Y,M)\ar[r]^{f^*} &\GHom_{\mit\Lambda}(X, M) \ar[r] & }\vspace{-6pt} $$
$$\hspace{-5pt}\xymatrixcolsep{22pt}\xymatrix{D\GHom_{\mit\Lambda}(\tau^-\hspace{-1.3pt} M,Z) \ar[r]^{Dg_*} &D\GHom_{\mit\Lambda}(\tau^- \hspace{-1.3pt} M,Y) \ar[r] &D\GHom_{\mit\Lambda}(\tau^- \hspace{-1.3pt} M,X) \ar[r] &0.}\vspace{-2pt}$$ 

Now, applying the exact functor $D=\Hom_k(-, k)$ yields the desired exact sequence stated in the lemma. The proof of the lemma is completed.

\vspace{4pt}

We shall denote by $\underline{\hspace{-.8pt}\GMod \hspace{-.6pt}}\hspace{1pt}\La$ and $\overline{\hspace{-.8pt}\GMod \hspace{-.6pt}}\hspace{1pt}\La$ the quotient categories of $\GrLa$ modulo the ideal $\mathcal P$ of morphisms factoring through graded projective modules and the ideal $\mathcal I$ of those factoring through graded injective modules, respectively. Given $M, N\in \GrLa$, we write $$\underline{\rm \hspace{-.5pt}GHom}_{\mathit\Lambda}(M, N)=\GHom_{\mathit\Lambda}(M, N)/\mathcal{P}(M, N)$$ and $$
\overline{\rm \hspace{-.5pt}GHom}_{\mathit\Lambda}(M, N)=\GHom_{\mathit\Lambda}(M, N)/\mathcal{I}(M, N).\vspace{2pt}$$ Moreover, put $\underline{\rm \hspace{-.5pt}GEnd}_{\mathit\Lambda}(M)=\underline{\rm \hspace{-.5pt}GHom}_{\mathit\Lambda}(M,M)$ and $\overline{\rm \hspace{-.5pt}GEnd}_{\mathit\Lambda}(M)=\overline{\rm \hspace{-.5pt}GHom}_{\mathit\Lambda}(M,M).$ 


\medskip

We are ready to obtain the promised Auslander-Reiten formulae 
as follows.

\begin{Theo}\label{AR-formu}

Let $\La=kQ/R$ be a graded algebra with $Q$ a locally finite quiver. 
\begin{enumerate}[$(1)$]

\item Given $M \!\in\m {\rm gmod}^{+,\hspace{.6pt} p\hspace{-2.5pt}}\La$ and $X \!\in\m {\rm GMod}\La$, we have a natural $k$-linear isomorphism $$D\hspace{1.5pt}{\rm \underline{\! GHom}}_{\mit\Lambda}(M,X)\cong {\rm GExt}_{\mit\Lambda}^1(X,\tau M).$$

\item Given $M \!\in\m {\rm gmod}^{-,\hspace{.6pt} i\hspace{-2.5pt}}\La$ and $X \!\in\m {\rm GMod}\La$, we have a natural $k$-linear isomorphism 
$$\hspace{25pt} D \hspace{.5pt}{\rm \overline{\m GHom}}_{\mit\Lambda}(X,M) \cong D^2\hspace{0.05em}{\rm GExt}^1_{\mit\Lambda}(\tau^{-}\hspace{-1.5pt} M, X).$$ 

\end{enumerate}\end{Theo}

\noindent{\it Proof.} (1) Let $M\in {\rm gmod}^{+,\hspace{.6pt} p\hspace{-2.5pt}}\La$ and $X\in {\rm GMod}\La$. \vspace{-1.5pt} By Proposition \ref{GM-proj}, there exists a short exact sequence $\xymatrixcolsep{20pt}\xymatrix{0\ar[r] & L \ar[r]^-{q} &P\ar[r]^{p} & X \ar[r] &0}$ in $\GMod\La$, where $P$ is graded projective. Applying $\GHom_{\mit\Lambda}(-,\tau M)$ yields an exact sequence  \vspace{-5pt}
$$\xymatrix{0\ar[r] &\GHom_{\mit\Lambda}(X,\tau M) \ar[r]
&\GHom_{\mit\Lambda}(P,\tau M) \ar[r]^-{q^{*}} & \GHom_{\mit\Lambda}(L,\tau M)} \vspace{-6pt}$$
$$\hspace{-159pt}\xymatrix{{}\ar[r] & {\rm GExt}_{\mit\Lambda}^1(X, \, \tau M) \ar[r] \hspace{1pt} &0,}$$ where 
$q^{*}=\GHom_{\mit\Lambda}(q,\tau M)$. Thus, ${\rm Coker}(q^{*})\cong{\rm GExt}_{\mit\Lambda}^1(X,\tau M),$ which is clearly natural in $X$. On the other hand, by Lemma \ref{Naka-fcon}, we have an exact sequence \vspace{-5pt}
$$\hspace{-20pt}\xymatrixcolsep{22pt}\xymatrix{0\ar[r] &\GHom_{\mit\Lambda}(X,\tau M)\ar[r]
&\GHom_{\mit\Lambda}(P,\tau M)\ar[r]^-{q^*} &\GHom_{\mit\Lambda}(L,\tau M)&}\vspace{-8pt}$$
$$\xymatrixcolsep{22pt}\xymatrix{\ar[r]^-\eta & D\GHom_{\mit\Lambda}(M,X) \ar[r]^-{D(p_*)} &D\GHom_{\mit\Lambda}(M, P) \ar[r] & D\GHom_{\mit\Lambda}(M,L) \ar[r] &0,}$$ where $p_{*}=\GHom_{\mit\Lambda}(M, p).$ Therefore, ${\rm Ker}(D(p_*))= {\rm Im}(\eta) \!\cong\! {\rm Coker}(q^{*}).$ It is not hard to see that this isomorphism is natural in $X$. Since $P$ is graded projective, we have an exact sequence \vspace{-4pt} $$\xymatrixcolsep{22pt}\xymatrix{\GHom_{\mit\Lambda}(M, P) \ar[r]^{p_*} & \GHom_{\mit\Lambda}(M, X) \ar[r] &\underline{\rm \hspace{-.5pt}GHom}_{\mit\Lambda}(M, X) \ar[r] &0.} \vspace{-3pt} $$
So, we have an isomorphism $D\underline{\rm GHom}_{\mit\Lambda}(M,X)\cong{\rm Ker}(D(p_*))$, which is clearly natural in $X$. As a consequence, we obtain a natural isomorphism 
$$ D\hspace{0.05em}{\rm \underline{GHom}}_{\mit\Lambda}(M,X)  \! \cong \!{\rm GExt}_{\mit\Lambda}^1(X,\tau M).
\vspace{-2pt}$$ 

(2) Let $M\in {\rm gmod}^{-,\hspace{.6pt} i\hspace{-2.5pt}}\La$ and $X\in {\rm GMod}\La$. \vspace{-1pt} By Proposition \ref{GM-pi}, there exists a short exact sequence $\xymatrixcolsep{20pt}\xymatrix{0\ar[r] & X \ar[r]^{q} &I\ar[r]^-{p} & L \ar[r] &0}\vspace{.5pt}$ in ${\rm GMod}\La$, where $I$ is graded injective. Applying $D^2\GHom_{\mit\Lambda}(\tau^{-\hspace{-1.5pt}} M,-)$, we obtain an exact sequence \vspace{-7pt} 
$$\xymatrix{0\ar[r] & D^2\GHom_{\mit\Lambda}(\tau^{-}\!M, X) \ar[r]
& D^2\GHom_{\mit\Lambda}(\tau^{-}\!M,I) \ar[r]^-{D^2(p_{*})} &D^2\GHom_{\mit\Lambda}(\tau^{-}\!M, L)} \vspace{-8pt}$$
$$\hspace{-192pt}\xymatrix{\ar[r] & D^2{\rm GExt}_{\mit\Lambda}^1(\tau^{-}\!M, X) \ar[r] & 0,}$$ where 
$p_{*}=\GHom_{\mit\Lambda}(\tau^{-}\!M,p)$. \vspace{1pt} Thus, ${\rm Coker}(D^2(p_*))\cong D^2{\rm GExt}_{\mit\Lambda}^1(\tau^{-}\!M, X),$ which is clearly natural in $X$. And by Lemma \ref{Naka-fcond}, we have an exact sequence \vspace{-8pt}
$$\xymatrixcolsep{22pt}\xymatrix{
0\ar[r] & D^2\GHom_{\mit\Lambda}(\tau^{-}\!M,X)\ar[r]^{} &D^2\GHom_{\mit\Lambda}(\tau^{-}\!M,I)\ar[r]^-{D^2(p_*)} &D^2\GHom_{\mit\Lambda}(\tau^{-}\!M,L)&} \vspace{-10pt}$$ 
$$\hspace{-5pt}\xymatrixcolsep{22pt}\xymatrix{\ar[r]^-\eta & D\GHom_{\mit\Lambda}(X,M) \ar[r]^-{D(q^*)} &D\GHom_{\mit\Lambda}(I,M) \ar[r] & D\GHom_{\mit\Lambda}(L,M)\ar[r] &0,}$$ 
where $q^*=\GHom_{\mit\Lambda}(q, M)$. Thus, ${\rm Ker}(D(q^*))={\rm Im}(\eta)\cong {\rm Coker}(D^2(p_{*})),$ which is natural in $X$. Since $I$ is graded injective, we have an exact sequence \vspace{-6pt} $$\xymatrix{\GHom_{\mit\Lambda}(I, M) \ar[r]^{q^*} & \GHom_{\mit\Lambda}(X, M) \ar[r] & 
\overline{\m \rm GHom\hspace{-.5pt}}_{\mit\Lambda}(X, M) \ar[r] &0.}
\vspace{-6pt}$$
This yields a natural isomorphism $D\hspace{0.5pt}{\overline{\m\rm GHom}}_{\mit\Lambda}(X,M) \cong {\rm Ker}(D(q^*)).$ As a consequence, we obtain a natural isomorphism \vspace{-1pt}
$$D\hspace{0.05em}{\rm \overline{GHom}}_{\mit\Lambda}(X,M) \cong D^2{\rm GExt}_{\mit\Lambda}^1(\tau^-\!M, X).\vspace{-1.5pt}$$ The proof of the theorem is completed.

\smallskip

\noindent {\sc Remark.} (1) We call the formula stated in Theorem \ref{AR-formu}(2) the {\it generalized Auslander-Reiten formula}.

\noindent (2) In case $Q$ is finite, Theorem \ref{AR-formu}(1) was established by Martinez-Villa for graded modules $M$ in ${\rm gmod}^{+,\hspace{.6pt} p\hspace{-2.5pt}}\La$ and $X$ in ${\rm gmod}\La$; see \cite[Page 42]{RMV2}.

\subsection{\sc Graded almost split sequences} To the best of our knowledge, there exists no existence theorem for almost split sequences starting with finitely copresented (graded) modules in the category of all (graded) modules.
over a general (graded) algebra. 
Using our generalized Auslander-Reiten formula and the result in \cite[(3.7)]{LiN}; see also \cite[(2.3)]{SLPC}, we are able to fill up this gap in our graded setting.

\begin{Theo}\label{AR-sequence}

Let $\La=kQ/R$ be a graded algebra with $Q$ a locally finite quiver. 

\begin{enumerate}[$(1)$]

\item If $M \!\in \m {\rm gmod}^{+,\hspace{.6pt} p\hspace{-2.5pt}}\La$ is indecomposable and not graded projective, \vspace{-.5pt} then there exists an almost split sequence $\xymatrixcolsep{20pt} \xymatrix{0 \ar[r] \hspace{-1.5pt}  & \tau M \hspace{-1.5pt} \ar[r] \hspace{-1.5pt}  & E \hspace{-1.5pt} \ar[r] &M \hspace{-1.5pt} \ar[r] \hspace{-1.5pt} &0\hspace{-1.5pt}}$ in  
${\rm GMod}\La$.

\vspace{-.5pt}

\item If $M\!\in\m {\rm gmod}^{-,\hspace{.6pt} i\hspace{-2.5pt}}\La$ is indecomposable and not graded injective, \vspace{-.5pt}  then there exists an almost split sequence $\xymatrixcolsep{20pt}\xymatrix{ 0\ar[r] \hspace{-1.5pt} & M \hspace{-1pt} \ar[r] \hspace{-1.5pt} &E \hspace{-.5pt} \ar[r] \hspace{-1.5pt} & \tau^{\hspace{-.5pt}-\hspace{-1.5pt}}M \hspace{-1pt}\ar[r] \hspace{-1.5pt} &0\hspace{-1.8pt}}$ in  ${\rm GMod}\La$.

\end{enumerate}\end{Theo}

\noindent{\it Proof.} We shall only prove Statement (2), since the proof for Statement (1) is similar and shorter.  
Let $M\in {\rm gmod}^{-,\hspace{.6pt} i\hspace{-2.5pt}}\La$ be indecomposable and not graded injective. Then, ${\rm \overline{GEnd}}_{\it\Lambda}(M)\ne 0.$ 
By Lemma \ref{AR-translation}, $\tau^{\hspace{-.5pt}-}\hspace{-1.5pt}M \!\m \in \m {\rm gmod}^{+\m,\hspace{.4pt} p\hspace{-3pt}}\La$ is indecomposable. Thus, $M$ and $\tau^{\hspace{-.5pt}-}\hspace{-1.5pt}M$ are strongly indecomposable; see (\ref{fpres}). 
And by Theorem \ref{AR-formu}(2), \vspace{1pt} we have a functorial isomorphism $\vspace{.5pt}\Phi: D^2{\rm GExt}^1_{\mit\Lambda}(\tau^-\hspace{-1pt} M,-) \to D\hspace{0.05em}{\rm \overline{GHom}}_{\mit\Lambda}(-, M)$. In particular, \vspace{1pt} $D^2{\rm GExt}^1_{\mit\Lambda}(\tau^{-}M, M)\cong D\overline{\rm GEnd}_{\mit\Lambda}(M)$, which is finite dimensional. So, ${\rm GExt}^1_{\mit\Lambda}(\tau^-M,M)\cong D^2{\rm GExt}^1_{\mit\Lambda}(\tau^-M,M)\ne 0$. As a consequence, ${\rm GExt}^1_{\mit\Lambda}(\tau^-M,M)$ has a nonzero socle as a left ${\rm GEnd}_\mathit\Lambda(M)$-module. 

Composing the canonical monomorphism $\vspace{1pt}{\rm GExt}^1_{\mit\Lambda}(\tau^{-}\! M,-) \!\to \! D^2{\rm GExt}^1_{\mit\Lambda}(\tau^{-}\! M,-)$ with $\vspace{.5pt}\Phi: D^2{\rm GExt}^1_{\mit\Lambda}(\tau^-\hspace{-1pt} M,-) \to D\hspace{0.05em}{\rm \overline{GHom}}_{\mit\Lambda}(-, M)$, we obtain a functorial monomorphism 
$\Psi:{\rm GExt}^1_{\mit\Lambda}(\tau^{-}N,-)\to D\hspace{0.05em}{\rm \overline{GHom}}_{\mit\Lambda}(-, M)$. In view of Theorem 2.3(2) in \cite{SLPC}, we have a desired almost split sequence as stated in Statement (2). The proof of the theorem is completed. 

\smallskip

\noindent {\sc Remark.} 
In case $Q$ is finite, an existence theorem in ${\rm gmod}\La$ for almost split sequences ending with finitely presented modules was obtained in \cite[(1.6.1), (1.7.1)]{RMV2}. 


\smallskip

Next, we shall study the existence of almost split sequences in the exact categories ${\rm gmod}^{+,\hspace{.6pt} p\hspace{-2.5pt}}\La$ and ${\rm gmod}^{-,\hspace{.6pt} i\hspace{-2.5pt}}\La$. As shown below, their almost split sequences are also almost split sequences in $\GrLa$; compare \cite[(3.6)]{BLP}. 

\begin{Theo}\label{fd-ass}

Let $\La=kQ/R$ be a graded algebra with $Q$ a locally finite quiver. 

\begin{enumerate}[$(1)$]

\vspace{-2pt}

\item If \vspace{-.5pt} $M\in {\rm gmod}^{+,\hspace{.6pt} p\hspace{-2.5pt}}\La$ is indecomposable not graded injective, then there exists an almost split sequence $\vspace{-.5pt} \xymatrixcolsep{18pt}\xymatrix{0\ar[r]& M \ar[r] & N \ar[r] & L\ar[r] & 0\hspace{-1.5pt}}$ in ${\rm gmod}^{+,\hspace{.6pt} p\hspace{-2.5pt}}\La$ if and only if $M\in {\rm gmod}^b\hspace{-2.5pt}\La;$ and in this case, the sequence is also almost split in $\GrLa$.

\vspace{1pt}

\item If \vspace{-.5pt} $M\in {\rm gmod}^{-,\hspace{.6pt} i\hspace{-2.5pt}}\La$ is indecomposable not graded projective, then there exists an almost split sequence $\vspace{-.5pt} \xymatrixcolsep{18pt}\xymatrix{0\ar[r]& L \ar[r] & N \ar[r] & M\ar[r] & 0\hspace{-1.5pt}}$ in ${\rm gmod}^{-,\hspace{.6pt} i\hspace{-2.5pt}}\La$ if and only if $M\in {\rm gmod}^b\hspace{-2.5pt}\La;$ in this case, the sequence is an almost split sequence in $\GrLa$.

\end{enumerate}

\end{Theo}

\noindent{\it Proof.} We shall only prove Statement (1). Let $M\in {\rm gmod}^{+,\hspace{.6pt} p\hspace{-2.5pt}}\La$ be indecomposable and not graded injective. If \vspace{-1pt} $M\in {\rm gmod}^{\hspace{.5pt}b\hspace{-2.5pt}}\La,$ then $M\in {\rm gmod}^{-,\hspace{.6pt} i\hspace{-2.5pt}}\La$, and by Theorem \ref{AR-sequence}(2), there exists an almost split sequence $\xymatrixcolsep{18pt}\xymatrix{0\ar[r]& M \ar[r] & E \ar[r] & \tau^-\hspace{-1.5pt} M\ar[r] & 0}$ in ${\rm GMod}\La$, where $\tau^-\hspace{-1.5pt} M\in  {\rm gmod}^{+,\hspace{.6pt} p \hspace{-2.5pt}}\La$. Since ${\rm gmod}^{+,\hspace{.6pt} p\hspace{-2.5pt}}\La$ is extension-closed in $\GrLa$; see (\ref{fpres}), this is an almost split sequence in ${\rm gmod}^{+,\hspace{.6pt} p\hspace{-2.5pt}}\La$. 
Conversely, suppose that $\xymatrixcolsep{20pt}\xymatrix{\hspace{-3pt}0\ar[r]& M \ar[r] & N \ar[r] & L \ar[r] & 0\hspace{-1.5pt}}$ is an almost split sequence in ${\rm gmod}^{+,\hspace{.6pt} p\hspace{-2.5pt}}\La$. Then, $L$ is indecomposable and not graded projective. In view of Theorem \ref{AR-sequence}, we can construct a commutative diagram \vspace{-5pt}
$$\xymatrix{0\ar[r]& M \ar[r]^f \ar[d]^u & N \ar[d]^v\ar[r]^g & L\ar[r] \ar[r] \ar@{=}[d] & 0\\ 0\ar[r]& X \ar[r]^q & Y \ar[r]^p& L \ar[r] & 0,} \vspace{-4pt}$$ where the lower row is an almost split sequence in $\GrLa$ with $X\in {\rm gmod}^{-,\hspace{.6pt} i\hspace{-2.5pt}}\La$. Note that there exists an integer $n$ such that $M=M_{\ge n};$ $N=N_{\ge n}$ and $L=L_{\ge n}.$ Since $X_{\ge n}$ is finite dimensional,  ${\rm gmod}^{+,\hspace{.6pt} p\hspace{-2.5pt}}\La$ contains  a commutative diagram 
\vspace{-4pt} $$\xymatrix{0\ar[r]& M \ar[r]^f \ar[d]^{u_{\ge n}} & N \ar[d]^{v_{\ge n}} \ar[r]^g & L\ar[r] \ar[r] \ar@{=}[d] & 0\\ 0\ar[r]& X_{\ge n} \ar[r]^{q_{\ge n}} & Y_{\ge n} \ar[r]^{p_{\ge n}} & L \ar[r] & 0.}
\vspace{-1pt}$$ Since $p$ is not a retraction and $L_i=0$ for all $i<n$, neither is $p_{\ge n}$. Thus, $v_{\ge n}$ is a section, and consequently, $u_{\ge n}$ is a monomorphism. So, $M$ is finite dimensional. As has been shown, $\vspace{-2pt} \xymatrixcolsep{20pt}\xymatrix{0\ar[r]& M \ar[r] & N \ar[r] & L \ar[r] & 0\hspace{-1.5pt}}$ is isomorphic to the almost split sequence $\vspace{-0pt}\xymatrixcolsep{20pt}\xymatrix{0\ar[r]& M \ar[r] & E \ar[r] & \tau^-\hspace{-1.5pt} M\ar[r] & 0}$ in $\GrLa$. The proof of the theorem is completed.


\subsection{\sc Subcategories having almost split sequences} In this subsection, we shall study when ${\rm gmod}^{+\m,\hspace{.6pt} p\hspace{-2.8pt}}\La$, ${\rm gmod}^{-\m,\hspace{.6pt} i\hspace{-2.8pt}}\La$ and ${\rm gmod}^{\hspace{.5pt}b\hspace{-2.5pt}}\La$ have almost split sequences.


\begin{Theo}\label{ass_pi}

Let $\La=kQ/R$ be a graded algebra with $Q$ a locally finite quiver. If $\La$ is locally left $($respectively, right$\hspace{.5pt})$ bounded, then ${\rm gmod}^{+\m,\hspace{.6pt} p\hspace{-2.8pt}}\La$ and ${\rm gmod}^{-\m,\hspace{.6pt} i\hspace{-2.8pt}}\La$ both have almost split sequences on the left $($respectively, right$\hspace{.5pt})$.


\end{Theo}

\noindent{\it Proof.} We only prove the first part of the statement. Let $\La$ be locally left bounded. By Lemma \ref{fg-KrS}(2), ${\rm gmod}^{+,\hspace{.6pt} p\hspace{-2.5pt}}\La = {\rm gmod}^{\hspace{.5pt}b\hspace{-2.5pt}}\La\subseteq {\rm gmod}^{-,\hspace{.6pt} i\hspace{-2.5pt}}\La.$ Consider an indecomposable and not Ext-injective module $M$ in ${\rm gmod}^{+,\hspace{.6pt} p\hspace{-2.5pt}}\La$.  \vspace{-0pt} Then, $M$ is finite dimensional and not graded injective. By Theorem \ref{fd-ass}(1), ${\rm gmod}^{+\m,\hspace{.6pt} p\hspace{-2.8pt}}\La$ has an almost split sequence starting with $M$. So, ${\rm gmod}^{+\m,\hspace{.6pt} p\hspace{-2.8pt}}\La$ has almost split sequences on the left.

Next, let $N\in {\rm gmod}^{-\m,\hspace{.6pt}i\hspace{-2.8pt}}\La$ be indecomposable and not Ext-injective. In particular, $N$ is not graded injective. By Theorem \ref{AR-sequence}(2), \vspace{-1pt} there exists an almost split sequence
$\xymatrixcolsep{22pt}\xymatrix{0\ar[r]& N \ar[r] & E \ar[r] &  \tau^-\hspace{-1.2pt} N \ar[r] & 0\hspace{-2pt} \vspace{-2pt}}$ in ${\rm GMod}\La$, where $\tau^-\hspace{-1.2pt} N$ lies in ${\rm gmod}^{+\m,\hspace{.6pt} p\hspace{-2.8pt}}\La\subseteq {\rm gmod}^{-,\hspace{.6pt} i\hspace{-2.5pt}}\La$. Since ${\rm gmod}^{-,\hspace{.6pt} i\hspace{-2.5pt}}\La$ \vspace{.5pt} is extension-closed in $\GrLa$; see (\ref{fpres}), this is an almost split sequence in ${\rm gmod}^{-,\hspace{.6pt} i\hspace{-2.5pt}}\La$. So, ${\rm gmod}^{-\m,\hspace{.6pt} i\hspace{-2.8pt}}\La$ has almost split sequences on the left. The proof of the theorem is completed.

\smallskip

\noindent{\sc Example.} Consider the locally right bounded graded algebra $\La=kQ/R$, where \vspace{-2pt}
$$\begin{tikzpicture}[-{Stealth[inset=0pt,length=4.5pt,angle'=35,round,bend]},scale=.6]
    \draw  (2.2,0.15) -- (3.7,0.6);
    \draw  (2.2,-0.15) -- (3.7,-0.6);
    \draw  (4.3,0.6) -- (5.8,0.15);
    \draw  (4.3,-0.6) -- (5.8,-0.15);
    \node at (2.8,0.65) {$\alpha$};
    \node at (2.8,-0.75) {$\beta$};
    \node at (1.95,0) {$1$};
    \node at (4,0.66) {$2$}; 
    \node at (4,-0.66) {$3$};
    \node at (5.1, 0.65) {$\gamma$};
    \node at (5.1,-0.75) {$\delta$};
    \node at (6.1,0) {$4$};
    \node at (0.6,0) {$Q:$};
    \draw (6.4,0) -- node[above]{} (7.9,0);
    \draw (8.5,0) -- node[above]{} (10,0);
    \node at (8.2,0) {$5$};
    \node at (10.3,0) {$6$};
    \draw (10.6,0)  -- (12.1,0);
    \node at (12.7,0){$\cdots$};  
 \end{tikzpicture}\vspace{-7pt} $$
and $R=\langle \gamma\alpha-\delta\beta\rangle$. 
By Theorem \ref{ass_pi}, ${\rm gmod}^{+,\hspace{.6pt} p\hspace{-2.5pt}}\La$ has almost split sequences on the right. On the other hand, $S_1$ has a minimal graded projective resolution \vspace{-2pt} $$\xymatrixcolsep{22pt}\xymatrix{0\ar[r] & P_4\oplus P_4\ar[r]
&P_2\oplus P_3\ar[r]
&P_1\ar[r] &S_1\ar[r] &0.} \vspace{-3pt}
$$ 
Thus, $\rad P_1\in {\rm gmod}^{+,\hspace{.5pt}p}\!\La$. Since $\rad^2P_1\cong P_4$ and ${\rm soc}(\rad P_1)=0$, $\rad P_1$ is indecompo\-sable and not graded injective. Being infinite dimensional, $\rad P_1$ is not the starting term of any almost split sequence in ${\rm gmod}^{+,\hspace{.6pt} p\hspace{-2.5pt}}\La$; see (\ref{fd-ass}). Thus, ${\rm gmod}^{+,\hspace{.6pt} p\hspace{-2.5pt}}\La$ does not have almost split sequences on the left. 

 

\vspace{2pt}

We conclude this section with the following statement, which generalizes the result stated in \cite[(3.5)]{GoG}. 

\begin{Theo}\label{rlb-ass}

Let $\La=kQ/R$ be a graded algebra with $Q$ a locally finite quiver. If $\La$ is locally $($left, right$\,)$ bounded, then ${\rm gmod}^{\hspace{.5pt}b\hspace{-2.5pt}}\La$ has almost split sequences $($on the left, on the right$\,)$.

\end{Theo}

\noindent{\it Proof.} If $\La$ is locally left or right bounded, then 
${\rm gmod}^{\hspace{.5pt}b\hspace{-2.5pt}}\La$ coincides with ${\rm gmod}^{+,\hspace{.6pt} p\hspace{-2.5pt}}\La$ or ${\rm gmod}^{+,\hspace{.6pt} p\hspace{-2.5pt}}\La$; see (\ref{lobo-abel}), which has almost split sequences on the left or right respectively; see (\ref{ass_pi}). 
The proof of the theorem is completed.

\section{Graded almost split triangles}

The objective of this section is to study the existence of almost split triangles in various derived categories of graded modules. The graded Nakayama functor constructed in Section 3 allows us to apply the results in \cite[Section 5]{LiN} for this purpose. Some of our results are analogous 
to those of Happel for finite dimensional ungraded algebras in \cite{Ha1, Ha3}. Throughout this section let $\La=kQ/R$ be a graded algebra, where $Q$ is a locally finite quiver and $R$ is a relation ideal in $kQ$. 

\subsection{\sc The existence in general} The bounded complexes of projective objects and those of injective objects in an abelian category play an essential role in the study of almost split triangles in their derived categories; see \cite{Ha1, 
Ha3, LiN}. 
\vspace{-1pt}

\begin{Lemma}\label{proj-htp}

Let $\La=kQ/R$ be a graded algebra with $Q$ a locally finite quiver.

\begin{enumerate}[$(1)$]

\vspace{-2.6pt}

\item The categories $K^b({\rm gproj}\La)$ and $K^b({\rm ginj}\La)$ are Hom-finite and Krull-Schmidt.

\vspace{.5pt}

\item The Nakayama functor induces two mutually quasi-inverse triangle equivalences $\nu: K^b({\rm gproj}\La)\to K^b({\rm ginj}\La)$ and $\nu^-\!: K^b({\rm ginj}\La)\to K^b({\rm gproj}\La).$

\end{enumerate}\end{Lemma}

\noindent{\it Proof.} (1) Since ${\rm gproj}\La$ is Hom-finite and Krull-Schmmidt; see (\ref{proj-KS}), $C^b({\rm gproj}\La)$ is a Hom-finite additive subcategory of $C^b({\rm gmod}\La)$, which is closed under direct summands. Therefore, $C^b({\rm gproj}\La)$ is Krull-Schmidt, and consequently, $K^b({\rm gproj}\La)$ is Hom-finite and Krull-Schmidt; see \cite[page 431]{SLiu}. Similarly, $K^b({\rm ginj}\La)$ is Hom-finite and Krull-Schmidt.

(2) The Nakayama functor $\nu: {\rm gproj}\La\to \GrLa$ induces two mutually quasi-inverse equivalences $\nu: {\rm gproj}\La\to {\rm ginj}\La$ and $\nu^-: {\rm ginj}\La\to {\rm gproj}\La$; see (\ref{N-Functor}). Applying them component-wise, we obtain two mutually quasi-inverse triangle equivalences $\nu: K^b({\rm gproj}\La)\to K^b({\rm ginj}\La)$ and $\nu^-\!: K^b({\rm ginj}\La)\to K^b({\rm gproj}\La).$ The proof of the lemma is completed.


\vspace{2pt}

It is well-known that $\!K^b({\rm gproj}\La)$ and $\!K^b({\rm ginj}\La)$ \vspace{.5pt} are full triangulated subcate\-gories of $\!D^b({\rm gmod} \La)$ and $\!D^b({\rm GMod} \La)$; see \cite[(10.4.7)]{Wei}, while $\!D^b({\rm gmod} \La)$ and $\!D^b({\rm GMod} \La)$ are full triangulated subcategories of $\!D({\rm gmod}\La)$ and $\!D({\rm GMod}\La)$ respectively; see \cite[(III.3.4.5)]{Mil}. Note, however, that $D({\rm gmod}\La)$ is not necessarily a triangulated subcategory of $\!D({\rm GMod}\La)$. 

\vspace{-1pt}

\begin{Theo}\label{GAST-1}

Let \vspace{.5pt} $\La=kQ/R$ be a graded algebra with $Q$ a locally finite quiver. If $P^\pdt \!\in \m K^b\m({\rm gproj}\La)$ is indecomposable, \vspace{-1pt} then \vspace{-.5pt} $D^{\hspace{.5pt}b\hspace{-0.6pt}}({\rm gmod}\La),\!$ $D({\rm gmod}\La),\!$ $D^{\hspace{.5pt}b\hspace{-0.8pt}}(\m{\rm GMod}\La)$ and $D(\m{\rm GMod}\La)$ each have an almost split triangle 
$\xymatrixcolsep{18pt}\xymatrix{\hspace{-2pt}\nu\m P^\pdt\hspace{.3pt} [-\m1]  \hspace{-1pt}\ar[r]  & \hspace{-1pt}L^\cdt \hspace{-2pt} \ar[r] &  P^\pdt \hspace{-2pt} \ar[r] & \nu P^\pdt\hspace{-1pt}.} \vspace{-4pt}$

\end{Theo}

\noindent{\it Proof.} Let $P^\pdt\in K^b\m({\rm gproj}\La)$ be indecomposable. By Lemma \ref{proj-htp}, $P^\pdt$ and $\nu P^\pdt$ are strongly indecomposable. Considering the Nakayama functors $\nu\hspace{-.5pt}:\hspace{-1pt} {\rm gproj}\La \!\to\! {\rm gmod}\La$ and $\nu\hspace{-.5pt}:\hspace{-1pt} {\rm gproj} \La\!\to\! \GrLa$; see (\ref{N-Functor}), we deduce from Theorem 5.8 in \cite{LiN} a desired almost split triangle in each of $D^{\hspace{.5pt}b\hspace{-0.8pt}}({\rm gmod} \La)$, $\!D({\rm gmod}\La)$, $D^{\hspace{.5pt}b\hspace{-0.8pt}}({\rm GMod} \La)$ and $\!D({\rm GMod}\La)$. The proof of the theorem is completed.



\subsection{\sc The existence in the locally noetherian case} The bounded derived cate\-gory of finitely generated graded modules over noetherian graded algebras is important in geometry; see, for example, 
\cite[(2.12.6)]{BGS}.

\begin{Lemma}\label{FGDC}

Let $\La=kQ/R$ be a graded algebra with $Q$ a locally finite quiver. If $\La$ is locally left or right noetherian, then $D^{\hspace{.5pt}b\hspace{-0.8pt}}({\rm gmod}^{+\m,\hspace{.4pt}b\hspace{-3pt}}\La)$ or $D^{\hspace{.5pt}b\hspace{-0.8pt}}({\rm gmod}^{-\m,\hspace{.4pt}b\hspace{-3pt}}\La)$ is a Hom-finite Krull-Schmidt full triangulated subcategory of $D^{\hspace{.5pt}b\hspace{-0.8pt}}({\rm gmod}\La)$, respectively. 

\end{Lemma}

\noindent{\it Proof.} Assume that $\La$ is locally left noetherian. Then,  ${\rm gmod}^{+\m,\hspace{.4pt}b\hspace{-3pt}}\La$ is an Ext-finite abelian subcategory of ${\rm gmod}\La$; see (\ref{fg-abel}). By Corollary B in \cite{LeC}, $D^{\hspace{.5pt}b\hspace{-0.8pt}}({\rm gmod}^{+\m,\hspace{.4pt}b\hspace{-3pt}}\La)$ is Hom-finite and Krull-Schmidt. Since ${\rm gmod}^{+\m,\hspace{.4pt}b\hspace{-3pt}}\La$ has enough graded projective modules; see (\ref{p-cover}), $D^{\hspace{.5pt}b\hspace{-0.8pt}}({\rm gmod}^{+\m,\hspace{.4pt}b\hspace{-3pt}}\La)$ is 
a full triangulated subcategory of $D^{\hspace{.5pt}b\hspace{-0.8pt}}({\rm gmod}\La)$; see (\ref{Der_add}).  
The proof of the lemma is completed. 

\smallskip

In case $\La$ is locally left and right noetherian, \vspace{.5pt} 
both $D^{\hspace{.5pt}b\hspace{-0.6pt}}({\rm gmod}^{\m+\m,\hspace{.4pt}b\hspace{-3pt}}\La)$ and $D^{\hspace{.5pt}b\hspace{-0.6pt}}({\rm gmod}^{\m-\m,\hspace{.4pt}b\hspace{-3pt}}\La)$ are full triangulated subcategories of $D^{\hspace{.5pt}b\hspace{-0.6pt}}({\rm gmod}\La)$.

\begin{Theo}\label{GAST-2}

Let $\La=kQ/R$ be a locally left and right noetherian graded algebra, where $Q$ is a locally finite quiver.

\begin{enumerate}[$(1)$]

\vspace{-1pt}

\item If \vspace{-1.5pt} $M^\cdt \!\in \! D^{\hspace{.5pt}b\hspace{-0.6pt}}({\rm gmod}^{\m+\m,\hspace{.4pt}b\hspace{-3pt}}\La)$ is indecomposable, then $D^{\hspace{.5pt}b\hspace{-0.6pt}}({\rm gmod}\La)$ has an almost split triangle $\xymatrixcolsep{22pt}\xymatrix{\hspace{-2pt} N^\cdt \!\m \ar[r] & L^\cdt \! \ar[r] & M^\pdt \! \ar[r] & N^\pdt[1]\hspace{-2pt}}\vspace{-3pt}$ if and only if $M^\cdt$ admits a finite graded projective resolution over ${\rm gproj}\La\hspace{.5pt};$ and in this case, $N^\cdt\in D^{\hspace{.5pt}b\hspace{-0.6pt}}({\rm gmod}^{\m-\m,\hspace{.4pt}b\hspace{-3pt}}\La)$.

\vspace{1pt}

\item If $M^\cdt \! \m \in \!\m D^{\hspace{.5pt}b\hspace{-0.6pt}}({\rm gmod}^{\m-\m,\hspace{.4pt}b\hspace{-3pt}}\La)$ is indecomposable, \vspace{-2pt} then $D^{\hspace{.5pt}b\hspace{-0.6pt}}({\rm gmod}\La)$ has an almost split triangle \hspace{-3pt}$\xymatrixcolsep{20pt}\xymatrix{M^\cdt \!\m \ar[r] & L^\cdt \! \ar[r] & N^\pdt \! \ar[r] & M^\pdt[1]}$ \vspace{-3pt} \hspace{-8pt} if and only if $M^\cdt$ admits a finite graded injective coresolution over ${\rm ginj}\La\hspace{.5pt};$ and in this case, $N^\cdt\in D^{\hspace{.5pt}b\hspace{-0.6pt}}({\rm gmod}^{\m+\m,\hspace{.4pt}b\hspace{-3pt}}\La)$. 

\end{enumerate} \end{Theo}

\noindent{\it Proof.} We shall only prove Statement (1). Consider an indecomposable complex $M^\cdt$ in $D^{\hspace{.5pt}b\hspace{-0.8pt}}({\rm gmod}^{+\m,\hspace{.4pt}b\hspace{-3pt}}\La).$ If $M^\cdt$ has a graded projective resolution $P^\pdt$ in $C^{\hspace{.5pt}b\hspace{-0.6pt}}({\rm gproj}\La)$, then $M^\cdt \m \cong \m P^\pdt$ in $D^{\hspace{.5pt}b\hspace{-0.6pt}}({\rm gmod}\La)$. By Theorem \ref{GAST-1}, \vspace{-1pt} $D^{\hspace{.5pt}b\hspace{-0.6pt}}({\rm gmod}\La)$ has an almost split triangle $\xymatrixcolsep{22pt}\xymatrix{\hspace{-3pt} \nu P^\pdt[-\m1] \ar[r] & L^\cdt \ar[r] & M^\cdt\ar[r] & \nu P^\pdt,}\vspace{-4pt}$ where $\nu P^\pdt[-\m1]  \in  D^{\hspace{.5pt}b\hspace{-0.6pt}}({\rm gmod}^{\m-\m,\hspace{.4pt}b\hspace{-3pt}}\La)$.

Conversely, \vspace{-1pt} suppose that $\xymatrixcolsep{22pt}\xymatrix{\!N^\cdt\ar[r] & L^\cdt \ar[r] & M^\cdt\ar[r] & N^\cdt[1]\!}\vspace{-3pt}$ is an almost split triangle in $D^{\hspace{.5pt}b\hspace{-0.8pt}}({\rm gmod}\La)$. Since ${\rm gmod}^{+\m,\hspace{.4pt}b\hspace{-2.8pt}}\La$ is abelian with enough graded projective mo\-dules, $M^\cdt$ has a graded projective resolution in $C^-({\rm gproj}\La)$; see \cite[(7.5)]{Gri}. Now, it follows from  Theorem 5.2 in \cite{LiN} that $M^\cdt$ has a 
graded projective resolution in $C^{\hspace{.5pt}b\hspace{-0.6pt}}({\rm gproj}\La)$. The proof of the theorem is completed.

\smallskip

\noindent{\sc Example.} Theorem \ref{GAST-2} holds for graded special multi-serial algebras. 




\subsection{\sc Existence in the locally bounded case} In this subsection, we shall concentrate on the bounded derived category of finite dimensional graded $\La$-modules. Although our results are analogous to those of Happel in \cite{Ha1, Ha3}, they do exhibit some particular features of the locally bounded graded setting. 

\begin{Prop}\label{GAST-lb}

Let $\La=kQ/R$ be a locally bounded graded algebra, where $Q$ is a locally finite quiver. 
If $M^\cdt\in D^{\hspace{.4pt}b\hspace{-.3pt}}({\rm gmod}^{\hspace{.4pt}b\hspace{-3pt}}\La)$ is indecomposable, then 

\begin{enumerate}[$(1)$]

\vspace{-3pt}

\item $D^{\hspace{.4pt}b\hspace{-.3pt}}({\rm gmod}^{\hspace{.4pt}b\hspace{-3pt}}\La)$ has an almost split triangle $\xymatrixcolsep{20pt}\xymatrix{\! N^\cdt \hspace{-1pt} \ar[r] & L^\cdt \hspace{-1pt}  \ar[r] & M^\cdt \hspace{-1pt} \ar[r] & N^\cdt[1]\!}$ \vspace{-1pt} if and only if $M^\cdt$ has a finite graded projective resolution over ${\rm gproj}\La\hspace{.5pt};$ 

\item $D^{\hspace{.4pt}b\hspace{-.3pt}}({\rm gmod}^{\hspace{.4pt}b\hspace{-3pt}}\La)$ has an almost split triangle $\xymatrixcolsep{20pt}\xymatrix{\! M^\cdt \hspace{-1pt} \ar[r] & L^\cdt \hspace{-1pt}  \ar[r] & N^\cdt \hspace{-1pt}  \ar[r] & M^\cdt[1]\!}$ \vspace{-1pt} if and only if $M^\cdt$ has a finite graded injective coresolution over ${\rm ginj}\La$. 

\end{enumerate} \end{Prop}

\noindent{\it Proof.} Since $\La$ is locally bounded, ${\rm gproj}\La$ and ${\rm ginj}\La$ are contained in ${\rm gmod}^{\hspace{.4pt}b\hspace{-3pt}}\La$. Thus, ${\rm gmod}^{\hspace{.4pt}b\hspace{-3pt}}\La$ is a Hom-finite Krull-Schmidt abelian $k$-category with enough projective objects and enough injective objects. In view of Theorem \ref{N-Functor}, we have a Nakayama functor $\nu: {\rm gproj}\La\to {\rm gmod}^{\hspace{.4pt}b\hspace{-3pt}}\La.$ Now, the result follows directly from Theorem 5.12 in \cite{LiN}. The proof of the proposition is completed. 

\vspace{2pt}

We are ready to obtain the sufficient and necessary conditions for the bounded derived category of finite dimensional graded modules to have almost split triangles.

\begin{Theo}\label{GAST-lb-gdim}

Let $\La=kQ/R$ be a locally bounded graded algebra, where $Q$ is a locally finite quiver. Then
$D^{\hspace{.5pt}b}({\rm gmod}^{\hspace{.5pt}b\hspace{-2.5pt}}\La)$ has almost split triangles on the right $($respectively, left$)$ if and only if every graded simple module in $\mmod \La$ is of finite graded projective $($respectively, injective$)$ dimension. 


\vspace{-2pt}

\end{Theo}

\noindent{\it Proof.} Since $\La$ is locally bounded,  
${\rm gmod}^{\hspace{.4pt}b\hspace{-3pt}}\La$ is a Hom-finite Krull-Schmidt abelian $k$-category with enough projective objects and enough injective objects. \vspace{.5pt} So, every complex in $C^b({\rm gmod}^{\hspace{.4pt}b\hspace{-3pt}}\La)$ has a graded projective resolution in $C^-({\rm gmod}^{\hspace{.4pt}b\hspace{-3pt}}\La)$ and a graded injective coresolution in $C^+({\rm ginj}\La)$; see \cite[(7.5)]{Gri}. Then, it is not hard to see that every complex in $C^b({\rm gmod}^{\hspace{.5pt}b\hspace{-2.5pt}}\La)$ has a finite graded projective resolution over ${\rm gproj}\La$ if and only if every module in ${\rm gmod}^{\hspace{.5pt}b\hspace{-2.5pt}}\La$ is of finite graded projective dimension, or equivalently, every graded simple module in ${\rm gmod}^{\hspace{.5pt}b\hspace{-2.5pt}}\La$ is of finite graded projective dimension. Dually, every complex in $C^b({\rm gmod}^{\hspace{.5pt}b\hspace{-2.5pt}}\La)$ has a finite graded injective coresolution over ${\rm ginj}\La$ if and only if every simple module in ${\rm gmod}^{\hspace{.5pt}b\hspace{-2.5pt}}\La$ is of finite graded injective dimension. Now, the statement follows immediately from Proposition \ref{GAST-lb}. The proof of the theorem is completed. 

\vspace{2pt}

\noindent{\sc Remark.} The bounded derived category of finite dimensional modules over a finite dimensional ungraded algebra has almost split triangles on either side if and only if it has almost split triangles on both sides; see  \cite[(1.5)]{Ha3}. As shown below, this is not the case in the locally bounded graded setting.

\vspace{2pt}

\noindent{\sc Example.} Consider the locally bounded graded algebra $\La=kQ/R$, where \vspace{-1pt}
$$\xymatrix{Q: & \cdots \ar[r] & n\ar[r]^-{\alpha_n} & n-1 \ar[r] & \cdots\ar[r] &2 \ar[r]^{\alpha_2} & 1 \ar[r]^{\alpha_1} & 0}\vspace{-2pt}$$
and $R=k\langle\alpha_{i}\alpha_{i+1}\mid i\ge 1\rangle$. Given $n\ge 0$, we see that $S_n$ is of graded projective dimension $n$ and of infinite graded injective dimension. By Theorem \ref{GAST-lb-gdim}, $D^{\hspace{.5pt}b}({\rm gmod}^{\hspace{.5pt}b\hspace{-2.5pt}}\La)$ has almost split triangles on the right but not on the left.

\section{\sc Graded representations of quivers}

In this section, we shall specialize to graded representations of an arbitrary locally finite quiver $Q$. In this case, we shall strengthen some of the results obtained in Sections 3 and 4, which are analogous to those for ungraded representations of strongly locally finite quivers stated in \cite{BLP}. 

\subsection{\sc Graded almost split sequences} We shall restrict our attention to finitely presented graded representations and finitely co-presented graded representations. 

\begin{Lemma}\label{hered_fpres}

Let $Q$ be a locally finite quiver. Then ${\rm gmod}^{\hspace{.5pt}+\hspace{-.5pt},\hspace{.6pt}p\hspace{-.5pt}}kQ$ and ${\rm gmod}^{\hspace{.5pt}-\hspace{-.5pt},\hspace{.6pt}i\hspace{-.5pt}}kQ$ are hereditary abelian $k$-categories, which are Hom-finite and Krull-Schmidt.
    
\end{Lemma}

\noindent{\it Proof.} By Lemma \ref{fpres}, ${\rm gmod}^{+,\hspace{.6pt}p}kQ$ and ${\rm gmod}^{-,\hspace{.6pt}i}kQ$ are Hom-finite Krull-Schmidt and extension-closed in $\GMod kQ$. Note that the category of all unitary left $kQ$-modules is here\-ditary; see \cite[(8.2)]{GaR}. In particular, ${\rm GExt}^2_{kQ}(M, N)=0$ for all $M, N\in {\rm GMod}\hspace{.5pt}kQ$. So, ${\rm GMod}\hspace{.5pt}kQ$ is hereditary. Since ${\rm GMod}\hspace{.5pt}kQ$ has enough projective objects; see (\ref{GM-pi}), a subobject of a projective object in ${\rm GMod}\hspace{.5pt}kQ$ is projective. 

Consider a morphism $f: P\to P'$ in ${\rm gproj}\hspace{.6pt}kQ$. Since ${\rm Im}(f)$ is graded projective, $P\cong {\rm Ker}(f)\oplus {\rm Im}(f)$. So, ${\rm Ker}(f)\in {\rm gproj}\hspace{.6pt}kQ$. By Proposition 2.1 in \cite{Aus2}, ${\rm gmod}^{+,\hspace{.6pt}p\hspace{-.5pt}}kQ$ is closed under kernels and cokernels, and consequently, it is abelian. Since ${\rm GMod}\hspace{.6pt}kQ$ is hereditary, so is ${\rm gmod}^{+,\hspace{.6pt}p\hspace{-1pt}}kQ$. Finally, by Lemma \ref{fpres}(1), ${\rm gmod}^{-,\hspace{.6pt}i\hspace{-.5pt}}kQ$ is also hereditary and abelian. The proof of the lemma is completed.


\vspace{2pt}

Note that $kQ$ is locally left bounded if and only if $Q$ has no infinite path with a starting point, and it is locally right bounded if and only if $Q$ has no infinite path with an end point. 
 
\begin{Theo}\label{ass_rep}

Let $Q$ be a locally finite quiver. The following statement holds.

\begin{enumerate}[$(1)$]

\vspace{-2pt}

\item The abelian category ${\rm gmod}^{+,\hspace{.5pt}p\hspace{-.5pt}}kQ$ has almost split sequences on the left if and only if $Q$ contains no infinite path with a starting point.

\item The abelian category ${\rm gmod}^{-,\hspace{.5pt}i\hspace{-.5pt}}kQ$ has almost split sequences on the right if and only if $Q$ contains no infinite path with an end point.

\vspace{1pt}

\item Both ${\rm gmod}^{+,\hspace{.5pt}p}kQ$ and ${\rm gmod}^{-,\hspace{.5pt}i\hspace{-.5pt}}kQ$ have almost split sequences if and only if $Q$ contains no infinite path. 

\end{enumerate} \end{Theo}

\noindent{\it Proof.} The sufficiency of Statement (1) follows from Theorem \ref{ass_pi}.
Suppose that $Q$ has an infinite path with a starting point, say it starts with an arrow $a\to b$. Then, $P_b$ is infinite dimensional and not graded injective. By Theorem \ref{fd-ass}(1), ${\rm gmod}^{+,\hspace{.5pt}p\hspace{-.5pt}}kQ$ has no almost split sequence starting with $P_b$. 
This establishes Statement (1). And Statement (2) follows dually. Finally, the necessity of Statement (3) follows from Statements (1) and (2). If $Q$ contains no infinite path, then ${\rm gmod}^{+,\hspace{.5pt}p}kQ={\rm gmod}^{-,\hspace{.5pt}i\hspace{-.5pt}}kQ={\rm gmod}^bkQ$, which has almost split sequences by Theorem \ref{rlb-ass}. The proof of the theorem is completed. 

\vspace{-15pt}

\subsection{\sc Graded almost split triangles} First of all, by Lemma \ref{hered_fpres} and Proposition \ref{Der_add}, both $D^{\hspace{.5pt}b\hspace{-0.6pt}}({\rm gmod}^{+\m,\hspace{.4pt}p\m}kQ)$ and $D^{\hspace{.5pt}b\hspace{-0.6pt}}({\rm gmod}^{-\m,\hspace{.4pt}i\m}kQ)$ are Hom-finite and Krull-Schmidt full triangulated subcategories of $D^{\hspace{.5pt}b\hspace{-0.8pt}}({\rm gmod} kQ)$. 


\vspace{2pt}

\begin{Theo}\label{GAST-1-hered}

Let $Q$ be a locally finite quiver. The following statements hold.

\begin{enumerate}[$(1)$]

\vspace{-5.5pt}

\item If \vspace{-1pt} $M^\cdt \!\in \! D^{\hspace{.5pt}b\hspace{-0.6pt}}({\rm gmod}^{+\m,\hspace{.4pt}p\hspace{-.5pt}}kQ)$ is indecomposable, then $D^{\hspace{.5pt}b\hspace{-0.6pt}}({\rm gmod}\hspace{.5pt} kQ)$ has an almost split triangle $\xymatrixcolsep{22pt}\xymatrix{\hspace{-2pt} N^\cdt \!\m \ar[r] & L^\cdt \! \ar[r] & M^\pdt \! \ar[r] & N^\pdt[1],\hspace{-2pt}}\vspace{-3pt}$ where $N^\cdt\in D^{\hspace{.5pt}b\hspace{-0.6pt}}({\rm gmod}^{\m-\m,\hspace{.4pt}i\hspace{-.5pt}}kQ)$.

\vspace{2pt}

\item If $M^\cdt \! \m \in \!\m D^{\hspace{.5pt}b\hspace{-0.6pt}}({\rm gmod}^{-\m, \hspace{.4pt}i}kQ)$ is indecomposable, \vspace{-1pt} then $D^{\hspace{.5pt}b\hspace{-0.6pt}}({\rm gmod}\hspace{.5pt}kQ)$ has an almost split triangle \hspace{-3pt}$\xymatrixcolsep{20pt}\xymatrix{M^\cdt \!\m \ar[r] & L^\cdt \! \ar[r] & N^\pdt \! \ar[r] & M^\pdt[1],}$ \vspace{-1pt} \hspace{-5pt} where $N^\cdt\in D^{\hspace{.5pt}b\hspace{-0.6pt}}({\rm gmod}^{+\m,\hspace{.4pt}p\hspace{-.5pt}}kQ)$. 

\end{enumerate} \end{Theo}

 \noindent{\it Proof.} We shall only prove Statement (1). Let $M^\cdt \in D^{\hspace{.5pt}b\hspace{-0.8pt}}({\rm gmod}^{+\m,\hspace{.4pt}p\hspace{-.5pt}}kQ)$ be indecomposable. Since ${\rm gmod}^{+\m,\hspace{.4pt}p\hspace{-.5pt}}kQ$ is hereditary; see (\ref{hered_fpres}), $M^\cdt$ has a graded projective resolution $P^\pdt\in C^{\hspace{.5pt}b\hspace{-0.6pt}}({\rm gproj}\La)$. Since $M^\cdt \cong P^\pdt$ in $D^{\hspace{.5pt}b\hspace{-0.6pt}}({\rm gmod}\La)$, by Theorem \ref{GAST-1}, there exists an almost split triangle $\xymatrixcolsep{22pt}\xymatrix{\hspace{-3pt} \nu P^\pdt[-\m1] \ar[r] & L^\cdt \ar[r] & M^\cdt\ar[r] & \nu P^\pdt \hspace{-2pt}}\vspace{-3pt}$ in $D^{\hspace{.5pt}b\hspace{-0.6pt}}({\rm gmod}\La)$, where $ \nu P^\pdt[-\m1] \! \in \! D^{\hspace{.5pt}b\hspace{-0.6pt}}({\rm gmod}^{-\m,\hspace{.4pt}i\hspace{-.5pt}}kQ)$. The proof of the theorem is completed.

\vspace{2pt}

The following statement is analogous to Theorem 7.11 in \cite{BLP}, which contains a partial converse of 

\begin{Theo}\label{GAST-3_hered}

Let $Q$ be a locally finite quiver. Then, the following statements are equivalent$\,:$

\begin{enumerate}[$(1)$]

\vspace{-2.5pt}

\item $D^b({\rm gmod}^{+,\hspace{.5pt}p\hspace{-.8pt}}kQ)$ has almost split triangles $(\hspace{-1pt}$on the left, on the right$\hspace{.8pt});$ 

\item $D^b({\rm gmod}^{-,\hspace{.5pt}i\hspace{-1pt}}kQ)$ has almost split triangles $(\hspace{.4pt}$on the left, on the right$);$

\item $Q$ has no infinite path $(\hspace{-1pt}$with a starting point, with an end point$\hspace{.8pt}).$ 
\end{enumerate}


\end{Theo}

\noindent{\it Proof.} In view of the duality $\mk D: {\rm gmod}^{+,\hspace{.5pt}p\hspace{-1pt}}kQ^{\rm o}\to {\rm gmod}^{-,\hspace{.5pt}i\hspace{-1pt}}kQ$, we shall only prove the equivalence of Statements (1) and (3). 
Suppose first that $D^b({\rm gmod}^{+,\hspace{.5pt}p\hspace{-.5pt}} kQ)$ has almost split triangles on the left. By Corollary 7.3(1) 
in \cite{BLP}, ${\rm gmod}^{+,\hspace{.5pt}p\hspace{-.5pt}} kQ$ has almost split sequences on the left; and by Theorem \ref{ass_rep}(1), $Q$ has no infinite path with a starting point. Suppose conversely that $Q$ has no infinite path with a starting point. Then, ${\rm gmod}^{+,\hspace{.5pt}p\hspace{-1pt}}kQ={\rm gmod}^{\hspace{.5pt}b \hspace{-.5pt}}kQ$; and by Theorem \ref{ass_rep}(1), ${\rm gmod}^{+,\hspace{.5pt}p\hspace{-1pt}} kQ$ has almost split sequences on the left. Moreover, by Proposition \ref{fdi_pres}, every indecomposable injective object in ${\rm gmod}^{+,\hspace{.5pt}p\hspace{-1pt}}kQ$ is isomorphic to a finite dimensional module $I_x$ in ${\rm ginj}\hspace{.5pt} kQ$ for some $x\in Q_0$, whose socle $S_x$ is simple and admits a projective cover $P_x$ in ${\rm gmod}^{+,\hspace{.5pt}p\hspace{-1pt}}kQ$. By Corollaries 2.2(2) and 7.3(1) in \cite{BLP}, $D^b({\rm gmod}^{+,\hspace{.5pt}p\hspace{-.5pt}} kQ)$ has almost split triangles on the left. 

Next, suppose that $D^b({\rm gmod}^{+,\hspace{.5pt}p\hspace{-.5pt}}kQ)$ has almost split triangles on the right. By Corollary 7.3(1) in \cite{BLP}, every $S_x$ with $x\in Q_0$ admits an injective envelope in ${\rm gmod}^{+,\hspace{.5pt}p\hspace{-.5pt}}kQ$. By Proposition \ref{lobo-abel}(1), 
$Q$ has no infinite path with an end point. Suppose conversely that $Q$ has no infinite path with an end point. 
By Proposition \ref{lobo-abel}, 
every $S_x$ with $x\in Q_0$ has an injective envelope $I_x$ in ${\rm gmod}^{+,\hspace{.5pt}p\hspace{-.5pt}}kQ$, and by Theorem \ref{ass_pi}, ${\rm gmod}^{+,\hspace{.5pt}p\hspace{-.5pt}}kQ$ has almost split sequence on the right. By Corollary 7.3(2) in \cite{BLP}, $D^b({\rm gmod}^{+,\hspace{.5pt}p\hspace{-.5pt}} kQ)$ has almost split triangles on the right. Finally, combining what has been shown, we see that $D^b({\rm gmod}^{+,\hspace{.5pt}p\hspace{-.5pt}} kQ)$ has almost split triangles if and only if $Q$ has no infinite path. The proof of the theorem is completed.




\begin{thebibliography}{99}

\smallskip

\bibitem{FWAF} {\sc F. W. Anderson and K. R. Fuller}, ``Rings and categories of module,'' Second Edition, Graduate Texts in Mathematics 13 (Springer-Verlag, New York, 1992).



\bibitem{Aus2} {\sc M. Auslander}, ``Coherent functors,'' Proceedings of the Conference on Categorical Algebra (La Jolla, California, 1965) (Springer, New York, 1966) 189-231.

\bibitem{Aus} {\sc M. Auslander}, ``Functors and morphisms determined by objects,'' Lecture Notes in Pure Appl. Math. 37 (Dekker, New York, 1978) 
1-244.



\bibitem{AuR2} {\sc M. Auslander and I. Reiten,} ``Representation theory of artin algebras IV,'' Comm. Algebra 5 (1977) 443-518.

\bibitem{AuR0} {\sc M. Auslander and I. Reiten,} ``Almost split sequences in dimension 2," Adv. Math. 66 (1987) 88-118.

\bibitem{AuR} {\sc M. Auslander and I. Reiten,} ``Almost split sequences for abelian group graded rings,'' J. Algebra 114 (1988) 29-39.


\bibitem{ARS} {\sc M. Auslander, I. Reiten and Smal\o}, ``Representation theory of artin algebras,'' Cambridge Studies in Advanced Mathematics 36 (Cambridge University Press, Cambridge, 1995).



%

\bibitem{BaL} {\sc R. Bautista and S. Liu}, ``Covering theory for linear categories with application to derived categories,'' J. Algebra 406 (2014) 173-225.


\bibitem{BLP} {\sc R. Bautista, S. Liu and C. Paquette}, ``Representation theory of strongly locally finite quivers,'' Proc. London Math. Soc. 106 (2013) 97-162.

%



\bibitem{BGS} {\sc A. Beilinson, V. Ginzburg and W. Soergel}, ``Koszul duality patterns in representation theory,'' J. Amer. Math. Soc. 9 (1996) 473-527.

%

%

\bibitem{BoG} {\sc K. Bongartz and P. Gabriel}, ``Covering spaces in representation theory,'' Invent. Math. 65 (1982) 331-378.






\bibitem{CEi} {\sc H. Cartan and S. Eilenberg}, ``Homological Algebra," (Princeton University Press, Princeton, New Jersey, 1956).

\bibitem{EKo} {\sc A. Elduque and M. Kochetov}, ``Gradings on simple Lie algebras,"
Mathematical Surveys and Monographs 189 (American Mathematical Society, 2013).

%




%
%

%




\bibitem{GaR} {\sc P. Gabriel and A. V. Roiter}, ``\hspace{1pt}Representations of Finite-Dimensional Algebras," Algebra VIII, Encyclopedia Mathematical Sciences 73 (Springer, Berlin, 1992).

%


\bibitem{GoG} {\sc R. Gordon and E. L. Green}, ``Representation theory of graded artin algebras," J. Algebra 76 (1982) 138-152.

\bibitem{Gre} {\sc E. L. Green}, ``Graphs with relations, coverings and group-graded algebras,'' Trans. Amer. Math. Soc. 279 (1983) 297-310.

\bibitem{GrS} {\sc E. L. Green and S. Schroll}, ``Multiserial and special multiserial algebras and their representations," Adv. Math. 302 (2016) 1111-1136.

%








\bibitem{Gri} {\sc P.-P. Grivel},``Cat\'egories d\'eriv\'ees et foncteurs d\'eriv\'es,'' in {\it Algebraic D-modules}, Perspective in Mathematics 2 (Academic Press Inc., Boston, 1987) 1-108.


\bibitem{Har} {R. Hartshorne}, ``Algebraic Geometry,'' Graduate Texts in Mathematics 52 (Springer-Verlag, New York, Heidelberg, 1977).


\bibitem{Ha1} {\sc D. Happel}, ``On the derived category of a finite-dimensional algebra,'' Comment. Math. Helv. 62 (1987) 339-389.


\bibitem{Ha3} {\sc D. Happel}, ``Auslander-Reiten triangles in derived categories of finite dimensional algebras,'' Proc. Amer. Math. Soc. 112 (1991) 641-648.

%

\bibitem{VHW} {\sc H.-J. Von H\"{o}hne and J. Waschb\"{u}sch}, ``Die struktur
n-reihiger Algebren," Comm. Algebra 12 (1984) 1187-1206.







\bibitem{HKKS} {\sc H. Krause,} ``Krull-Remak-Schmidt categories and projective covers,'' Expo. Math. 33 (2015)
535-549.

\bibitem{HKra} {\sc H. Krause,} ``A short proof for Auslander's defect formula,'' Linear Algebra Appl. 365 (2003) 267-270.

\bibitem{HKMS} {\sc H. Krause and M. Saor\'in,} ``On minimal approximations of modules,'' 
Contemp. Math. 229 (Amer. Math. Soc., Providence, RI, 1998) 227-236.

\bibitem{LeC} {\sc J. Le and X. Chen}, ``Karoubianness of a triangulated category,'' J. Algebra 310 (2007) 452-457.

\bibitem{Len} {\sc H. Lenzing}, ``{\hskip 1pt}Hereditary Categories", {\it Handbook of tilting theory}, London Mathematical Society Lectures Note Series 332 (Cambridge University Press, Cambridge, 2007) 105-146.



\bibitem{SLiu} {\sc S. Liu}, ``Auslander-Reiten theory in a Krull-Schmidt category,'' Sao Paulo J. Math. Sci. 4 (2010) 425-472. 

%

\bibitem{LiM} {\sc S. Liu and J.-P. Morin}, ``The strong no loop conjecture for special biserial algebras,'' Proc. Amer. Math. Soc. 132 (2004) 3513-3523.

\bibitem{LiN} {\sc S. Liu and H. Niu}, ``Almost split sequences in tri-exact categories,'' J. Pure Appl. Algebra 226 (2022) 1-31.

\bibitem{SLPC} {\sc S. Liu, P. Ng and C. Paquette,} ``Almost split sequence and approximation,'' Algebr. Represent. Ther. 16 (2013) 1809-1827.


%

\bibitem{Mac} {\sc S. MacLane}, ``Homology,'' (Springer-Verlag, Berlin, Heidelberg, New York, 1971).

\bibitem{Mat} {\sc H. Matsumura}, ``Commutative Ring Theory," \hspace{-2pt} Cambridge University Press, Cambridge, 1989.





\bibitem{RMV2} {\sc R. Martinez-Villa,} \hspace{-5pt} ``Graded self-injective and Koszul algebras," \hspace{-3pt} J. Algebra 215 (1999) 34-72.

\bibitem{Mar} {\sc R. Martinez-Villa,} ``Introduction to Koszul algebras," Rev. Un. Mat. Argentina 48 (2007) 67-95.





\bibitem{Mil} {\sc D. Milicic}, ``Lectures on derived categories," www.math.utah.edu/\~{}milicic/Eprints/dercat.pdf




\bibitem{NCF2} {\sc C. Nastasescu and F. van Oystaeyen,} ``Methods of Graded Rings,'' Lecture Notes in Mathematics 1836 (Springer-Verlag, Berlin, Heidelberg, 2004).






%

%

%

%


%

%





\bibitem{Wei} {\sc C. A. Weibel}, ``An Introduction to Homological Algebra", Cambridge Studies in Advanced Mathematics 38 (Cambridge University Press, Cambridge, 1994).

\end{thebibliography}
\end{document}